# Times Square sampling: an adaptive algorithm for free energy estimation

Cristian Predescu,[1,†] Michael Snarski,[1] Avi Robinson-Mosher,[1] Duluxan Sritharan,[1]
Tamas Szalay,[1] and David E. Shaw[1,2,†]

[1] D. E. Shaw Research, New York, NY 10036, USA.

[2] Department of Biochemistry and Molecular Biophysics, Columbia University,
New York, NY 10032, USA.

† To whom correspondence should be addressed.

Cristian Predescu

E-mail: Cristian.Predescu@DEShawResearch.com

Phone: (212) 478-0433

Fax: (212) 845-1433

David E. Shaw

E-mail: David.Shaw@DEShawResearch.com

Phone: (212) 478-0260

Fax: (212) 845-1286

# Abstract

Estimating free energy differences, an important problem in computational drug discovery and in a wide range of other application areas, commonly involves a computationally intensive process of sampling a family of high-dimensional probability distributions and a procedure for computing estimates based on those samples. The variance of the free energy estimate of interest typically depends strongly on how the total computational resources available for sampling are divided among the distributions, but determining an efficient allocation is difficult without sampling the distributions. Here we introduce the *Times Square sampling* algorithm, a novel on-the-fly estimation method that dynamically allocates resources in such a way as to significantly accelerate the estimation of free energies and other observables, while providing rigorous convergence guarantees for the estimators. We also show that it is possible, surprisingly, for on-the-fly free energy estimation to achieve lower asymptotic variance than the maximum-likelihood estimator MBAR, raising the prospect that on-the-fly estimation could reduce variance in a variety of other statistical applications.

# Introduction

An outstanding problem in computational statistics is how to efficiently sample multi-modal, high-dimensional probability distributions in order to estimate ratios of their normalizing constants. One important instance of this problem is the calculation of free energy differences, which has applications in various fields, including computational drug discovery. The probability distributions to be sampled in order to estimate free energy differences are Gibbs distributions, which give the probability that a physical system adopts a certain state; the free energy difference between a pair of Gibbs distributions is proportional to the logarithm of the ratio of their normalizing constants. Sampling Gibbs distributions is typically computationally intensive, and often intractable with existing methods.

Sampling algorithms such as Markov chain Monte Carlo (MCMC) and molecular dynamics (MD) explore the state space of a simulated system through local moves and can get trapped in regions of high probability density (which are also called "modes"). To facilitate sampling of the many modes, techniques such as replica exchange[1–3] and simulated tempering[4,5] are commonly used. These methods "enhance" the sampling by introducing a parameterized family of overlapping distributions, and the problem becomes one of sampling each distribution within the family, which is typically still a formidable computational task. If enough uncorrelated samples can be collected, free energy differences within the family can be estimated using a technique such as the weighted histogram analysis method (WHAM),[6] or its "binless" version, the multistate Bennett acceptance ratio (MBAR).[5,7–10]

A major drawback of enhanced sampling methods such as replica exchange and simulated tempering is that their sampling efficiency is strongly dependent on the underlying free energy landscape, which is, by the nature of the estimation problem, typically unknown a priori. Each distribution to be sampled does not contribute equally to the variance (i.e., the mean-squared error) of the free energy difference estimate of interest, and appropriately allocating limited sampling resources among these distributions requires that their relative contributions be estimated. These estimates are, however, themselves determined by underlying samples, creating a cycle that makes it difficult to improve the efficiency of the free energy estimation procedure.

To improve the efficiency of free energy estimation, methods that blend sampling and estimation have been developed in recent years, in which on-the-fly estimates guide sampling. Using the method of stochastic approximation, Z. Tan has shown that one can use simulated tempering to estimate free energy differences on the fly, and moreover do so with the smallest possible asymptotic variance (i.e., in the limit of infinitely many samples) when the allocation of resources is fixed.[11] A related algorithm, accelerated weighted histogram (AWH), estimates on the fly not only free energy differences, but also ensemble averages for other functions defined on the state space (i.e., observables). This allows the sampler to estimate an asymptotic allocation of resources that lowers the variance of the desired free energy difference estimate.[12,13] These advances underline the value of incorporating estimate information into the sampling procedure in order to reduce the computational effort spent sampling less relevant parts of the free energy landscape. Without establishing explicit conditions that guarantee convergence of free energy estimators, however, it is difficult to develop resource allocation schemes that can further reduce variance or accelerate convergence.

Here we present Times Square sampling (TSS), an algorithm that provides a straightforward recipe to estimate ensemble averages on the fly for any function defined on the state space, and enables the use of this information to allocate computational resources to an intermediate distribution within the parameterized family where more sampling is needed. TSS gives simple conditions that establish how resources can be allocated based on the state of the system, the parameter of the intermediate distribution being sampled, or the current free energy estimates, while guaranteeing convergence. This enabled us to introduce, for example, a *visit-control mechanism*, which modulates the likelihood of visiting upcoming distributions based on whether they have been under- or over-visited, that in our experiments greatly sped up convergence of the free energy estimates.

Importantly, we also prove that an additional property of TSS that we observed, "self-adjustment" (a term coined by Z. Tan for a property by which the likelihood of sampling upcoming distributions is modulated based on the current free energy estimates[11]), makes it possible for on-the-fly estimation algorithms to achieve lower asymptotic variance than the maximum-likelihood estimator MBAR. Self-adjustment enables this improvement—a surprising one, given that maximum-likelihood estimators have the lowest variance when independent samples are available[14]—by taking advantage of on-the-fly information, which is unavailable to MBAR. This raises the intriguing prospect that it may be possible to develop on-the-fly estimation procedures that improve over maximum-likelihood estimation in a variety of other statistical applications.

We also describe features of TSS relevant to its practical implementation. Performing energy evaluations for all intermediate distributions can be prohibitively expensive. We thus include a windowing system, which limits the number of energy evaluations required to update the estimates. Another issue is that convergence can be slow if the initial free energy estimates are sufficiently inaccurate. We thus introduce a history-forgetting mechanism to TSS, which in many cases significantly reduces the number of steps it takes before the statistical variance becomes the dominant contributor to the error. We summarize our implementation of TSS in Section 2 (and provide a full implementation in the Appendix), and in Section 3 we demonstrate its use by applying it to a set of biomolecular free energy calculations using MD as the sampler.

## 1. Mathematical background and theoretical derivation of TSS

This section introduces mathematical notation used throughout the paper and explores theoretical aspects of the TSS algorithm. In Section 1.1 we describe the simulated tempering algorithm, the sampling procedure on which TSS is based; Section 1.2 introduces stochastic approximation, which is used in TSS for on-the-fly estimation; Section 1.3 brings together the sampling and estimation, defines the TSS algorithm, and states theorems about its convergence; finally, Section 1.4 provides an example amenable to explicit calculations to study the effects of the mechanisms introduced in Section 1.3 on the free energy estimators.

Our starting point is a compact state space $\mathcal{S} \subset \mathbb{R}^N$, elements of which are (for instance) positions of atoms; a compact parameter space $\Lambda \subset \mathbb{R}^d$, which represents a quantity to be varied (such as temperature); and a family of probability distributions $\{\rho_\lambda(x)dx, \lambda \in \Lambda\}$ defined on $\mathcal{S}$ and parameterized by $\lambda \in \Lambda \subset \mathbb{R}^d$. The distributions are defined in terms of dimensionless potential functions $H_\lambda: \mathbb{R}^N \to \mathbb{R}$ (which we will call *Hamiltonians*). These functions are assumed to be continuously differentiable in $(x, \lambda)$, and they specify the density $\rho_\lambda(x)$ with respect to the Lebesgue measure $dx$ by

$$\rho_\lambda(x) = \frac{e^{-H_\lambda(x)}}{Z_\lambda^\star}, \qquad x \in \mathcal{S}, \lambda \in \Lambda, \tag{1}$$

where $Z_\lambda^\star = \int_{\mathcal{S}} e^{-H_\lambda(x)} dx$ are unknown normalizing constants called "partition functions" in the context of chemical physics. It is desired to estimate ratios of partition functions $Z_\lambda^\star / Z_{\lambda'}^\star$ or, equivalently, "free energy differences" $F_\lambda^\star - F_{\lambda'}^\star$ for $\lambda \neq \lambda'$, where we will call $F_\lambda^\star = -\log Z_\lambda^\star$ the (dimensionless) "free energy". (It is important to note that choosing a good parameterized family $\{\rho_\lambda, \lambda \in \Lambda\}$ is a very difficult task in itself; the problem we address in this manuscript is that of optimizing the sampling and estimation for a given family.)

Enhanced sampling methods such as replica exchange and simulated tempering assume $\Lambda \subset \mathbb{R}^d$, select a discretization $\{\lambda_1, \ldots, \lambda_K\} \subset \Lambda$, and sample $\{\rho_{\lambda_k}, k \in [K]\}$, where $[K] = \{1, \ldots, K\}$. The index $k$ is called the *rung*. For simplicity, going forward we write $\rho_k$ for $\rho_{\lambda_k}$ and $H_k$ for $H_{\lambda_k}$ for each $k \in [K]$. In simulated tempering, the method on which TSS is built, one specifies a "rung density" $\gamma = (\gamma_1, \ldots, \gamma_K)$ associated with the discretization $\{\lambda_1, \ldots, \lambda_K\}$, and then samples the joint distribution $\gamma_k \rho_k(x) \, dx \, \mathcal{C}(dk)$, where $\mathcal{C}(dk)$ is the counting measure on the set $[K]$ (i.e., the measure that assigns a unit mass to every singleton set $\{k\} \subset [K]$). To avoid degeneracies, we fix some $\varepsilon > 0$ and define $\mathcal{P}^\varepsilon(K)$ to be the set of probability densities $\phi$ on $[K]$ satisfying $\phi_k \geq \varepsilon$ for each $k \in [K]$, and we only consider rung densities $\gamma \in \mathcal{P}^\varepsilon(K)$. We will see that $\gamma$ can be used to perform resource allocation by guiding the sampling in the rung space.

Throughout, we assume that both $\{\rho_\lambda, \lambda \in \Lambda\}$ and a discretization are given, with $\Lambda \subset \mathbb{R}^d$ (the size of $d$ is irrelevant to the mathematical analysis, but in practice is limited by computational resources to a small integer: less than 5, for direct discretization). To avoid multi-index notation we simply consider an enumeration $\{\lambda_1, \ldots, \lambda_K\} \subset \Lambda$ of states, but the reader should note that the indexing scheme ignores the original topology, and the latter can be important at times. We expect, for example, that whenever $\lambda_k$ and $\lambda_j$ are sufficiently close, the distributions $\rho_k$ and $\rho_j$ exhibit good statistical overlap, as measured by the Rényi entropy divergence of order ½. The list of states $\{\lambda_1, \ldots, \lambda_K\}$ is obtained by discretization, meaning that $K$ grows, e.g., as $100^d$ with $d$. In practice, the exponential growth constrains $d$ to modest values, making the $\Lambda$ space "smaller" than the state space $\mathcal{S}$, which contains $10^4$–$10^7$ degrees of freedom in typical applications.

Consequently, Monte Carlo sampling of the $\Lambda$ space can be done more efficiently than sampling of the state space, a situation we refer to as a "two-timescale sampling" scenario.

## 1.1 *The simulated tempering invariant distribution*

In this section, we describe how simulated tempering is used to sample a family of joint distributions of the form $\gamma_k \, \rho_k(x) \, dx \, \mathcal{C}(dk)$, and we impose assumptions on both the invariant distribution and the Markov transition kernel corresponding to the simulated tempering process. We then explain how the collected samples can be used to estimate the free energy differences, a task which is typically considered separately from the sampling procedure.

The joint invariant density, which is our primary object of interest, is defined for any pair of "hypothesized" parameters $\gamma \in \mathcal{P}^\varepsilon(K)$ and $Z \in \mathbb{R}_{>0}^K$ (or, equivalently, $F \in \mathbb{R}^K$) by

$$p_{\gamma,Z}(x,k) = \frac{\gamma_k Z_k^\star/Z_k}{\sum_{\ell \in [K]} \gamma_\ell \, Z_\ell^\star/Z_\ell} \rho_k(x) = \frac{\gamma_k \, e^{-H_k(x)}/Z_k}{\sum_{\ell \in [K]} \gamma_\ell \, Z_\ell^\star/Z_\ell}. \qquad (2)$$

We will denote an expectation with respect to $p_{\gamma,Z}(x,k) \, dx \, \mathcal{C}(dk)$ by $\mathbb{E}_{\gamma,Z}$, where we recall that $\mathcal{C}$ is the counting measure on $[K]$. The marginal densities of $p_{\gamma,Z}(x,k)$ are

$$p_{\gamma,Z}^{[K]}(k) = \frac{\gamma_k Z_k^\star/Z_k}{\sum_{\ell \in [K]} \gamma_\ell \, Z_\ell^\star/Z_\ell}, \qquad p_{\gamma,Z}^S(x) = \sum_{k \in [K]} \frac{\gamma_k Z_k^\star/Z_k}{\sum_{\ell \in [K]} \gamma_\ell \, Z_\ell^\star/Z_\ell} \rho_k(x). \qquad (3)$$

Although $Z^\star$ (or $F^\star$) is unknown, it is possible to sample the density (2) for any choice of $\gamma$ and $Z$ using, for instance, the simulated tempering sampler (defined below). (Although we use the simulated tempering sampler in this work for the sake of concreteness, we note that the TSS algorithm is compatible with other samplers; aside from regularity conditions, the only important property of a sampler for it to be used with TSS is that the sampler be ergodic for the distribution with density (2).)

The simulated tempering algorithm samples (2) for any fixed $(\gamma, Z) \in \mathcal{P}^\varepsilon(K) \times \mathbb{R}^K$ by composing the pair of Markov transition kernels of densities

$$P_{\gamma,Z}^k(x', k'|x, k) = \delta(x' - x) \frac{\gamma_{k'} e^{-H_{k'}(x)}/Z_{k'}}{\sum_{\ell \in [K]} \gamma_\ell e^{-H_\ell(x)}/Z_\ell} \tag{4}$$

$$P_{\gamma,Z}^x(x'', k''|x', k') = \delta(k'' - k')\mathcal{T}_{k'}^{(n)}(x''|x'). \tag{5}$$

The above densities are written with respect to the product of Lebesgue and counting measures on $\mathcal{S} \times [K]$; $\delta$ denotes the Dirac density for the continuous variable or the Kronecker symbol for the discrete variable. The conditional distribution $\mathcal{T}_{k'}(x''|x')dx''$ is assumed to be a Feller transition kernel, $\mathcal{T}_{k'}^{(n)}(x''|x')dx''$ the same kernel iterated $n$-steps, and for each fixed $k'$ we assume that the unique Markov process corresponding to the kernel converges in distribution to the unique invariant distribution $\rho_{k'}(x)dx$. In practice, the Markov process is simulated on a computer with a local-move-based sampler such as MCMC or MD, and $\mathcal{T}_{k'}(x''|x')dx''$ represents the sampling distribution of a single time step of this sampler, which we only need assume to be ergodic for $\rho_k(x)dx$ for each $k \in [K]$. As noted above, the realities of chemical physics applications typically lead to two time scales for the combined sampling, with (5) mixing slowly and being computationally expensive, and (4) being achieved comparatively cheaply and rapidly between Markov chain updates of $x$. From now on $P_{\gamma,Z} = P_{\gamma,Z}^k \circ P_{\gamma,Z}^x$ will denote the composition of the two transition kernels.

For a fixed pair $(\gamma, Z)$, the empirical distribution of the rung trajectory $\mathcal{K}^t$ using the sampler $P_{\gamma,Z}$ would converge to the marginal distribution $p_{\gamma,Z}^{[K]}(k)\mathcal{C}(dk)$. In other words, under certain ergodicity assumptions,

$$\lim_{t\to\infty} \frac{1}{t} \sum_{s=1}^{t} 1_{\{k\}}(\mathcal{K}^s) = \frac{\gamma_k Z_k^\star/Z_k}{\sum_{\ell\in[K]} \gamma_\ell Z_\ell^\star/Z_\ell} = \frac{\gamma_k e^{F_k - F_k^\star}}{\sum_{\ell\in[K]} \gamma_\ell e^{F_\ell - F_\ell^\star}}, \qquad (6)$$

with $1_{\{k\}}(k')$ denoting the indicator function of the set $\{k\}$. Because the hypothesized free energies $F$ are known, in principle (6) allows for the determination of the exact free energy differences, by solving for $F_k^\star - F_1^\star, k = 2, \ldots, K$. Unless the hypothesized and exact free energies are close to each other, however, the visitation of some of the rungs is a rare event, which leads to very high error in the free energy difference estimates. The free energies can also be estimated using the variable $\mathcal{X}^t$ by means of importance sampling, where $\mathcal{X}^t$ is drawn from the mixture density $p_{\gamma,Z}^S(x) \propto \sum_{\ell\in[K]} \gamma_\ell e^{-H_\ell(\mathcal{X}^s)}/Z_\ell$ defined in (3):

$$\lim_{t\to\infty} \frac{1}{t} \sum_{s=1}^{t} \frac{\gamma_k e^{-H_k(\mathcal{X}^s)}/Z_k}{\sum_{\ell\in[K]} \gamma_\ell e^{-H_\ell(\mathcal{X}^s)}/Z_\ell} = \frac{\gamma_k Z_k^\star/Z_k}{\sum_{\ell\in[K]} \gamma_\ell Z_\ell^\star/Z_\ell} = \frac{\gamma_k e^{F_k - F_k^\star}}{\sum_{\ell\in[K]} \gamma_\ell e^{F_\ell - F_\ell^\star}}. \qquad (7)$$

The estimator (7) is expected to have lower error than the counting estimator (6), though it also suffers from the exponential sensitivity to the hypothesized free energies $F$.

Finally, we note that sampling (2) with a hypothesized parameter $Z_k$ and rung density $\gamma \in \mathcal{P}^\varepsilon(K)$ can be thought of as sampling (2) with the exact (unknown) parameter $Z_k^\star$ along with a *tilted* rung density $\tilde{\gamma}_k \propto \gamma_k Z_k^\star/Z_k$ (i.e., with $\gamma_k$ being scaled by $Z_k^\star/Z_k$). This interplay between the hypothesized parameter $Z_k$ and the rung density $\gamma$ allows us to use (6) and (7) together to mutual advantage, resulting in a mechanism we call *visit control*, which we will introduce in Section 1.3. There remains, however, the issue of exponential sensitivity to the hypothesized free energies $F$, which is a problem for both estimators (6) and (7); we will address this problem in Section 1.2.

*1.2 Stochastic approximation and iterative importance sampling*

To remedy the main disadvantage of simulated tempering, which is that a reasonably accurate guess of the free energies is required for good sampling, one can instead estimate the free energies on the fly and adapt them throughout the simulation. Such a development enables one to make use of the simulated tempering transition kernels (4) and (5) even if the initial guess is not reasonably accurate. In this section, we describe stochastic approximation, which is an appropriate mathematical framework to handle such an on-the-fly procedure, and explain how it has been applied to estimating free energy differences using samples generated by simulated tempering. In the next section, we will extend the stochastic approximation procedure to estimate two additional types of quantities, and introduce the TSS algorithm as a triple of stochastic approximation recursions. We note that we originally derived TSS differently, using the idea of iterative importance sampling (IIS); we will also define IIS in this section, as it provides a useful perspective on the algorithm.

Stochastic approximation originates with the seminal paper of Robbins and Monro,[15] and can be described as a stochastic root-finding procedure. In this section, we outline the procedure with Markovian noise and explain its use for the on-the-fly estimation of averages. In recent years, stochastic approximation has received significant attention in the fields of MCMC and MD (see refs. 16–18 and the survey refs. 19 and 20; see also refs. 12 and 13 for a similar algorithm. For more complete references we refer the reader to refs. 21 and 22).

Let $\mathbb{Y} \subset \mathbb{R}^{n_y}$ and $\Theta \subset \mathbb{R}^{n_\theta}$ denote subsets of Euclidean space, and for each $\theta \in \Theta, y \in \mathbb{Y}$, let $P_\theta(y, dy')$ denote a Markov transition kernel with unique stationary distribution $p_\theta(y) dy$, where $dy$ is taken to be the Lebesgue measure for simplicity. Suppose we wish to find a root $\theta^\star$ of $g(\theta) = 0$, where $g: \Theta \to \mathbb{R}^{n_\theta}$ is defined by an expectation with respect to the stationary distribution $p_\theta(y)dy$:

$$g(\theta) = \mathbb{E}_\theta[G(Y;\theta)] = \int_\mathbb{Y} G(y;\theta)p_\theta(y)dy. \tag{8}$$

The function $g(\theta)$ is referred to as the "mean field" of the observable $G(y;\theta)$.

Fix an initial guess $(\mathcal{Y}^0, \theta^0) \in \mathbb{Y} \times \Theta$, and for $t \geq 0$ define the process

$$\mathcal{Y}^{t+1} \sim P_{\theta^t}(\mathcal{Y}^t, \cdot) \tag{9}$$

$$\theta^{t+1} = \theta^t + \frac{1}{t+1} \Gamma\, G(\mathcal{Y}^{t+1}; \theta^t), \tag{10}$$

where for each $t \geq 0$, $\Gamma$ is an $n_\theta \times n_\theta$ matrix called the *gain matrix* and $t^{-1}$ is the *gain factor*. Together, equations (8), (9), and (10) define the stochastic approximation procedure with Markovian noise. In the limit $t \to \infty$, the average behavior of the iterates $\theta^t$ is described by the ordinary differential equation (ODE) expressed in terms of the mean field $g$,

$$\dot{\theta}(t) = \Gamma g(\theta(t)). \tag{11}$$

(Here, the overhead dot notation indicates a time derivative.) This "mean field ODE" can be derived heuristically from (10) by rearranging the left-hand side to form a difference quotient, taking limits $t \to \infty$, and replacing the right-hand side by its average (8). The theory of stochastic approximation makes this heuristic rigorous, by leveraging stability properties of the ODE to establish convergence of the stochastic iterates $\theta^t$ (see Section 4, Part 1 of the SM). One can choose a more general gain factor of $t^{-\beta}$ for $\beta \in (0.5, 1]$, but $\beta = 1$ is needed to obtain the minimal asymptotic variance of the estimate $\theta^t$. With $\beta = 1$, it has further been shown that the asymptotically optimal gain matrix is $\Gamma = -\left(J_g(\theta^\star)\right)^{-1}$, where $J_g(\theta)$ denotes the Jacobian of $g$ evaluated at $\theta$ (see, for example, Section 3.2.3, Proposition 4 of ref. 22). Under certain regularity conditions on the function $G$, the transition kernels $P_\theta$, and the gain $\Gamma$, one can establish the almost sure convergence of $\theta^t$ to $\theta^\star$ (see ref. 23 or 24).

For some fixed $\gamma \in \mathcal{P}^\varepsilon(K)$, the stochastic approximation procedure with Markovian noise can be applied to free energy estimation by taking $\mathcal{Y}^t = (\mathcal{X}^t, \mathcal{K}^t)$, the free energy $F$ as the parameter $\theta$, and the kernel to be the simulated tempering kernel $P_{\gamma,Z}$ for some fixed $\gamma \in \mathcal{P}^\varepsilon(K)$. As we have seen in the previous section, the common right-hand side of (6) and (7) uniquely identifies the free energies (up to an additive constant) if set equal to $\gamma_k$. Those equalities (which determine the free energies) can be identified as the roots of the following mean field

$$g_k(F) \equiv 1 - \frac{e^{F_k - F_k^\star}}{\sum_{\ell \in [K]} \gamma_\ell e^{F_\ell - F_\ell^\star}} = 0, \quad \forall k \in [K]. \tag{12}$$

The left-hand sides of (6) and (7) provide the following two observables $G_k(\cdot; F)$ that by averaging against $p_{\gamma,Z}(x,k)\, dx\, \mathcal{C}(dk)$ produce (12); namely, $g_k(F)$ can be written as $\mathbb{E}_{\gamma,Z}[G_k(\mathcal{K}; F)] = 0$ if

$$G_k(j; F) = 1 - \frac{1}{\gamma_k} 1_{\{k\}}(j) \tag{13}$$

or $\mathbb{E}_{\gamma,Z}[G_k^F(\mathcal{X}; F)] = 0$ if

$$G_k^F(x; F) = 1 - \frac{e^{F_k - H_k(x)}}{\sum_{\ell \in [K]} \gamma_\ell e^{F_\ell - H_\ell(x)}}. \tag{14}$$

Tan applied the Robbins-Monro algorithm to these (and other related) observables and made the crucial observation that the optimal gain matrix is trivial, being given by $\Gamma = I_{K \times K}$ (cf. Theorem 1 in ref. 11), which from (10) gives the recursions

$$F_k^{t+1} = F_k^t + \frac{1}{t+1}\left(1 - \frac{1}{\gamma_k}1_{\{k\}}(\mathcal{K}^t)\right) \tag{15}$$

and

$$F_k^{t+1} = F_k^t + \frac{1}{t+1}\left(1 - \frac{e^{F_k^t - H_k(\mathcal{X}^{t+1})}}{\sum_{\ell \in [K]} \gamma_\ell \, e^{F_\ell^t - H_\ell(\mathcal{X}^{t+1})}}\right), \tag{16}$$

respectively. He noted that (15) has been derived before by Liang et al.[16] as Stochastic Approximation Monte Carlo (SAMC). The recursions (15) and (16) are instances of a class of algorithms Tan has named Self-Adjusted Mixture Sampling (SAMS).

We illustrate how self-adjustment affects convergence of free energy differences, by applying (15) to a simple example. Consider two equally weighted identical distributions defined by $\mathcal{S} = [0,1]$, $H_1(x) = H_2(x) = 0$ and $\gamma_1 = \gamma_2 = 1/2$. We wish to estimate the free energy difference $F_2^\star - F_1^\star$, which is equal to zero because $F_k^\star = -\log Z_k^\star = -\log \int_\mathcal{S} e^{-H_k(x)} dx$ and so $F_2^\star = F_1^\star = 0$ (see (1)). Let $\Delta F^t = F_2^t - F_1^t$, with $F_k^t$ defined by (15) for $k = 1,2$. Because $[K] = \{1,2\}$, we have $1_{\{1\}}(j) + 1_{\{2\}}(j) = 1$, and we compute

$$\Delta F^{t+1} = \Delta F^t + \frac{2}{t+1}(2 \cdot 1_{\{1\}}(\mathcal{K}^t) - 1). \tag{17}$$

We recall that the stationary density for $P_{\gamma,Z}$ is proportional to $\gamma_\ell \, e^{-H_\ell(\mathcal{X}^s)}/Z_\ell$ (see (2)). If we use the asymptotic transition kernel $P_{\gamma,Z^\star}$, then $\mathcal{K}^t$ is drawn independently from the *exact* asymptotic distribution, and $\Delta F^t$ defines the $t^{th}$ partial sum of twice the random harmonic series, which converges to a non-constant distribution.[24] Curiously then, at the level of typical stationary Monte Carlo sampling, the estimator $\Delta F^t$ is *not* consistent. On the other hand, $P_{\gamma,Z^t}^k(\mathcal{K}^{t+1} = 2) = \left(1 + e^{-\Delta F^t}\right)^{-1}$, so if $\mathcal{K}^{t+1}$ is drawn from $P_{\gamma,Z^t}^k$, whenever $\Delta F^t > 0$ it is more likely that $\mathcal{K}^{t+1} = 2$ and hence by (17) the estimate $\Delta F^t$ is more likely to decrease, while if $\Delta F^t < 0$ it is more likely that $\mathcal{K}^{t+1} = 1$ and hence the estimate $\Delta F^t$ is more likely to increase. It is this self-adjustment provided by the fluctuations in the stationary density $p_{\gamma,Z^t}$ (through its dependence upon $Z_2^t/Z_1^t \propto e^{-\Delta F^t}$) that ensures the almost sure convergence of $\Delta F^t$ to 0. We will see (in the example provided in Section 1.4) that the self-adjustment property holds for (16) as well. In the current example, however, (16), which is a variance-reduced form of (15) (see Theorem 2 in ref. 11), has zero variance because the two distributions $e^{-H_k(x)}/Z_k^\star dx$, $k = 1,2$ are identical.

Now we introduce the idea of iterative importance sampling (IIS), which turns out to be asymptotically equivalent to (16), and provides useful perspective on the TSS algorithm introduced in the following section. First we note that the simulated tempering importance sampling estimator (7) behaves very differently from (14). In particular, when dealt samples drawn from a stationary distribution, (7) consistently estimates the free energy difference *regardless* of the nature of the mixture distribution from which the samples are drawn. With fixed but arbitrary hypothesized $(\gamma, Z)$ and samples $\mathcal{X}^s$ drawn from the mixture density $p_{\gamma,Z}^S(x) \propto \sum_{\ell \in [K]} \gamma_\ell e^{-H_\ell(\mathcal{X}^s)}/Z_\ell$, eq. (7) can be re-arranged as

$$Z_k^t \equiv \frac{1}{t}\sum_{s=1}^{t} \frac{e^{-H_k(\mathcal{X}^s)}}{\sum_{\ell \in [K]} \gamma_\ell e^{-H_\ell(\mathcal{X}^s)}/Z_\ell} \to \frac{Z_k^\star}{\sum_{\ell \in [K]} \gamma_\ell Z_\ell^\star/Z_\ell}, \quad \text{as } t \to \infty. \tag{18}$$

As mentioned in the previous section, the efficiency of this method is acceptable only if the hypothesized partition functions are close to the exact ones. The form of (18), though, suggests that this issue of exponential sensitivity can be resolved by adjusting the hypothesized $Z_k$ iteratively by setting it to the previous estimate of $Z_k^\star$. More explicitly, we define the following iterative importance sampling (IIS) estimator

$$Z_k^t \equiv \frac{1}{t}\sum_{s=1}^{t} \frac{e^{-H_k(\mathcal{X}^s)}}{\sum_{\ell \in [K]} \gamma_\ell e^{-H_\ell(\mathcal{X}^s)}/Z_\ell^{s-1}}. \tag{19}$$

The IIS estimator converges with optimal asymptotic efficiency, as follows under certain conditions from Tan's work. The estimator takes the recursive form

$$Z_k^{t+1} = Z_k^t + \frac{1}{t+1}\left(\frac{e^{-H_k(X^{t+1})}}{\sum_{\ell\in[K]}\gamma_\ell\, e^{-H_\ell(X^{t+1})}/Z_\ell^t} - Z_k^t\right), \quad (20)$$

which in terms of the free energy $F_k = -\log Z_k$ reads

$$F_k^{t+1} = F_k^t - \log\left(1 + \frac{1}{t+1}\left[\frac{e^{F_k^t - H_k(X^{t+1})}}{\sum_{\ell\in[K]}\gamma_\ell\, e^{F_\ell^t - H_\ell(X^{t+1})}} - 1\right]\right). \quad (21)$$

The iteration (21) (and hence (19)) is asymptotically equivalent to (16), as can be inferred by considering the expansion $\log(1+x) = x + O(x^2)$ and neglecting the $O(x^2)$ term, and all three recursions (15), (16), and (21) feature the self-adjustment property demonstrated in the example above. Despite the asymptotic equivalence between (16) and (21), we will adopt the IIS formulation (19) (and hence the recursion (21)) because it makes transparent the interpretation of (20) as an importance sampling estimator when the adaptation of $Z^t$ is frozen (i.e., when the expressions defining $Z_k^t$ in (18) and (19) are the same). With this interpretation, the global convergence of (20) when samples are drawn from stationary mixtures follows from the Strong Law of Large Numbers, whereas the same convergence is generally more difficult to establish for stochastic approximation recursions like (16), and does not hold for estimators driven purely by self-adjustment (such as (15)). The same perspective on importance sampling will be used again in Section 2.3 to derive a history-forgetting mechanism.

In the following section, we will write the IIS estimator in stochastic approximation form by using (10) directly in terms of the partition function $Z = e^{-F}$, and we will favor $Z$ as the coordinate to be used for derivations, as they are the natural coordinates for the IIS form (19). On a computer, the state of the parameter will be stored as $F$ rather than $Z$, in order to cope with the limited range for the exponent of floating-point number formats, and so in our implementation of the TSS algorithm we will use (21).

## 1.3   The Times Square Sampling algorithm

In the previous section we reviewed how stochastic approximation can be used to estimate free energy differences on the fly. Despite this advance, on-the-fly estimation of free energy differences remains difficult in molecular simulations, owing to the immense complexity of the space being sampled, so further efficiency improvements are highly desirable. In this section we introduce the TSS algorithm, which provides a mechanism to optimally estimate auxiliary averages—namely, averages against $p_{\gamma,Z}(x,k)dx\mathcal{C}(dk)$ and its associated conditional and marginal distributions—on the fly, and a means by which to use these averages to guide the sampling process, thereby improving the estimation of free energy differences. (This approach of estimating auxiliary averages to improve the estimation of a quantity of interest is in the same spirit as the developing field of adaptive Monte Carlo.[20]) After briefly motivating the use of such auxiliary averages, we will describe in Section 1.3.1 three theoretical features of TSS, two of which use on-the-fly estimation of the auxiliary averages to guide the sampling of the rungs, while the third relates to the impact of the rung sampling on the asymptotic variance. In Section 1.3.2, we introduce the TSS algorithm itself as an extended triple of stochastic approximation recursions, and state mathematical results that illustrate how the features we have described improve the estimation of free energy differences.

To motivate the introduction of auxiliary averages, we note that the IIS estimator introduced in equation (19) and the unnormalized stationary density $\gamma_k e^{-H_k(x)}/Z_k$ both depend on the rung density $\gamma$, the selection of which affects the asymptotic variance of the free energy difference estimates $F_k^t - F_j^t$ (as defined by (21)). The specification of an optimal rung density $\gamma$ is outside the intended scope of this paper, but (as discussed in more detail below) choosing $\gamma_k$ equal to a constant for all $k \in [K]$ is unlikely to be optimal. Rather, it is natural to define $\gamma$ in terms of averages of functions of $x$ against the various $\rho_k(x)dx$ (which is the distribution $p_{\gamma,Z}(x,k)dx\mathcal{C}(dk)$ conditioned on $k$). These averages are (like the free energies) unknown a priori, and estimating them will require an additional set of stochastic approximation recursions. In Section 1.3.1, we will show how any rung density defined as a deterministic function of such averages can be optimally estimated on the fly.

The idea of using the rung density to guide the sampling is further extended in TSS. We use $\pi \in \mathcal{P}^\varepsilon(K)$ to denote a general rung density that depends not only on the averages that determine $\gamma$, but also on a third set of averages, this time against the marginal distributions $p_{\pi,Z}^{[K]}(k)\mathcal{C}(dk)$, which are averages of functions of $k$; estimating these averages on the fly in turn yields a third set of stochastic approximation recursions. Asymptotically, $\pi$ will be equal to $\gamma$, but there will also be a transient regime in which $\pi$ affects the *rate* at which the free energy estimates approach the equilibrium.

With the introduction of $\pi$, the recurrence (20) written in the stochastic approximation form (10) becomes

$$G_k^Z(x; \pi, Z) = \frac{e^{-H_k(x)}}{\sum_{\ell \in [K]} \pi_\ell\, e^{-H_\ell(x)}/Z_\ell} - Z_k, \tag{22}$$

with mean field (see (8))

$$g_k^Z(\pi, Z) = \mathbb{E}_{\pi,Z}\big[G_k^Z(\mathcal{X}; \pi, Z)\big] = \frac{Z_k^\star}{\sum_{\ell \in [K]} \pi_\ell\, Z_\ell^\star/Z_\ell} - Z_k. \tag{23}$$

This is the first of three stochastic approximation observables that will be defined in this section. A common source of difficulty for the three observables is that producing good estimates of the mean field (e.g., (23)) from the observable (e.g., (22)) is the very large dimensionality of the configuration space $\mathcal{S}$ of molecular systems. To ameliorate these difficulties, two design goals are suggested for an algorithm intended to sample the joint density $\gamma_k \rho_k(x)$ and estimate the free energy differences and auxiliary averages simultaneously. First, one should seek to minimize the asymptotic variance of the estimates, which requires that the stochastic approximation algorithm be applied with an optimal gain matrix. Second, because the space $\Lambda$ is tractable (whereas $\mathcal{S}$ is not), correlations between the rung process $\mathcal{K}^t$ and the free energy estimates $F^t$ should be leveraged to improve the quality of the free energy estimation by making the most of the samples $\mathcal{X}^t$. We will emphasize these design goals when we present the theoretical features of TSS in Section 1.3.1; in Section 1.3.2, we will then gather the required on-the-fly estimates into a triple of stochastic approximation recursions and present the TSS algorithm, along with theoretical results that demonstrate the improvements resulting from the features we have described.

### 1.3.1 Features of TSS

In this subsection we will describe three theoretical features of TSS, each of which makes a distinct contribution to efficiently sampling the joint density $\gamma_k \rho_k(x)$. The first two features involve on-the-fly estimation of auxiliary averages, and in each case the stochastic approximation form is chosen to have an identifiable optimal gain matrix, which enables the asymptotic variance of the estimates to be minimized. The third feature relates to the effect that the mixing properties of the sampler have on the asymptotic variance of the free energy difference estimates.

The first feature is the choice of asymptotic rung density $\gamma \in \mathcal{P}^\varepsilon(K)$ and its on-the-fly estimation, which will yield the second of the three stochastic approximation recursions of TSS. Although we will choose a concrete form for $\gamma$, our arguments below imply that any $\gamma \in \mathcal{P}^\varepsilon(K)$ defined as a deterministic function of averages against $\rho_k(x)dx$ can be estimated. Recall that $\gamma = (\gamma_1, \ldots, \gamma_K)$ is a discretization of a density originally defined over the continuous variable $\lambda \in \Lambda$. That density is often taken to be constant (and so $\gamma_k = 1/K$), but this choice is unlikely to be optimal; an optimal density must transform as a pseudotensor under coordinate transformations having $\Lambda$ as their image (its definition must be such that the density is multiplied by the absolute value of the Jacobian of the transformation). Lindahl et al.[13] have investigated the use of Riemannian metrics for the sampling of $\Lambda$, including the determinant of the Fisher information matrix

$$\gamma_k \propto \det\left(\left\{\mathbb{E}_k\left[\frac{\partial H}{\partial \lambda_i}\frac{\partial H}{\partial \lambda_j}\right] - \mathbb{E}_k\left[\frac{\partial H}{\partial \lambda_i}\right]\mathbb{E}_k\left[\frac{\partial H}{\partial \lambda_j}\right]\right\}_{1 \leq i,j \leq d}\right)^{1/2}, \qquad (24)$$

which we take as the default TSS reference rung density (see also ref. 25). Note, though, that other choices are possible, as the function $H$ in the formula above can be replaced by other functions with preservation of the coordinate transformation law. The formula expresses $\gamma_k$ as a deterministic smooth function of a collection of auxiliary averages $\mu^\star = \{\mu^\star_{ki} = \mathbb{E}_k[\psi_i(\mathcal{X})], i = 1, \ldots, M, k \in [K]\}$ of functions $\psi_i(x)$ against the distributions $\rho_k(x)dx$. We also note that (24) illustrates a lack of dependence on the partition functions $Z^\star$.

In general we do not know the exact averages $\mu^\star$ that would allow us to use the target rung density $\gamma_k(\mu^\star)$, so we instead use $\gamma_k(\mu)$, where $\mu = \{\mu_{ki}; i = 1, \ldots, M, k \in [K]\}$ is a collection of hypothesized parameters. These parameters $\mu$ represent our estimates of $\mu^\star$ obtained through the stochastic simulation algorithm, which leads to the second of the three of stochastic approximation equations. The functional dependence given by (24) may indicate that there is no $\varepsilon > 0$ such that $\gamma(\mu) \in \mathcal{P}^\varepsilon(K)$ for all possible values of $\mu$. To prevent this situation, we re-define $\gamma(\mu)$ through the $\epsilon_\gamma$-regularization

$$\gamma_k(\mu) := \frac{(1-\epsilon_\gamma)\gamma_k(\mu) + \epsilon_\gamma \mathcal{M}(\mu)}{\sum_{\ell \in [K]}[(1-\epsilon_\gamma)\gamma_\ell(\mu) + \epsilon_\gamma \mathcal{M}(\mu)]}, \quad \forall \, 1 \leq k \leq K, \tag{25}$$

where $\mathcal{M}(\mu) = \max_{k \in [K]} \gamma_k(\mu)$ and $\gamma(\mu)$ is unregularized if appearing on the right-hand side. The value $\epsilon_\gamma \in (0,1]$ is taken to be small relative to 1. The $\epsilon_\gamma$-regularization implies that $\gamma(\mu) \in \mathcal{P}^\varepsilon(K)$ for all $\mu$ and $\varepsilon \leq \epsilon_\gamma K^{-1}$.

Instead of $\mu$, we work in terms of coordinates collectively represented by $\xi = Z\mu = e^{-F}\mu$ and explicitly given by $\xi_{km} = Z_k \mu_{km} = e^{-F_k}\mu_{km}$, for all $k \in [K]$ and $m = 1, \ldots, M$. To estimate the averages $\mu_{km}^\star$ we define the family of observables

$$G_{km}^\xi(x; \pi, Z, \xi) = \frac{e^{-H_k(x)}\psi_m(x)}{\sum_{\ell=1}^K \pi_\ell e^{-H_\ell(x)}/Z_\ell} - \xi_{km}, \tag{26}$$

with mean field $g_{km}^\xi(\pi, Z, \xi) = \mathbb{E}_{\pi,Z}\left[G_{km}^\xi(\mathcal{X}; \pi, Z, \xi)\right]$ given by

$$g_{km}^\xi(\pi, Z, \xi) = \frac{Z_k^\star \mathbb{E}_k[\psi_m(X)]}{\sum_{\ell \in [K]} \pi_\ell Z_\ell^\star/Z_\ell} - \xi_{km}. \tag{27}$$

A notable property of the mean fields (23) and (27) is their *independence* from $\pi$ when evaluated at $Z = Z^\star$. As a first consequence, the unique root of $g_{km}^\xi(\pi, Z^\star, \xi) = 0$ is $\xi_{km}^\star = Z_k^\star \mu_{km}^\star$ for all $k \in [K]$, and $m = 1, \ldots, M$. As a second consequence (which anticipates Proposition 1), the optimal gain matrix for the stochastic approximation algorithm in $\xi$ is (as we will see) the identity matrix.

Second, we introduce an important feature that operates at the level of the mean field $g_{\pi,Z}^Z$, defined in (23), rather than at the level of fluctuations about the equilibrium $Z^\star$: by controlling the mean field, this feature helps estimates reach a neighborhood of the equilibrium more efficiently, and will yield the third stochastic approximation recursion. The mean field is changed through a modification (or "tilt") in the rung density, and the tilt is chosen to bias the sampling toward under-visited rungs, which contributes to our second design goal. If instead we were to use $\gamma(\mu)$ directly in the makeup of the sampling distribution $p_{\pi,Z}(x,k)\,dx\,\mathcal{C}(dk)$ (which, to our knowledge, is how it is commonly used in both simulated tempering and adaptive Monte Carlo), then when $Z$ is far away from $Z^\star$, the marginal rung occupancies would severely tilt away from the desired reference $\gamma(\mu)$, which would thereby hinder the sampling of the rung process and the estimation of the free energy differences, slowing their approach to equilibrium. The $F$-dependence on the right-hand side of (36) illustrates the extreme sensitivity to errors in the free energy. The errors tend to be large, given that free energies are extensive properties (i.e., properties that scale linearly with the dimensionality of $\mathcal{S}$). Observe that, if instead of $\pi$ in (36), we were to use

$$\pi_k^\circ(Z,\xi) = \frac{Z_k/Z_k^*}{\sum_{\ell=1}^{K}\gamma_\ell(\mu)\,Z_\ell/Z_\ell^*}\gamma_k(\mu), \tag{28}$$

then we would actually recover $\gamma(\mu)$ as the marginal even with $Z$ being away from equilibrium. Of course, $\pi_k^\circ$ is not available, since it assumes knowledge of the exact partition functions. TSS proposes to use (6)—its sibling (7) has already been used up in the form of (22)—as a control variate measuring and penalizing deviations of $\pi$ from $\gamma(\mu)$, with the penalty vanishing at the equilibrium $Z^\star$. To achieve this, we introduce the *tilt* parameters $o = (o_1, \ldots, o_K)$, which will be used to scale $\gamma(\mu)$ through a deterministic function $\pi(Z,\xi,o): \mathbb{R}_{\geq 0}^K \times \mathbb{R}^{KM} \times \mathbb{R}_{\geq 0}^K \to \mathcal{P}^\varepsilon(K)$, and which are estimated using the observable

$$G_k^o(j;\pi,Z,\xi,o) = \frac{1}{\gamma_k(\mu)}\mathbf{1}_{\{k\}}(j) - o_k. \tag{29}$$

The mean field

$$g_k^o(\pi(Z,\xi,o),Z,\xi,o) = \mathbb{E}_{\pi,Z}\left[\frac{1}{\gamma_k(\mu)}1_{\{k\}}(\mathcal{K}) - o_k\right] = \frac{1}{\gamma_k(\mu)}\frac{\pi_k(Z,\xi,o)\,Z_k^*/Z_k}{\sum_{\ell=1}^K \pi_\ell(Z,\xi,o)\,Z_\ell^*/Z_\ell} - o_k \quad (30)$$

is different from the special form assumed by (23) and (27) in that its dependence upon $Z_k^*$ is not solely through the normalization constant, except for special cases, including the case $\pi = \pi°$. Proposition 1, which will be stated in the next subsection, demonstrates that the optimal gain matrix is trivially available for this special case. For the sake of efficiency, it is important to at least approximate $\pi°$, and TSS uses the relative empirical frequencies $o$ to achieve this desiderate through a feedback mechanism we call *visit control*.

Visit control relies on the two-timescale separation mentioned in Section 1, which we now express more concretely to mean that the ODE $\dot{o}_k(t) = g_k^o(Z(t),\xi(t),o(t))$, which governs the average behavior of the estimates away from equilibrium, evolves at its steady state with respect to the tilt parameters $o$. More explicitly, for any deterministic, smooth function $\pi: \mathbb{R}_{\geq 0}^K \times \mathbb{R}^{KM} \times \mathbb{R}_{\geq 0}^K \to \mathcal{P}^\varepsilon(K)$, and fixed $Z \in \mathbb{R}_{\geq 0}^K$ and $\xi \in \mathbb{R}^{KM}$, the steady-state condition for the tilts $o = (o_1,\ldots,o_K)$ is obtained by setting $g_k^o(\pi(Z,\xi,o),Z,\xi,o) = 0$ (defined in (30)), which defines $o_k^\star(\pi,Z)$. We cannot use the steady-state condition to determine $o$, however, since $g_k^o(\pi,Z,\xi,o)$ depends not only on the estimate $Z$ but also on the unknown partition functions $Z^\star$. Instead, we use stochastic approximation with the observable $G_k^o(j;\pi,Z,\xi,o)$ (defined in (29)) to provide an estimate of the tilt parameters $o$, which we then *assume* are at the steady state, in accordance with the two-timescale separation. In general our estimate $o$ will not be exactly at the steady-state ($o_k \neq o_k^\star(Z,\xi)$), so the steady state condition $g_k^o(\pi(Z,\xi,o),Z,\xi,o) = 0$ will be used to define *new* state variables $Z_k^\circ$. After re-arranging, this steady-state condition becomes

$$\frac{\gamma_k(\mu)o_k}{\sum_{\ell=1}^K \gamma_\ell(\mu)o_\ell} = \frac{\pi_k(Z,\xi,o)\,Z_k^\circ/Z_k}{\sum_{\ell=1}^K \pi_\ell(Z,\xi,o)\,Z_\ell^\circ/Z_\ell}. \quad (31)$$

(In the above equation, the variable $\mu$ is a function of $\xi$. Here and elsewhere the functional dependence will not be given explicitly.) When the tilt parameters are actually at their steady state ($o_k = o_k^\star(Z, \xi)$) then $Z_k^\circ = Z_k^\star$, but otherwise the variables $Z_k^\circ$ are simply functions of $(Z, \xi, o)$, which we refer to as visit control partition functions, and they define the visit control free energies by way of $F_k^\circ = -\log(Z_k^\circ)$.

Having obtained this auxiliary estimate $Z_k^\circ$ for $Z_k^\star$, we can now provide an approximation for $\pi^\circ$ in (28). TSS proposes

$$\pi_k^{\text{TSS}}(Z, \xi, o) = \frac{(Z_k/Z_k^\circ)^{\eta/(\eta+1)}}{\sum_{\ell=1}^K \gamma_\ell(\mu)(Z_\ell/Z_\ell^\circ)^{\eta/(\eta+1)}} \gamma_k(\mu), \tag{32}$$

for some $\eta > 0$, thus defining the feedback mechanism we called visit control. By combining (31) and (32) we derive the explicit form

$$\pi_k^{\text{TSS}}(Z, \xi, o) = \frac{\gamma_k(\mu)/(o_k)^\eta}{\sum_{\ell \in [K]} \gamma_\ell(\mu)/(o_\ell)^\eta}. \tag{33}$$

The case $\eta = 0$ corresponds to no visit control (as in SAMC and SAMS), whereas the other extreme $\eta = \infty$ corresponds to "infinitely strong" visit control, for which $\pi^{\text{TSS}}$ is formally equal to $\pi^\circ$, but which is unattainable because it would require $o_k = 1$, by way of (31). In practice we use $\eta \in (0, \infty)$, with larger $\eta$ creating a stronger bias to under-visited rungs. Additionally, we re-define $\pi_k^{\text{TSS}}(Z, \xi, o)$ through the $\epsilon_\pi$-regularization

$$\pi_k^{\text{TSS}}(Z, \xi, o) := (1 - \epsilon_\pi)\pi_k^{\text{TSS}}(Z, \xi, o) + \epsilon_\pi \gamma_k(\mu), \tag{34}$$

where the symbol $\pi_k^{\text{TSS}}(Z, \xi, o)$ appearing in the right-hand side of (34) denotes the right-hand side of (33) and $\gamma_k(\mu)$ denotes the right-hand side of (25). We note that

$$\pi_k^{\text{TSS}}(Z, \xi, o) \geq \epsilon_\pi \gamma_k(\mu) \geq \epsilon_\pi \epsilon_\gamma K^{-1}, \tag{35}$$

so $\pi^{\text{TSS}}(Z, \xi, o) \in \mathcal{P}^\epsilon(K)$ for all $(Z, \xi, o)$ and $\epsilon \leq \epsilon_\pi \epsilon_\gamma K^{-1}$. With this choice, the $\epsilon_\pi$-regularization is certain to leave $\gamma_k(\mu)$ unchanged if applied to it. The regularization is needed in order to ensure a minimum visiting propensity for rungs that could otherwise acquire very small unregularized weights through large fluctuations of $o$.

The visit-control distribution $\pi^{\text{TSS}}$ is used in the denominator of (22), thereby affecting the mean field (23). We show in Proposition 2 that with $o$ at its steady state $o^\star(Z, \xi)$ and any $\eta > 0$, the mean fields (23) converge faster to the equilibrium $(Z_1^\star, \ldots, Z_K^\star)$ (up to a multiplicative constant), as measured by the rate of decrease of a naturally associated Lyapunov function defined later in (43). The effect of the visit-control mechanism disappears at the equilibrium $(Z^\star, \xi^\star, \mathbf{1})$, where we note that $\mathbf{1}$ is the steady state of $o$ when $Z = Z^\star$.

Finally, we present the third of the theoretical features, the *self-adjustment* property (which we first illustrated in an example in Section 1.2). This property is not an explicit use of an auxiliary average to guide the sampling, but is instead an implicit relationship between the observable for the partition functions (22) and the invariant density (2), which allows for in some cases substantial variance reduction. Although the choice of $\gamma$ determines the asymptotic distribution of the process $\mathcal{K}^t$, and thus modifies the asymptotic variance of the free energy difference estimates, variance reduction from self-adjustment is instead tied to the mixing properties of the transition kernel $\mathcal{P}_{\pi,Z}$. Recall that self-adjustment drives convergence of the free energy estimator using only the form of the stationary distribution of $\mathcal{P}_{\pi,Z}$ for fixed $\pi, Z$; this was first discussed in Section 1.2, and will be explicitly illustrated in the uniform distributions example in Section 1.4. The stationary distribution for $\mathcal{P}_{\pi,Z}$ is the TSS sampling distribution $p_{\pi,Z}(x, k)\, dx\, \mathcal{C}(dk)$, and from the right-hand side of (6), the marginal probabilities for rung occupancy are

$$\frac{Z_k^\star/Z_k}{\sum_{\ell\in[K]} \pi_\ell Z_\ell^\star/Z_\ell}\pi_k = \frac{e^{F_k-F_k^\star}}{\sum_{\ell\in[K]} \pi_\ell e^{F_\ell-F_\ell^\star}}\pi_k. \tag{36}$$

The marginal distribution (36) favors selecting the rung $k$ with the greatest discrepancy $F_k - F_k^\star$, weighted by $\pi$. This self-adjustment property inherent to the weighted counting estimator (15) carries over to its variance-reduced counterparts (16) and (21), and in both cases the variance of the free energy estimates is tied to the integrated autocorrelation time of the process $(\mathcal{X}^t, \mathcal{K}^t)$. To preserve the effect of self-adjustment and keep the variance as low as possible, we should make certain the contribution of the process $\mathcal{K}^t$ to the integrated autocorrelation time is as small as possible.

One way of preventing rung selection from becoming a slow process is to intersperse multiple updates of $\mathcal{K}^t$ in between updates to the free energy estimates. In the limiting case in which independent samples from $\rho_k(x)dx$ can also be cheaply generated and the kernel $P_{\pi,Z}$ applied multiple times between free energy updates, the impact of variance reduction from self-adjustment will be shown to be substantial enough that the on-the-fly estimator outperforms post-sampling estimators. Even optimal estimators such as MBAR have their variance bounded below by the Cramér-Rao bound;[26] the self-adjusted on-the-fly estimation procedure, which samples from a different distribution $p_{\pi,Z^t}(x,k)\mathcal{C}(dk)dx$ every $t \geq 0$, does not.

More concretely, if we fix the rung density $\pi \in \mathcal{P}^\varepsilon(K)$, given $Z^t$ at time $t \geq 0$, and we draw $(\mathcal{X}^{t+1}, \mathcal{K}^{t+1})$ independently from $p_{\pi,Z^t}(x,k)dx\mathcal{C}(dk)$, then the asymptotic variance of the free energy differences $F_k^t - F_j^t$ is strictly smaller than that of the maximum-likelihood estimator MBAR. A complete statement is in Proposition 3. The rate of sampling the space $\Lambda$ will be raised again in Section 2.1, where we introduce an important windowing system that restricts the movement of the process $\mathcal{K}^t$.

*1.3.2 Stochastic approximation triple and theoretical results*

In this section, we present the basic TSS algorithm. The form of the presentation is useful for theoretical investigation, and after placing the estimation of $\theta = (Z, \xi, o)$ within the framework of stochastic approximation outlined in Section 1.2, we are able to state a convergence theorem and provide theoretical results on the optimality of the stochastic approximation procedure, the benefits of visit control and the variance reduction due to self-adjustment. A more practical implementation of TSS, which builds on the presentation in this section, will be presented in Section 2.

The parameter space of $\theta$ is $\Theta = \mathbb{R}^K \times \mathbb{R}^{KM} \times \mathbb{R}^K_{\geq 0}$, and the state space is $\mathbb{Y} = \mathcal{S} \times \Lambda$, with $\mathcal{Y}^t = (\mathcal{X}^t, \mathcal{K}^t)$. Because $\pi$ depends deterministically on $\theta$, but not on $(x, k)$, the expectations in (23), (27), and (30) remain unchanged whether $\pi \in \mathcal{P}^\varepsilon(K)$ is fixed or $\pi = \pi(\theta)$ for fixed $\theta \in \Theta$. We will continue to write $\mathbb{E}_{\pi,Z}$ and $p_{\pi,Z}$ when there is no dependence of $\pi$ on $\theta$ or when it is irrelevant, and we will write $\mathbb{E}_\theta$ when we wish to emphasize the dependence on $\theta$; similarly, we will write $p_\theta$ for $p_{\pi(\theta),Z}$ and $P_\theta(x', k'|x, k)$ for $P_{\pi(Z,\xi,o),Z}(x', k'|x, k)$. We note that $\pi^{\text{TSS}}$ does not depend on $Z$, but this lack of dependence is immaterial for the proofs, so we allow $\pi$ to be a general deterministic function of $(Z, \xi, o)$. Also recall our preference for using the coordinates $(Z, \xi, o) \equiv (e^{-F}, e^{-F}\mu, o)$ for deriving the stochastic approximation equations, and for using the coordinates $(F, \mu, o)$ for storing information on a computer. The observables (22), (26), and (29) lead to the following equations for TSS: For each $k \in [K]$ and $m \in 1, \ldots, M$,

$$Z_k^{t+1} = Z_k^t + \frac{1}{t+1}\left(\frac{e^{-H_k(\mathcal{X}^{t+1})}}{\sum_{\ell=1}^K \pi_\ell(\theta^t)\, e^{-H_\ell(\mathcal{X}^{t+1})}/Z_\ell^t} - Z_k^t\right) \tag{37}$$

$$\xi_{km}^{t+1} = \xi_{km}^t + \frac{1}{t+1}\left(\frac{e^{-H_k(\mathcal{X}^{t+1})}}{\sum_{\ell=1}^K \pi_\ell(\theta^t)\, e^{-H_\ell(\mathcal{X}^{t+1})}/Z_\ell^t}\psi_m(\mathcal{X}^{t+1}) - \xi_{km}^t\right) \tag{38}$$

$$o_k^{t+1} = o_k^t + \frac{1}{t+1}\left(\frac{1}{\gamma_k(\theta^t)}\mathbf{1}_{\{k\}}(\mathcal{K}^{t+1}) - o_k^t\right). \tag{39}$$

Alternatively, if expressed in terms of $(F \quad \mu \quad o)^\top$,

$$F_k^{t+1} = F_k^t - \log\left(1 + \frac{1}{t+1}\left[\frac{e^{F_k^t - H_k(\mathcal{X}^{t+1})}}{\sum_{\ell=1}^K \pi_\ell(\theta^t) e^{F_\ell^t - H_\ell(\mathcal{X}^{t+1})}} - 1\right]\right) \tag{40}$$

$$\mu_{km}^{t+1} = e^{F_k^{t+1} - F_k^t}\left[\mu_{km}^t + \frac{1}{t+1}\left(\frac{e^{F_k^t - H_k(\mathcal{X}^{t+1})}}{\sum_{\ell=1}^K \pi_\ell(\theta^t) e^{F_\ell^t - H_\ell(\mathcal{X}^{t+1})}} \psi_m(\mathcal{X}^{t+1}) - \mu_{km}^t\right)\right] \tag{41}$$

$$o_k^{t+1} = o_k^t + \frac{1}{t+1}\left(\frac{1}{\gamma_k(\theta^t)}\mathbf{1}_{\{k\}}(\mathcal{K}^{t+1}) - o_k^t\right). \tag{42}$$

Note that owing to (40), equation (41) produces zero variance for $\psi_m(x) = const$.

**Proposition 1.** The duple of recursions defined by (37) and (38) has the minimal asymptotic variance among all gain matrices, regardless of the rung-sampling probability $\pi$. If in addition $\pi = \pi°$, then the triple formed by adjoining (39) also has minimal asymptotic variance.

The proof can be found in Section 1 of the Supplementary Materials. Because the stochastic approximation algorithm with optimal gain matrix is invariant under coordinate changes, the analysis of the mean fields and their Jacobians (carried out in Proposition 1; see Section 1.2 of the SM) is valid in any coordinate system—and thus applies to both (37)–(39) and (40)–(42)— once the optimality is demonstrated. Away from the point of stability, the choice of coordinates obviously matters.

We can now describe the general TSS algorithm.

*Algorithm 1*

0) Fix $(\mathcal{X}^0, \mathcal{K}^0) \in \mathcal{S} \times [K]$ and any $\theta^0 \in \Theta$.

1) For $t \geq 0$, generate $\mathcal{K}^{t+1}$ from $P_{\theta^t}^k(\cdot | \mathcal{X}^t)$, and then $\mathcal{X}^{t+1}$ from $P_{\theta^t}^x(\cdot | \mathcal{K}^{t+1})$.

2) Update $\theta^t = (F^t, \mu^t, o^t)$ to $\theta^{t+1}$ using equations (40), (41), (42) and go to step 1.

The main object we use to establish a convergence theorem for Algorithm 1 is the relative entropy

$$V_\gamma(Z) = D_{\mathrm{KL}}(\gamma || r) = -\sum_{k \in [K]} \gamma_k \log\left(\frac{r_k}{\gamma_k}\right), \tag{43}$$

with

$$r_k = \gamma_k \frac{Z_k^\star}{Z_k} \bigg/ \left(\sum_{\ell \in [K]} \gamma_\ell \frac{Z_\ell^\star}{Z_\ell}\right). \tag{44}$$

The relative entropy is defined for any $\gamma \in \mathcal{P}^\varepsilon(K)$ and all $Z \in \mathbb{R}_{\geq 0}^K$; $V_\gamma(Z)$ is non-negative and equal to zero only when $Z_k = c\, Z_k^\star$ for some common multiplicative constant $c > 0$. This function serves as a global Lyapunov function for the mean field (23) *uniformly* in $\pi \in \mathcal{P}^\varepsilon(K)$, which allows us to establish convergence of the system $\dot{Z}(t) = g_k^Z(\pi, Z(t))$ in the more general case that allows for dependence of $\pi$ on $\theta$.

Despite the availability of a global Lyapunov function, it is non-trivial to provide general conditions under which the sequence $\{Z^t, t \geq 0\}$ is globally convergent with probability 1. The main theoretical approaches used to establish almost-sure convergence globally are known as "reprojection" techniques;[27,28] these are discussed in the Supplementary Materials. In our statement of Theorem 1, we assume that the space $\Theta$ is compact.

**Theorem 1.**

Let $\pi$ be a non-vanishing continuously differentiable function of $\theta$, with continuous derivative on the compact set $\Theta$, and assume $o^\star = \mathbf{1}$ is the unique globally asymptotically stable equilibrium point of $\dot o_k = \pi_k(Z^\star, \xi^\star, o)/\gamma_k(\xi^\star) - o_k(t)$. For $t \geq 0$, let $(\mathcal{X}^t, \mathcal{K}^t)$ be the stochastic process produced as in Algorithm 1. Then we have

I. $F_k^t - F_\ell^t \to F_k^\star - F_\ell^\star$ as $t \to \infty$ almost surely;

II. $\mu_{km}^t \to \mathbb{E}_k[\psi_m(\mathcal{X})]$ as $t \to \infty$ almost surely;

III. $\pi_k^t \to \gamma_k(\mu^\star)$ as $t \to \infty$ almost surely, and $(\mathcal{X}^t, \mathcal{K}^t) \Rightarrow p_{\gamma(\mu^\star), F^\star}(x, k) dx \mathcal{C}(dk)$, where $\Rightarrow$ denotes convergence in distribution (see Supplementary Materials).

We now present Propositions 2 and 3, which relate to visit control and self-adjustment, respectively; an analytically tractable example that illustrates the effect of both these mechanisms follows in Section 1.4. The proofs of Propositions 2 and 3 can be found in the Supplementary Materials.

**Proposition 2.**

Fix an arbitrary rung density $\gamma: \mathbb{R}^{KM} \to \mathcal{P}^\varepsilon(K)$ and let $o^\star(Z, \xi)$ denote the steady state of the mean field of $o$ for fixed $Z, \xi$ and variable $\eta \geq 0$. Let $Z(t)$ denote the solution of the ODE $\dot Z(t) = g_k^Z(\pi^{\text{TSS}}(Z(t), \mu(t), o^\star(Z(t), \xi(t))), Z(t))$, where $g_k^Z(\pi, Z)$ is defined in (23) and $\pi^{\text{TSS}}(Z, \mu, o)$ is defined in (33). Then, with $V_\gamma(Z)$ defined in (43),

$$\left.\frac{d}{dt}V_\gamma(Z)\right|_{\eta'} \leq \left.\frac{d}{dt}V_\gamma(Z)\right|_{\eta} \leq 0, \quad \forall t \geq 0, \eta' > \eta \geq 0. \tag{45}$$

Both inequalities are strict unless $Z$ is a scalar multiple of $Z^\star$, in which case the time derivatives are equal to 0.

**Proposition 3.**

Let $\pi \in \mathcal{P}^\varepsilon(K)$, and suppose MBAR uses a proportion $\pi_k > 0$ of independent samples from $\rho_k(x)dx$, whereas for each fixed $F \in \mathbb{R}^K$ and the same $\pi$, TSS uses independent samples from $p_{\pi,F}(x,k)dx\,\mathcal{C}(dk)$. Then, for any $i,j \in [K]$, the asymptotic variance of $F_i - F_j$ for TSS is smaller than that of MBAR.

## 1.4 Example

The example in this section is intended to illustrate the self-adjustment procedure and the corresponding result of Proposition 3, as well as the visit-control mechanism and the corresponding result of Proposition 2. Consider two overlapping uniform distributions of equal width (Figure 1), formally defined by $H_1(x) = \infty \cdot 1_{[-1+\delta,\delta]^c}(x)$ and $H_2(x) = \infty \cdot 1_{[-\delta,1-\delta]^c}(x)$, where $A^c$ denotes the complement of a set $A \subset \mathbb{R}$. Both uniform distributions have a width of 1, so $F_i^\star = 0$ for $i = 1,2$, and overlap on a region of width $2\delta$. We use $\gamma_1 = \gamma_2 = 1/2$, and we rewrite the linearized estimator (16) in terms of $\Delta^t = F_2^t - F_1^t$, which yields

$$\Delta^{t+1} = \Delta^t + \frac{1}{t+1} \cdot \frac{e^{-H_1(x^{t+1})} - e^{\Delta^t - H_2(x^{t+1})}}{\frac{1}{2}e^{-H_1(x^{t+1})} + \frac{1}{2}e^{\Delta^t - H_2(x^{t+1})}}. \tag{46}$$

(We use the linearized estimator here because it is better suited to algebraic manipulations, but we note that the properties we illustrate also hold for the TSS estimator (21).) We first illustrate the variance reduction from self-adjustment, which concerns the local fluctuations of $\Delta^t$, and so we take $\pi(Z, \xi, o) = \gamma(\mu)$ (i.e., independent of $o$). We consider two pairs of sampling and estimation schemes. The first is normally considered to be the best way of estimating the free energy differences: Simulate independent samples $\mathcal{X}^1, \ldots, \mathcal{X}^t$ from $\rho_1(x)dx$ and $\mathcal{X}^{t+1}, \ldots, \mathcal{X}^{2t}$ from $\rho_2(x)dx$ (in equal proportions by the choice of $\gamma$) and process the samples using MBAR (or, in this case, simply BAR). The asymptotic variance of the estimator is $\Sigma_{\text{MBAR}} = 2(1 - 2\delta)/\delta = 2\,p_\delta/\delta$, where $p_\delta = (1 - 2\delta)$ (see Section 3.3 of the SM).

The second scheme, which features self-adjustment, picks any initial $(\mathcal{X}^0, \mathcal{K}^0)$ and initial estimate $\Delta^0$, and generates $(\mathcal{X}^t, \mathcal{K}^t)$ using the transition kernel $P^N_{\gamma, Z^t}$ at time $t \geq 0$, where $Z^t = (1, e^{-\Delta^t})$, $P_{\gamma, Z}$ is the simulated tempering transition kernel, and $P^N_{\gamma, Z} = P_{\gamma, Z} \circ \ldots \circ P_{\gamma, Z}$ is its $N$th step iteration. As shown in Part 1, Section 3.3 of the SM, the asymptotic variance of the TSS estimator is

$$\Sigma^N_{\text{TSS}} = 4(1 - 2\delta)\left(1 + 2\frac{(1-2\delta)^N}{1-(1-2\delta)^N}\right) = 4p_\delta + 8\frac{p_\delta^{N+1}}{1-p_\delta^N}. \tag{47}$$

Observe that $\Sigma_{\text{MBAR}} = O(\delta^{-1})$, but for any $\delta > 0$, there is $N$ large enough that $\Sigma^N_{\text{TSS}} = O(1)$; in the limit $\delta \to 0$, the effect of self-adjustment can completely eliminate variance due to the small but finite overlap $2\delta$. In fact, whenever $N \geq 2$, $\Sigma^N_{\text{TSS}} < \Sigma_{\text{MBAR}}$, and post-processing the samples with MBAR would increase the variance of the estimate. In this example, all it takes to outperform MBAR is an extra iteration of the simulated tempering sampler between estimator updates.

Next, we illustrate how the visit control mechanism affects the mean field (23) of the free energy differences. Equation (46) clearly illustrates the exponential slowdown that may occur when the gain factor is taken to be proportional to $t^{-1}$, but this choice of gain factor is highly desirable given the lack of efficiency in sampling $\mathcal{S}$ in typical problems. If $\Delta^0$ deviates by a large quantity from the equilibrium value (in this case 0)—if $\Delta^0 = 100$, for instance—then the contribution of $H_1$ to (46) is suppressed for as long as $\Delta^t$ remains large (larger than ~10). During this time, $\Delta^t$ is decreased by $-2/(t+1)$ to produce $\Delta^{t+1}$; due to the slow, logarithmic rate of divergence of the harmonic series, however, it takes $\propto \exp(|\Delta^0|/2)$ steps to reach equilibrium.

The standard way of handling this difficulty is by increasing the gain factor, typically to $t^{-\beta}$ for $\beta \in (0.5,1]$.[29] This choice biases the Monte Carlo samples entering the free energy estimators toward the most recent samples in a substantial way; with this approach, not only is the rate of decay of the root-mean-square error slower than the celebrated $t^{-1/2}$, but the approach is unattractive for situations where the samples in $\mathcal{S}$ have long correlation times. Tan has proposed a two-stage algorithm (see equation 15 of ref. 11) in which a more aggressive gain factor is used only up to a certain "burn-in" time $t_0$, but it is not clear how the idea can be expanded for setups in which there are many competing times $t_0$ that must be automatically established. Another general way of coping with the problem of relaxing from large initial deviations is provided by proximal operators, as discussed in ref. 30.

TSS handles the problem of the exponential slowdown with the visit-control mechanism, which works in part by increasing the magnitude of the mean field of the free energy estimate (i.e., by modifying $g(\theta)$ in (11) itself) when $\Delta^t$ is far from the equilibrium $\Delta^\star$ ($Z^t$ and $Z^\star$ in the general case). To illustrate how this works, we compare the mean field with and without visit control: Using the expression for the mean field $g_k(F)$ in (12), the mean field corresponding to (46) is determined by $g_2(F) - g_1(F)$, where $F_2 - F_1 = \Delta$; this gives

$$\dot{\Delta}(t) = \dot{\Delta}(\Delta(t), o(t)) = \frac{1 - e^{\Delta(t)}}{\pi_1(o(t)) + \pi_2(o(t))e^{\Delta(t)}}. \tag{48}$$

Not using visit control is equivalent to saying that $\pi = \gamma(\mu)$, which implies $|\dot{\Delta}(t)| \leq \max(\gamma_1^{-1}, \gamma_2^{-1}) = 2$. This boundedness of the mean field prevents the stochastic discretization (46) from compensating for the small step size $dt = 1/(t+1)$, exponentially slowing the rate of convergence of the free energy difference estimates $\Delta^t$.

We will now show the effect of visit control (using the unregularized form: i.e., $\pi = \pi^{TSS}(Z, \xi, o)$, as defined in (33)) on the mean field. In order to evaluate $\pi = \pi^{TSS}(Z, \xi, o)$, we must invoke the two-timescale assumption introduced at the beginning of Section 1, which allows us to treat $o$ as being at its $\Delta$-dependent steady state (which we denote $o^\star(Z, \xi)$ in the general case). In the present case, with $o$ expressed in terms of $\Delta$ (by setting (30) equal to 0 and solving for $\Delta$), we find $o_k \propto e^{F_k/(\eta+1)}$, and evaluating $\pi^{TSS}$ (as defined in (33)) gives

$$\pi^{TSS}(\Delta) = \left( \frac{1}{1 + e^{-\beta\Delta}}, \frac{e^{-\beta\Delta}}{1 + e^{-\beta\Delta}} \right)^\top, \tag{49}$$

where $\beta = \eta/(\eta + 1)$. Using $\pi = \pi^{TSS}(\Delta)$ in (48), we find that

$$\dot{\Delta}(t) = -2 \sinh\left(\frac{\Delta}{2}\right) \frac{\cosh\left(\frac{\beta\Delta}{2}\right)}{\cosh\left(\frac{(1-\beta)\Delta}{2}\right)} = g(\Delta(t)). \tag{50}$$

Note that $g(\Delta) = -g(-\Delta)$, and that $g(\Delta) \to \pm\infty$ exponentially as $\Delta \to \mp\infty$. The dynamics have thus improved, pulling back to 0 exponentially fast as $|\Delta| \to \infty$, which compensates for the small step size and thus avoids the exponential slowdown.

## 2  Implementation

Despite the theoretical features laid out in Section 1.3, TSS is not yet able to handle the demands of many problems of interest to the chemical or biological physicist. The main reason is that the visit control mechanism, while satisfactorily addressing the exponential slowdown of estimation arising from the $O(t^{-1})$ gain in simple examples (as in Section 1.4), is only a partial solution in more complex problems. In Section 2, we complement the theoretical features with practical ones: We describe a system of partitioning the space $\Lambda$, called windowing, that is needed for the stability of the estimates and the computational efficiency of the algorithm; we address the issue of the sampler's dependence on initial conditions using a modified version of what is known as "burn-in" suitable for use in on the fly estimators; and we modify the algorithm to use multiple replicas (i.e., multiple parallel samplers that contribute to the same free energy estimate) to take advantage of distributed computational resources. In Section 2.4, we then use a flow-chart to sequentially guide the reader through each step of the algorithm, showing how the features introduced in this section fit together and comparing the relative computational cost of different parts of the algorithm.

To explain these key algorithmic features, we must first understand the nature of the slowdown that TSS as described to this point cannot overcome: Even with visit control, the empirical distribution of $\mathcal{K}^t$ may become skewed to a specific region of $\Lambda$ early in the simulation when the estimates are poor, and as we discussed in Section 1, a skewed sampling of $\Lambda$ may substantially hinder estimation. The skewed sampling of $\Lambda$ is exacerbated by the tight coupling of the $\mathcal{X}$ and $\mathcal{K}$ dynamics through the composition of samplers $P_{\pi,Z}^k \circ P_{\pi,Z}^x$. First, the coupling means that the $x$-conditioned density, sampled by $P_{\pi,Z}^k$, tends to be highly localized around an $x$-dependent mode, reducing mobility of $\lambda_{\mathcal{K}^t}$ and preventing proper exploration of the full space $\Lambda$. Second, the $\lambda_k$-conditioned density, sampled by $P_{\pi,Z}^x$ (typically molecular dynamics), is highly complex and difficult to equilibrate given present computational resources, implying that the early history of a TSS simulation is dependent upon hard-to-control initial conditions, with unknown decay times of the dependence.

There are a few approaches to overcoming the coupled problems of sampling $P_{\pi,Z}^k$ (which can be directly monitored) and sampling $P_{\pi,Z}^x$ (which exhibits unknown decay times and cannot be directly monitored), though none are silver bullets. Notably, the AWH algorithm uses a recursion that is very similar to (21), but modulates the weight of samples in the recursion through a "covering criterion," which (roughly speaking) changes the gain in the recursion depending on how many times the space $\Lambda$ has been covered by $\lambda_{\mathcal{K}^t}$. Such an approach is effective, but may require some trial and error, which is costly when running long-timescale MD simulations.

The features of our TSS implementation that tackle these coupled problems are a windowing system and a history-forgetting system. The windowing system (introduced in Sections 2.1 and 2.2) addresses the issue of exponential slowdown as it arises from the skewed distribution of $\mathcal{K}^t$ by subdividing the $\Lambda$ space into a collection of loosely connected individual TSS dynamics, each constrained to its own, sufficiently small domain, called a *window*. The windowing system makes the algorithm robust: When $\mathcal{K}^t$ is restricted to a set of the form $\{k: |\lambda_k - \lambda| < \delta\}$ for some $\lambda \in \Lambda, \delta > 0$, the variance of the free energy difference estimate $F_k - F_j$ goes to 0 and skew in the empirical distribution of $\mathcal{K}^t$ is eliminated as $\delta \to 0$ for all $k, j$ in that set. At the same time, the nature of the windowing system makes the algorithm very scalable, in large part by enabling the efficient use of multiple replicas (as defined in Section 2.3).

The history-forgetting mechanism (introduced in Section 2.3) addresses the issue of long decay times that arises in sampling $P_{\pi,Z}^x$ by continuously adapting the estimates to use a fixed fraction $\alpha \in (0,1)$ of the samples collected. In this way, the samples $\{\mathcal{X}^\tau: 0 \leq \tau < \alpha t\}$ do not participate explicitly in the formation of estimators for the current time $t$. This adaptation entails redefining the estimators, and in doing so we also modify TSS to use multiple replicas, in order to take further advantage of the wide availability of distributed computing resources. The algorithm is then summarized in Section 2.4, and we offer some guidance as to its use and implementation.

## 2.1 The windowed parametric invariant distribution

In this section we introduce the windowing system, which decomposes the global sampling of $\Lambda$ into a collection of local sampling problems. This section parallels Section 1.1, in that it introduces a family of extended probability densities depending on a pair of parameters $\pi$ and $F$, and for each such pair defines ergodic sampling dynamics for which the corresponding density is invariant. These local dynamics entail local estimates of $\theta = (F, \mu, o)$, which will then have to be stitched together to provide global estimates across all of $\Lambda$; in Section 2.2, we provide a procedure to determine global estimates based on the collection of window-local estimates.

The windowing system decomposes the parameter space $\Lambda$ into a collection $\{W_j; 1 \leq j \leq J\}$ of connected overlapping windows $W_j \subset \Lambda$, and separate TSS iterations are performed within each window. To link the windows together, it introduces window dynamics that allow the sampler to select a different window among those containing $\lambda_k$ for the current rung $k$. Instead of sampling a pair $(\mathcal{X}, \mathcal{K})$, we will thus now sample a triple $(\mathcal{X}, \mathcal{K}, \mathcal{J})$, with the specification that the stationary density conditioned on $\mathcal{J} = j$ is equal to the original stationary density $p_{\pi,F}(x, k)$ defined in (2). We provide a more detailed description for the sake of concreteness, but note that other variants are possible. The important attributes that must be respected are that the windows cover the entire parameter space $\Lambda$ and that the window dynamics are positively recurrent.

The windows are structured so that each $\lambda \in \Lambda$ belongs to exactly two windows, with the additional constraint that the symmetric matrix $\mathcal{W}$ with $i,j$ entry equal to 1 if $W_i \cap W_j \neq \emptyset$ and otherwise equal to 0 be irreducible (Figure 2). If we define the set $\text{win}(k) = \{j \in [J]: \lambda_k \in W_j\}$, it follows that the cardinality of $\text{win}(k)$ is 2 for each rung $k \in [K]$. For each $j \in [J]$, quantities such as $Z_{j;\cdot}$ and $\pi_{j;\cdot}$ (which are specified with a bullet in place of the second subscript) denote the sets of per-window quantities $\{Z_{j;k}; k \in [K]: \lambda_k \in W_j\}$ and $\{\pi_{j;k}; k \in [K]\}$, respectively. Quantities with no subscript denote the totality of components over all windows: $Z \equiv \{Z_{j;k}; j \in [J], k: \lambda_k \in W_j\}$, for example. Similar notation will be used for the per-window estimates in Sections 2.2 and 2.3. All per-window rung densities, such as $\pi_{j;\cdot}$, are defined so that $\pi_{j;k} \geq \varepsilon > 0$ for all $\lambda_k \in W_j$ and 0 for $\lambda_k \notin W_j$. The densities are also normalized, and thus $\sum_{k \in [K]} \pi_{j;k} = \sum_{k: \lambda_k \in W_j} \pi_{j;k} = 1$.

Next, we define the window process $\mathcal{J}^t$, which at any fixed time $t \geq 0$ is equal to one of the two possible indices in $\text{win}(\mathcal{K}^{t-1})$. The window $W_{\mathcal{J}^t}$ will be called the *active window*. The dynamics of $\mathcal{J}^t$ are determined by the transition density

$$\bar{P}^j_{\pi,Z}(x',k',j'|x,k,j) = \delta(x-x')\delta(k-k')\left(1 - 1_{\{j\}}(j')\right) 1_{\text{win}(k)}(j'), \tag{51}$$

which, given $k$, selects the other window in $\text{win}(k)$ with probability 1. (If during the early stages of the simulation the estimates $\pi_{j;\cdot}, F_{j;\cdot}$ are not yet defined for some window $W_j$, neither the window nor rung is updated—that is, neither of the transition kernels (51) or (52) are used—and instead $(k',j')$ is set to $(k,j)$; the rest of the algorithm continues normally, and $\pi_{j;\cdot}$ and $F_{j;\cdot}$ become defined in the following cycle.) This antithetic selection can of course be generalized to a random selection in which the other window is selected with probability $p \in (0,1]$, but the need for efficient rung sampling (as discussed in Section 1.3) motivates us to take the most aggressive choice $p = 1$. The transition density for the triple $(\mathcal{X}, \mathcal{K}, \mathcal{J})$ is defined by $\bar{P}_{\pi,Z} = \bar{P}^k_{\pi,Z} \circ \bar{P}^j_{\pi,Z} \circ \bar{P}^x_{\pi,Z}$, where

$$\bar{P}^k_{\pi,Z}(x',k',j'|x,k,j) \propto \delta(j'-j)\mathbf{1}_{W_j}(\lambda_{k'})P^k_{\pi_{j;\cdot},Z_{j;\cdot}}(x',k'|x,k) \tag{52}$$

has been restricted to only propose rungs within the active window, and $\bar{P}^x_{\pi,Z}(x',k',j'|x,k,j) = \delta(j'-j)P^x_{\pi,Z}(x',k'|x,k)$ is unchanged. These dynamics and the irreducibility condition on the symmetric matrix $\mathcal{W}$ ensure each window will be visited asymptotically, and together with (51) determine the parametric invariant density to take the form

$$\bar{p}_{\pi,Z}(x,k,j) = \bar{p}_{\pi,Z}(j|x,k)\bar{p}_{\pi,Z}(x,k) = \frac{1}{2}\mathbf{1}_{W_j}(\lambda_k)\bar{p}_{\pi,Z}(x,k), \tag{53}$$

where $\bar{p}_{\pi,Z}(j|x,k)$ denotes a conditional density and $\bar{p}_{\pi,Z}(x,k)$ denotes a marginal density. Equation (53) states that at stationarity each rung $k$ is assigned one of the two windows containing $\lambda_k$, with equal probability.

To complete the description of the parametric invariant density, we derive an explicit dependence of $\bar{p}_{\pi,Z}(x,k)$ on $\pi,Z$. To this end, set $\bar{p}_j \equiv \bar{p}_{\pi,Z}(j)$ (the marginal density of $\mathcal{J}$), and let $\bar{\mathbb{P}}_{\pi,Z}$ denote the distribution associated to the density $\bar{p}_{\pi,Z}$ (with respect to $dx\,\mathcal{C}(dj)\mathcal{C}(dk)$; see Section 1.1). The marginal density $\bar{p}_{\pi,Z}(x,k)$ is relatable to the complementary marginal $\bar{p}_j$ by way of the conditional probability decomposition

$$\bar{p}_{\pi,Z}(x,k) = \sum_{j\in[J]} \bar{p}_j\,\bar{p}_{\pi,Z}(x,k|j) = \sum_{j=1}^{J} \bar{p}_j \frac{\pi_{j;k}e^{F_{j;k}-F^\star_k}}{\sum_{\ell:\lambda_\ell\in W_j}\pi_{j;\ell}e^{F_{j;\ell}-F^\star_\ell}}\rho_k(x). \tag{54}$$

Multiplying (54) by $\mathbf{1}_{W_i}(\lambda_k)$, integrating over $x \in \mathcal{S}$, and summing over $k$ implies, for $(\mathcal{X},\mathcal{K},\mathcal{J})$ distributed according to $\bar{\mathbb{P}}_{\pi,Z}$ and each $i \in [J]$,

$$\bar{\mathbb{P}}_{\pi,Z}(\lambda_\mathcal{K} \in W_i) = \sum_{j=1}^{J} \bar{p}_j \bar{\mathbb{P}}_{\pi,Z}(\lambda_\mathcal{K} \in W_i|\mathcal{J}=j) = \sum_{j=1}^{J} \bar{p}_j \frac{\sum_{k:\lambda_k\in W_i\cap W_j}\pi_{j;k}e^{F_{j;k}-F^\star_k}}{\sum_{k:\lambda_k\in W_j}\pi_{j;k}e^{F_{j;k}-F^\star_k}}. \tag{55}$$

A different decomposition produces, for each $i \in [J]$,

$$\bar{p}_i = \mathbb{\bar{P}}_{\pi,Z}(\mathcal{J} = i | \lambda_{\mathcal{K}} \in W_i) P(\lambda_{\mathcal{K}} \in W_i) + \mathbb{\bar{P}}_{\pi,Z}(\mathcal{J} = i | \lambda_{\mathcal{K}} \notin W_i) \mathbb{\bar{P}}_{\pi,Z}(\lambda_{\mathcal{K}} \notin W_i)$$
$$= \frac{1}{2} \mathbb{\bar{P}}_{\pi,Z}(\lambda_{\mathcal{K}} \in W_i) + 0 \cdot \mathbb{\bar{P}}_{\pi,Z}(\lambda_{\mathcal{K}} \notin W_i), \tag{56}$$

which together with (55) leads to the equation

$$\frac{1}{2} \sum_{j=1}^{J} \frac{\sum_{k:\lambda_k \in W_i \cap W_j} \pi_{j;k} e^{F_{j;k} - F_k^\star}}{\sum_{k:\lambda_k \in W_j} \pi_{j;k} e^{F_{j;k} - F_k^\star}} \bar{p}_j = \bar{p}_i, \quad \forall \, 1 \le i \le J. \tag{57}$$

Note that, for each $i \in [J]$, the family of disjoint subsets $\{W_i \cap W_j; 1 \le j \le J, j \ne i\}$ form a partition of $W_i$. It follows then that (57) asks for the invariant right-eigenvector of a left-stochastic matrix, which determines $\bar{p}_j$ (using $\sum_{j \in [J]} \bar{p}_j = 1$) and completes the description of the invariant distribution $\bar{p}_{\pi,Z}(x, k, j) dx\, \mathcal{C}(dj) \mathcal{C}(dk)$. This distribution is analogous to the invariant distribution $p_{\pi,Z}(x, k) dx \mathcal{C}(dk)$ defined in Section 1.1, but now with windows.

## 2.2 Obtaining global estimates from local, per-window estimates

In this section we define the dependence of $\pi$ and $Z$ on the collection of per-window estimates $\theta \equiv \{\theta_{j;k}; j \in [J], k: \lambda_k \in W_j\}$, which requires defining the rung weights $\pi_k$, the partition functions $Z_k$, and the visit control partition functions $Z_k^\circ$ globally from the local estimates $\theta$ (defined for $j \in [J], k: \lambda_k \in W_j$). There are two results in this section: The first is the formulation in Section 2.2.1 of a convex optimization problem that allows us to produce, for a given $\theta$, a set of global visit control free energies $\{F_k^\circ, k \in [K]\}$ that are used for sampling, analogously to those defined in the case of a single window in Section 1.3. These free energies are useful for driving the sampling, but are too noisy for the purpose of providing a final estimate. To address this issue, in Section 2.2.2 we also formulate a second, simpler optimization problem that provides a way to unite the free energy estimates across windows with a lower noise (albeit in a way only applicable near convergence) for the purpose of reporting a set of global free energies $\{F_k^{TSS}, k \in [K]\}$ to the user, which we call the *reported* free energies. The per-window recursions that define $\theta_{j;k}^{t+1}$ for $j \in [J]$ and $k: \lambda_k \in W_j$ based on $(\mathcal{X}^t, \mathcal{K}^t, \mathcal{J}^t)$ and $\theta^t$ will be given in Section 2.3.

We first describe an element common to both the visit control and reported free energies, which is the asymptotic rung density $\gamma \equiv \{\gamma_{j;k}; j \in [J], k: \lambda_k \in W_j\}$. Although the densities $\gamma_{j;\cdot}(\mu)$ could be defined to vary with $j$ in an arbitrary fashion, in practice we take all these densities to be conditional densities of the common global density $\gamma(\mu)$ of components $\gamma_k(\mu)$ considered in equation (24) of Section 1. In what follows we only need the un-normalized functional dependence of $\gamma_k(\mu)$ on $\mu$ (e.g., of the kind expressed by the right-hand side of (24)), and its locality property that $\gamma_k(\mu)$ depends only on row $k$ of the original $K \times M$ matrix $\mu$. This locality implies that, for each $k \in [K]$, the function $\gamma_k$ initially defined on matrices of size $K \times M$ can be naturally restricted to a function defined on matrices of size $1 \times M$ (i.e., rows), and hence we write $\gamma_k(\mu_{j;\cdot})$ to mean this restricted function evaluated on the $k^{th}$ row of $\mu_{j;\cdot}$. This property allows for a simple definition of $\gamma_{j;k}$ as a function of the $J$ matrices $\{\mu_{j;\cdot}, j \in [J]\}$ of sizes $K_j \times M$, where $K_j$ is the number of rungs in $W_j$. We take

$$\gamma_{j;k}(\mu) = \frac{(1-\epsilon_\gamma)\gamma_k(\mu_{j;\cdot}) + \epsilon_\gamma \mathcal{M}(\mu)}{\sum_{\ell:\lambda_\ell \in W_j}[(1-\epsilon_\gamma)\gamma_\ell(\mu_{j;\cdot}) + \epsilon_\gamma \mathcal{M}(\mu)]} 1_{W_j}(\lambda_k), \quad \forall\, 1 \leq k \leq K, \tag{58}$$

where $\mathcal{M}(\mu) = \max_{k \in [K], j \in \text{win}(k)} \gamma_k(\mu_{j;\cdot})$ and $\epsilon_\gamma \in (0,1]$ is a positive regularization parameter that is small relative to 1. The parameter $\theta$ will thus be defined for each $j \in [J]$ and $k: \lambda_k \in W_j$ by a tuple of the form $\theta_{j;k} = (F_{j;k}, \mu_{j;k,1}, \ldots \mu_{j;k,M}, o_{j;k})$.

### 2.2.1 Global visit control free energies

To guide the reader in the definition of the global visit control free energies as a function of the per-window estimates $\theta = \{\theta_{j;k} : j \in [J], k : \lambda_k \in W_j\}$, we first sketch the general idea. In the single-window case we identified the visit control free energies explicitly in terms of the tilts, eliminating the need to refer to $F°$; in the case of multiple windows, however, we will need $F°$ as an intermediate in order to define visit control globally. The key step is to invoke the two-timescale separation described in Section 1.3 and impose the steady-state condition (31), thereby relating the tilts $o$ and density $\gamma$ to the density $\pi$ and visit control $F°$, the latter two being used for sampling. This relationship allows us to obtain the marginal density for $\mathcal{J}$, and in turn the marginal density for $\mathcal{K}$. Finally, upon specifying the dependence of $\pi_{j;k}$ on $\theta$ through our particular choice of visit control (32), the steady state relation (31) can be recast as a convex optimization problem, which leads to our definition of visit control free energies $\{F_k°, k \in [K]\}$. Notably, the complexity of the optimization problem scales with the number of windows, not the number of rungs (the latter typically being greater), which is important for efficiency.

To impose the steady-state condition, we begin by assuming that the tuple $(\mathcal{X}, \mathcal{K}, \mathcal{J})$ is drawn from the sampling density $\bar{p}_{\pi,F}(x, k, j)$ defined in (53) for fixed $\pi, F$, and that the tilts $o_{j;\cdot}$ are at their steady state for each $j \in [J]$. We recall that this is exactly how visit control was defined in Section 1.3. The steady state (31), averaged over all windows, gives the relation

$$\sum_{j=1}^{J} p_j \frac{\gamma_{j;k}(\mu) o_{j;k}}{\sum_{\ell : \lambda_\ell \in W_j} \gamma_{j;\ell}(\mu) o_{j;\ell}} = \sum_{j=1}^{J} p_j \frac{\pi_{j;k} e^{F_{j;k} - F_k°}}{\sum_{\ell : \lambda_\ell \in W_j} \pi_{j;\ell} e^{F_{j;\ell} - F_\ell°}}, \quad \forall\, 1 \leq k \leq K. \tag{59}$$

with to-be-determined marginal rung probabilities $p_j, j \in [J]$ and global visit control free energies $F_k°, k \in [K]$. As in the single-window case, when $o$ is at its steady state for any fixed $\pi, F$, the relation (59) correctly implies $F_k° = F_k^\star$; otherwise, (59) *defines* the visit control free energy $F_k^\circ$ up to an additive constant. The equation itself restates the conditional probability decomposition $\overline{\mathbb{P}}_{\pi,Z}(\mathcal{K} = k) = \sum_{j \in [J]} \overline{\mathbb{P}}_{\pi,Z}(\mathcal{K} = k|\mathcal{J} = j)\overline{\mathbb{P}}_{\pi,Z}(\mathcal{J} = j)$. Because the equality expressed by the equation holds at all times, we note (by multiplying by $1_{W_i}(\lambda_k)$ and summing over all $k \in [K]$, as in (54)–(57)) that the probabilities $\{p_j; 1 \le j \le J\}$ are the solution to the stochastic eigenvector problem

$$\frac{1}{2} \sum_{j \in [J]} \frac{\sum_{k:\lambda_k \in W_i \cap W_j} \gamma_{j;k}(\mu) o_{j;k}}{\sum_{k:\lambda_k \in W_j} \gamma_{j;k}(\mu) o_{j;k}} p_j = p_i, \qquad \forall \, 1 \le i \le J, \tag{60}$$

which can be solved by, for example, the QR method.[31] When $o$ is at its steady state, $F_k^\circ - F_k^\star$ is constant over all $k \in [K]$, and in this case the $\bar{p}_j$ determined from (57) are equal to the $p_j$ determined from (60).

We now use the specific form of $\pi_{j;k}^{\text{TSS}}$ as given by (32) (applied separately to each individual window) to relate the left-hand side of (59) to its right-hand side, as follows

$$\pi_{j;k}^{\text{TSS}}(\theta) = \frac{\gamma_{j;k}(\mu) e^{\frac{\eta}{\eta+1}(F_k^\circ - F_{j;k})}}{\sum_{\ell:\lambda_\ell \in W_j} \gamma_{j;\ell}(\mu) e^{\frac{\eta}{\eta+1}(F_\ell^\circ - F_{j;\ell})}}, \qquad \forall \, j \in [J], k \in \{\ell: \lambda_\ell \in W_j\}, \tag{61}$$

for some $\eta > 0$ and unknown visit control free energies $F_k^\circ, k \in [K]$. Combining everything, the relation (59) gives, for each $k \in [K]$ and unknowns $\{F_k^\circ, k \in [K]\}$

$$\begin{aligned} q_k &\equiv \sum_{j \in \text{win}(k)} p_j \frac{\gamma_{j;k}(\mu) o_{j;k}}{\sum_{\ell:\lambda_\ell \in W_j} \gamma_{j;\ell}(\mu) o_{j;\ell}} \\ &= \sum_{j=1}^J p_j \frac{\gamma_{j;k}(\mu) e^{(F_{j;k} - F_k^\circ)/(\eta+1)}}{\sum_{\ell:\lambda_\ell \in W_j} \gamma_{j;\ell}(\mu) e^{(F_{j;\ell} - F_\ell^\circ)/(\eta+1)}}, \qquad \forall \, 1 \le k \le K. \end{aligned} \tag{62}$$

The first equality in (62) defines $q_k$ as the left-hand side of (59), whereas the second equality defines a system of equations with unknowns $F_k^\circ$.

As for the case of a single window, with $\pi_{j;k}^{TSS}(\theta)$ and $\gamma_{j;k}(\mu)$ on the rhs of the following equation given by (61) and (58), respectively, we re-define $\pi_{j;k}^{TSS}(\theta)$ through the $\epsilon_\pi$-regularization

$$\pi_{j;k}^{TSS}(\theta) := (1 - \epsilon_\pi)\pi_{j;k}^{TSS}(\theta) + \epsilon_\pi \gamma_{j;k}(\mu), \qquad \forall k: \lambda_k \in W_j, \tag{63}$$

where $0 < \epsilon_\pi \leq 1$. We note that the $\epsilon_\pi$-regularization leaves $\gamma_{j;\cdot}(\mu)$ unchanged, if applied to it. The solution of (62) would be readily provided by

$$F_k^\circ = (\eta + 1)\log\left(\frac{1}{q_k} \sum_{j \in \text{win}(k)} p_j \gamma_{j;k}(\mu) e^{(F_{j;k} - f_j)/(\eta+1)}\right), \qquad \forall\, 1 \leq k \leq K, \tag{64}$$

if only we knew the window free-energy offsets

$$f_j = (\eta + 1)\log\left(\sum_{k: \lambda_k \in W_j} \gamma_{j;k}(\mu) e^{(F_{j;k} - F_k^\circ)/(\eta+1)}\right), \qquad \forall\, 1 \leq j \leq J. \tag{65}$$

We note that any solution $F^\circ$ to (62) also minimizes the convex function

$$\sum_{j=1}^{J} p_j \log\left(\sum_{k:\lambda_k \in W_j} \gamma_{j;k}(\mu) e^{(F_{j;k} - F_k^\circ)/(\eta+1)}\right) + \frac{1}{\eta + 1}\sum_{k=1}^{K} q_k F_k^\circ. \tag{66}$$

Except for an additive term equal to the Shannon entropy $-\sum_{k=1}^{K} q_k \log(q_k)$, the last expression is equal to

$$\frac{1}{\eta + 1}\sum_{j=1}^{J} p_j f_j + \sum_{k=1}^{K} q_k \log\left(\sum_{j \in \text{win}(k)} p_j \gamma_{j;k}(\mu) e^{(F_{j;k} - f_j)/(\eta+1)}\right), \tag{67}$$

as follows by substituting (64) and (65) into (66). The minimization is carried out subject to the constraint $\sum_{j=1}^{J} p_j f_j = 0$, since (67) does not change if all offsets $f_j$ are modified by the addition of one same number (details are provided in the Supplementary Materials, Part 2, Section 5.2.2). Furthermore, only the components with $p_j \neq 0$ are uniquely determined. We set $f_j = 0$ for the others, since their precise values do not affect the algorithm. We note once again that the number of windows $J$ is often much smaller than the number of rungs $K$, which helps mitigate the cost associated with this global step of the algorithm. Finally, we observe that the visit control free energies are heavily dependent on the tilts, which are estimated by counting events. Counting estimators are notoriously noisy, and thus we recommend that visit control free energies only be used for sampling, and that different estimates (to be defined in the following section) be reported to the user.

### 2.2.2 Global reported free energies

Windowed TSS evaluates a separate free energy $F_{j;k}^t$ for each window $W_j$ containing $\lambda_k$, raising the question of how to combine the different estimates available for the same rung in a single global set of estimates. The visit control mechanism provides one possible answer in the form of the unique set of global free energy estimates $\{F_k^\circ; k \in [K]\}$, but these estimates are rather noisy, since they depend on the many tilts $\{o_{j;k.}; j \in [J], k: \lambda_k \in W_j\}$, which are estimated by counting events. We can reduce the noise by assuming the algorithm has already converged, so the tilts have approached their asymptotic value of 1. The stationarity of the ODE for $o_{j;k}$ (30) at the equilibrium value of 1 implies

$$o_{j;k} = \frac{e^{(F_{j;k}-F_k^\circ)/(\eta+1)}}{\sum_{\ell:\lambda_\ell \in W_j} \gamma_{j;\ell}(\mu) e^{(F_{j;\ell}-F_\ell^\circ)/(\eta+1)}} = 1, \quad \forall j \in [J], k: \lambda_k \in W_j. \tag{68}$$

This many equations cannot be accommodated just by selecting $\{F_k^\circ; k \in [K]\}$, which forces us to select $\eta = \infty$. With the tilts set to unity, (60) becomes

$$\sum_{j \in [J]} \left( \sum_{k:\lambda_k \in W_i \cap W_j} \gamma_{j;k}(\mu) \right) p_j = p_i, \quad \forall i \in [J], \tag{69}$$

which defines the stochastic eigenvector $(p_1, \ldots, p_J)$. This eigenvector can be used to define a global estimate of the asymptotic rung density through the weighted average

$$\gamma_k^{\text{TSS}} \equiv \sum_{j \in \text{win}(k)} p_j \gamma_{j;k}(\mu), \qquad \forall\, k \in [K]. \tag{70}$$

We refer to $\gamma_k^{\text{TSS}}$ as the reported rung density. Similarly, by letting $\eta \to \infty$ in (64) with $q_k = \gamma_k^{\text{TSS}}$, we obtain the TSS reported free energies

$$F_k^{\text{TSS}} = \frac{1}{\gamma_k^{\text{TSS}}} \sum_{j \in \text{win}(k)} p_j \gamma_{j;k}(\mu) \left( F_{j;k} - f_j^{\text{TSS}} \right), \qquad \forall\, 1 \leq k \leq K. \tag{71}$$

The window free energy offsets $f_j^{\text{TSS}}$ are deduced by letting $\eta \to \infty$ in the system of equations obtained by setting the gradient of (67) to 0 and retaining the dominant terms. The offsets can be shown to solve the linear system

$$\sum_{j=1}^{J} (\delta_{ij} - t_{ij}) f_j^{\text{TSS}} = \sum_{k:\lambda_k \in W_i} \gamma_{i;k}(\mu) \left( F_{i;k} - \sum_{j \in \text{win}(k)} \frac{p_j \gamma_{j;k}(\mu)}{\gamma_k^{\text{TSS}}} F_{j;k} \right), \tag{72}$$

where $\delta_{ij}$ is the Kronecker symbol and $t_{ij}$ are the entries of the right-stochastic matrix

$$t_{ij} = \sum_{k:\lambda_k \in W_i \cap W_j} \gamma_{i;k}(\mu) \frac{p_j \gamma_{j;k}(\mu)}{\gamma_k^{\text{TSS}}}. \tag{73}$$

As in the case of the reported free energies, once the offsets $f_j^{\text{TSS}}$ are determined, the quantities $F_k^{\text{TSS}}$ can be reported to the user by way of (71). We note that solution of (72) is not uniquely determined, and that one equation must be dropped out and replaced by $\sum_{j \in [J]} p_j f_j^{\text{TSS}} = 0$. Additional details relating to the implementation of (72) within the TSS code are available in the Supplementary Materials, Part 2, Section 5.2.3.

## 2.3 History forgetting and multiple replicas: Multireplica stochastic approximation with burn-in

In this section we describe two features of our implementation of TSS inspired by traditional statistical practices, namely burn-in and use of multiple parallel samplers, modified to be compatible/useful for on-the-fly estimation. We adapted burn-in as a *history forgetting mechanism*, which also enables the estimation of asymptotic error bars, by making use of the fact that TSS is formulated as an IIS estimator (as in (19)). The same formulation also enables the use of multiple conditionally independent samplers, which we call *multireplica* TSS.

We first describe the history-forgetting mechanism. At any time $t$, the history of the simulation is divided into two parts: the early history, corresponding to times $s$ earlier than $\lfloor \alpha t \rfloor$ (the largest integer smaller than or equal to $\alpha t$), and the recent history, corresponding to times $s$ later than $\lfloor \alpha t \rfloor$. Here, $\alpha \in [0,1)$ is a fixed fraction that determines the proportion of history that is "forgotten." For the single-window IIS estimator (19), we can neglect the early history, and recast the estimator as follows

$$Z_k^t \equiv \frac{1}{t - \lfloor \alpha t \rfloor} \sum_{s=\lfloor \alpha t \rfloor+1}^{t} \frac{e^{-H_k(\mathcal{X}^s)}}{\sum_{\ell \in [K]} \gamma_\ell \, e^{-H_\ell(\mathcal{X}^s)} / Z_\ell^{s-1}}. \tag{74}$$

Note that, when $\alpha = 0$, (74) is exactly the IIS estimator (19). The influence of the early samples $\mathcal{X}^s$ on $Z_k^t$ is now indirect, and can decay faster than the law $t^{-1}$ otherwise allows. The problem with (74), though, is that it requires storing the entire state history in order to be able to evaluate the neglected part. We circumvent this shortcoming through the introduction of a system of *epochs*, which we define as follows:

Given a base multiplier $\phi > 1$, consider the exponentially increasing sequence of integers defined recursively by $\tau_0 = 0$, $\tau_1 = 1$, and $\tau_{l+1} = \lceil \phi \tau_l \rceil$ for $l \geq 2$. Here, $\lceil \phi \tau_l \rceil$ is the smallest integer greater than or equal to $\phi \tau_l$. For any integer $l \geq 1$, the semi-open interval $(\tau_{l-1}, \tau_l]$ is called the $l$-th epoch. For any real number $s \geq 0$, we let $n(s)$ be the smallest integer $l \geq 1$ such that $s \leq \tau_l$. The definition of $n(s)$ is such that $s$ lies in the epoch $n(s)$. The current time $t$ marks the beginning of a new epoch if and only if $t = \tau_{n(t)-1} + 1$, and the end of a new epoch if and only if $t = \tau_{n(t)}$. For any $\alpha \in [0,1)$, the epoch index $n(\alpha t)$ increases by 1 at every cycle $t$ for which $\lceil \alpha t \rceil$ crosses into a new epoch. Such crossing events happen if and only if simultaneously $\alpha(t-1) > 0$ and $n(\alpha t) > n(\alpha(t-1))$. We define the recent history of the simulation to be the sequence of state variables generated at the times $\tau_{n(\alpha t)-1} + 1$, $\tau_{n(\alpha t)-1} + 2, \ldots, t$.

Having established this notation, we can now modify the stochastic approximation triple for use with the history forgetting mechanism, and simultaneously incorporate the use of multiple parallel samplers. In order to use the parallel samplers, we need to define $\theta^t$ using information across all $R \geq 1$ replicas, indexed by $r \in [R] = \{1, \ldots, R\}$, which will allow us to define the rung density $\pi^t = \pi(\theta^t)$ across all windows. Each replica has its own state $(\mathcal{X}, \mathcal{K}, \mathcal{J})$ and runs in parallel with all other replicas; the only data shared between the replicas are the sets $\pi^t = \{\pi_{j;k}^t, j \in [J]; k: \lambda_k \in W_j\}$ and $F^t = \{F_{j;k}^t, j \in [J]; k: \lambda_k \in W_j\}$. We let $(\mathcal{X}^t(r), \mathcal{K}^t(r), \mathcal{J}^t(r))$ denote the tuple for the $r^{th}$ replica at time $t \geq 0$, which uses the transition kernel $\bar{P}_{\pi^{t-1}, F^{t-1}}$ (defined above in (52)) conditionally independently from all other replicas given $\pi^{t-1}, F^{t-1}$ for $t \geq 1$. A first step in defining the free energies with history forgetting is to compute and store the population counters

$$\mathcal{N}_j^{l,t} \equiv \sum_{s=\tau_{l-1}+1}^{\tau_l \wedge t} \sum_{r=1}^{R} 1_{\{j\}}(\mathcal{J}^s(r)), \tag{75}$$

representing the number of visits paid to window $W_j$ during the $l$-th epoch by the $R$ replicas. The number of visits paid to the window in the recent history of the simulation is

$$\mathcal{N}_j^t \equiv \sum_{l=n(\alpha t)}^{n(t)} \mathcal{N}_j^{l,t}. \tag{76}$$

Only the epochs $l$ with $n(\alpha t) \leq l \leq n(t)$ are needed at the current time. The current epoch is typically incomplete, so the population counters with $l = n(t)$ change with time. During any time $t$ that marks the beginning of a new epoch, new counters $\mathcal{N}_j^{n(t),t}$ are instantiated on the memory store. During any time $t$ for which $\alpha t$ crosses into a new epoch, the chronologically earliest counters $\mathcal{N}_j^{l,t}$, the ones with $l = n(\alpha t) - 1$, are removed from the memory store. The number of counters committed to the store thus remains bounded during the simulation.

The next step for defining the free energies with history forgetting is to define estimates within each epoch, which we call the epoch free energies. For convenience, we define the importance sampling ratios:

$$\mathcal{R}_{j;k}^t(r) \equiv \frac{e^{-H_k(X^t(r))}}{\sum_{\ell:\lambda_\ell \in W_j} \pi_{j;\ell}^{t-1} e^{F_{j;\ell}^{t-1} - H_\ell(X^t(r))}} 1_{\{j\}}(\mathcal{J}^t(r)), \tag{77}$$

with the denominator set to 1 if $\mathcal{N}_j^t = 0$. For all $n(\alpha t) \leq l \leq n(t)$, $k \in [K]$ and $j \in \text{win}(k)$, we then define the epoch free energies $F_{j;k}^{l,t}$ by

$$e^{-F_{j;k}^{l,t}} \equiv \frac{1}{\mathcal{N}_j^{l,t}} \sum_{s=\tau_{l-1}+1}^{\tau_l \wedge t} \sum_{r=1}^{R} \mathcal{R}_{j;k}^s(r), \tag{78}$$

where we use the convention that $F_{j;k}^{l,t} = 0$ if $\mathcal{N}_j^{l,t} = 0$. Despite the notation, at time $t$ all epoch free energies remain unchanged, except for those with $l = n(t)$, which can be computed as follows: for all $j \in \cup_{r=1}^{R} \{\mathcal{J}^t(r)\}$ and $k$ such that $\lambda_k \in W_j$

$$F_{j;k}^{l,t} = F_{j;k}^{l,t-1} - \ln\left(1 + \frac{1}{\mathcal{N}_j^{l,t}}\left[\mathcal{N}_j^{l,t-1} - \mathcal{N}_j^{l,t} + e^{F_{j;k}^{l,t-1}}\sum_{r=1}^{R}\mathcal{R}_{j;k}^{t}(r)\right]\right). \tag{79}$$

The overall free energy estimators with history forgetting $F_{j;k}^t$ are then obtained by straightforward averaging, as follows:

$$e^{-F_{j;k}^t} \equiv \frac{1}{\mathcal{N}_j^t}\sum_{l=n(\alpha t)}^{n(t)}\mathcal{N}_j^{l,t}e^{-F_{j;k}^{l,t}}. \tag{80}$$

Recursions for $\mu_{j;ki}^{l,t} = \mathbb{E}_k[\psi_i(\mathcal{X})]$ and $o_{j;k}^{l,t}$ analogous to (79) are constructed along similar lines and combined over epochs as in (80). Together, these three recursions define the tuple $\theta_{j;k}^{l,t} = \left(F_{j;k}^{l,t}, \mu_{j;k}^{l,t}, o_{j;k}^{l,t}\right)$ and the corresponding quantity averaged over epochs $\theta_{j;k}^t = \left(F_{j;k}^t, \mu_{j;k}^t, o_{j;k}^t\right)$ (see the Supplementary Materials, Part 2, Sections 5.1.1–5.1.4).

## 2.4 Putting it all together

In this section, we will use the illustration in Figure 3 to discuss the order in which the elements of TSS we have introduced are employed and highlight some of their computational characteristics. We will then provide practical guidance regarding the selection of windows and other parameters.

The simulation begins with each replica having some state $(x, k, j)$. As shown in the upper-left-hand rectangular boxes in Figure 3, each replica first swaps to the other window containing its current rung $k$, after which the MD dynamics are performed by each replica separately, and the Hamiltonians are evaluated for each window occupied by a replica. Although the evaluation of the $H_k(x)$ can be computationally expensive (especially if $\Lambda \subset \mathbb{R}^d$ and $d > 1$), the most computationally expensive part of these boxes is generally the MD step.

Moving on to the stochastic approximation updates (upper-right-hand box), we note that the epoch-based, per-window estimator updates to $F_{j;:}^{l,t-1}$ (as defined in (79)), and the analogous updates to $\mu_{j;:}^{l,t-1}$, $o_{j;:}^{l,t-1}$, and $\mathcal{N}_{j}^{l,t-1}$ (as defined in the SM, Sections 5.1.1–5.1.4), are computationally inexpensive operations relative to the cost of the MD step and the evaluation of the Hamiltonians for the active windows. The updates of $f_j^{t-1}$ and $\pi_{j;:}^{t-1}$ are more computationally expensive than the epoch-based, per-window estimator updates, because they require the minimization of (67) to obtain the window free energy offsets $\{f_j^{t-1}; j \in [J]\}$. Notably, the convex optimization problem (67) is a global operation, and its cost increases with the number of windows.

As shown in Figure 3, each element of the state $(x, k, j)$ is sampled sequentially, and neither the sampling of $x$ nor $j$ is dependent on the estimates $(F^{t-1}, \pi^{t-1})$. One motivation for structuring the algorithm in this way is to take advantage of this independence, by computing the estimates (as outlined in Sections 2.1–2.3) on an external device simultaneously with the MD step and the evaluation of the energies $\{H_k(x), k \in W_{j'}\}$, and then receiving the estimates in order to perform the $k$ dynamics, which is dependent on $(F^{t-1}, \pi^{t-1})$. The energies $\{H_k(x), k \in W_{j'}\}$ are sent to the external device to update the estimates for the next cycle.

On certain computing architectures, such as the special-purpose molecular dynamics supercomputer Anton 3,[32] the latency of communication with the external device (the time to send and receive a small packet of data) can be large relative to the MD step. To mitigate the cost of this communication, one can leverage the observation that self-adjustment reduces the variance of the free energy estimates, with the maximal improvement happening when the newly generated $(x', k')$ state is stochastically independent of previous $(x, k)$ states. We can approximate this stochastic independence by iterating a composition of the $x$- and $k$-transition kernels $\nu \geq 1$ times during a TSS cycle, without resampling the window index $j'$. The send and receive operations occur only once per cycle, and the $\nu$ iterations use the same fixed estimate $F^{t-1}, \pi^{t-1}$ received right before the first update to the state of $k$ within a cycle. The strategy is effective because the variables $(\mathcal{X}, \mathcal{K})$ are often strongly correlated, and would thus be poorly sampled if the $x$- and $k$-transition kernels in a cycle were not alternated before being composed.

Once the energy evaluations are available, the $k$ dynamics (the center-left box following the send operation in Figure 3) can be performed for each replica, an inexpensive process compared to the computation of the energies. The updates to $k$ do not require updates from the upper-right-hand box at every cycle, and thus do not incur any communication latency; they can be executed on the device that performs the sampling immediately after the energies are evaluated.

Finally, before proceeding to the next cycle $t + 1$, the epoch estimators must be updated. The epochs allow for the computation of error bars by means of a jackknife estimate. Although the computation of the error bars is computationally expensive (because for each epoch one needs to solve an optimization problem like the one for the reported free energies), this expense can be safely mitigated by computing the error bars infrequently, as they do not drive the simulation. (The full procedure is described in Part 2, Section 6 of the SM.)

We now provide some guidance on choosing parameters that have been introduced in TSS: the visit control parameter $\eta \geq 0$, the history forgetting parameter $\alpha \in [0,1)$ and base multiplier $\phi > 1$, the regularization parameters $\epsilon_\gamma$ and $\epsilon_\pi \in (0,1]$, as well as the selection of windows.

To begin, we recall that, at the extremes, the parameter $\eta \geq 0$ has the following binary behavior: When $\eta = 0$, visit control is "off," and otherwise it is "on." In practice, we find that $\eta = 2$ works well for a broad range of systems. A larger value of $\eta$ may further accelerate convergence for some systems (though we recall that the choice does not affect the asymptotic variance of the estimates), but with increased risk of rung moves that lead to failure of the MD component, whereas a value of $\eta$ smaller than 1 is never the optimal choice in our experience. Regardless of the particular value $\eta$ is set to, the most important choice is to turn visit control on (i.e., to set $\eta > 0$). In Section 3, a range of values of $\eta$ will be simulated, to illustrate the utility of $\eta > 0$ and the different rates at which the estimates reach their equilibrium value.

For the history forgetting parameter $\alpha$, which determines the proportion of history that is "forgotten," we recommend using the value 19%. We find this rate to be within a range of reasonable values for many systems, with the particular number chosen because the asymptotic statistical efficiency of using history forgetting relative to using the full trajectory, $\sqrt{1-\alpha} = 0.9$, is easy to remember. To determine $\phi$, it is convenient to think in terms of number of epochs $n_{\text{epochs}}$ and use the relation $\phi = \alpha^{-1/n_{\text{epochs}}}$; in general, we choose $n_{\text{epochs}} = 32$, which corresponds to $\phi \approx 1.0533$. Increasing the number of epochs will improve the precision of the error bars asymptotically, but will also require longer simulation times to ensure that the epochs are sufficiently decorrelated to use the jackknife error estimator.

The small regularization parameters $\epsilon_\gamma, \epsilon_\pi \in (0,1]$ ensure that the visitation frequency is at least a fraction $\epsilon_\gamma \epsilon_\pi$ of the frequency of the most visited rung. If, for instance, there is a window $W_j \subset \Lambda$ for which $H_k(x)$ varies little with $\lambda_k \in W_j$, then the unregularized density $\pi$ (which approaches $\gamma$ asymptotically) may assign very little weight to that window (as it should). The downside is that the window may not be visited sufficiently often in the recent history of the simulation, which can increase the errors of the overall free energy estimates, despite the window contributing little to the overall estimate. To accommodate such degenerate cases, we recommend setting $\epsilon_\gamma = 0.01$ and $\epsilon_\pi = 0.001$, which we have found to be suitable for a wide range of simulations.

In most cases, TSS is robust to the choice of window size (see, for instance, Section 3), and the user should observe good convergence of the free energy estimates. In cases where convergence is poor, the behavior of the rung sampling process may help indicate whether a poor choice of window size is the source of the problem. If the process is becoming diffusive (i.e., the rung process $\mathcal{K}^t$ is limited to moves between adjacent rungs, which can significantly increase the variance of the free energy estimates), or if the optimization problem becomes prohibitively expensive, the windows may be too small. One can reduce the computational expense resulting from a large number of windows, without increasing window size, by computing $f_i$ less frequently than on each cycle. On the other hand, the windows may be too large if one observes poor pre-asymptotic behavior of the estimates: notably, the exponential slowdown discussed in Section 1.4. Overly large windows also allow the selection of unphysical rungs early in the simulation, when the estimates are poor, which can lead to failures in the MD step. Reducing the size of the windows will eliminate these problems, with the added benefit of reducing the computational expense of evaluating the Hamiltonians.

## 3    Numerical example

Here we illustrate the use of TSS for sampling and estimation on an example problem relevant to computational drug discovery, and in particular how certain TSS parameter choices that influence sampling in the $\Lambda$ space affect free energy estimation in practice. Although the example we have chosen appears to be complex compared to typical models simulated with MCMC, it is a relatively straightforward instance of a problem that arises frequently in drug discovery:[33] the calculation of free energy differences in aqueous solution.

## 3.1 Description of the chemical system

We simulated a collection of eight amino acids (*A*, *C*, *F*, *H*, *I*, *L*, *M*, and *N*) that have a common substructure (*X*), which is alanine (*A*) with three hydrogens removed (see Figure 4). To compute the free energy differences among the different structures using TSS, we pick a parameterized family of distributions $\rho_\lambda, \lambda \in \Lambda$ in which the distribution of each amino acid structure corresponds to a rung $\lambda \in \Lambda$, and for which $\rho_\lambda$ varies smoothly with $\lambda$ between the structures. There are standard approaches (broadly called "alchemical free energy perturbation") in the MD literature for parameterizing such families of distributions, although we note that determining parameterizations that yield low-variance free energy estimates is itself an active area of research.[34,35]

Broadly, the alchemical free energy perturbation parameterization approach we use is as follows: To compute the free energy difference between the X and A substructures, for instance, a one-dimensional parameter $\lambda \in [0,1]$ is chosen such that the interaction energies of the three hydrogen atoms are "off" (equal to zero) when $\lambda = 0$ and "on" when $\lambda = 1$. As $\lambda$ gradually moves from 0 to 1, terms in the Hamiltonian $H(x, \lambda)$ corresponding to bonded interactions (such as bond stretching or angle bending) and non-bonded interactions (such as the van der Waals and electrostatic interactions) are progressively scaled up and become equal to the strength they normally have in the A structure once $\lambda = 1$. We note that the usual van der Waals interaction potential has a singularity at the origin (i.e., when distance between two atoms is zero), and we instead use a "softcore" van der Waals potential, which avoids this singularity (as described in Part 3, Section 7.2 of the SM). Additional information on the schedule and simulation details are provided in Part 3, Section 7.1 of the SM.

Each of the eight substructures is linked to the common $X$ substructure through such a parameterization, and we take $\Lambda$ to be a graph $\mathcal{G} = (\mathcal{V}, \mathcal{E})$ consisting of nine vertices $\mathcal{V} = \{X, A, C, F, H, I, L, M, N\}$ and eight edges $\mathcal{E} = \{(X, A), \ldots, (X, N)\}$. As shorthand, we refer to edge $(X, C)$ simply as "edge $C$".

We now discretize the $\Lambda$ space. We use a collection of 400 rungs, with 50 consecutive rungs for each edge $E_i \in \mathcal{E}$, $i = 0, \ldots, 7$. We identify each $E_i$ with the set $\{\lambda_{50i+k} \in \Lambda : k = 0, \ldots, 49\}$, which corresponds to the parameter set that guides the morph according to the family of distributions $\rho_\lambda(x)dx$, $\lambda \in E_i$; the morph $X \to A$ for $E_1$, for instance, is guided by the addition of the three missing hydrogens. The free energy difference of edge $E_i$ is defined as $F_{50i+49} - F_{50i}$. An illustration of the discretization of $\Lambda$ can be found in Figure 5.

## 3.2  *Numerical simulations and results*

In the next three sections, we study the effect on the free energy estimators of three features of TSS that affect the sampling in the $\Lambda$ space: the windowing system, the visit-control mechanism, and the self-adjustment property (the former introduced in Section 2.1, and the latter two in Section 1.3). All the free energy estimates are reported in units of $k_B T$. We performed each simulation using 32 replicas (see Section 2.3 for discussion of replicas), with a per-sampler time of 250 ns (aggregate 8 μs). The simulation length is much longer than needed to obtain useful estimates (as we will see below, for the best-performing setup the estimates stay within $0.5\ k_B T$ of their final values after just 1 ns of simulation). Long simulations do, however, enable a thorough examination of the effects of modifying the parameters, and also ensure that we obtain good estimates of the errors. The errors are estimated using the jackknife estimator (Section 2.3 and Part 2, Section 6 of the SM) with 32 epochs and a history forgetting parameter $\alpha = 0.19$.

We used the regularization parameters $\epsilon_\gamma = 0.01$ and $\epsilon_\pi = 0.001$ and, except for the alternatives tested in Section 3.2.2, a default value of $\eta = 2$ for the visit control parameter. Similarly, we used a default windowing system illustrated in Figure 5, except for the simulations in which we tested alternative systems (Section 3.2.1). We used an MD time step of 2 fs, with each TSS update being performed at a default interval of 9.6 ps (unless otherwise specified). The transition kernel $\mathcal{T}_k(x'|x)$ was determined by the force field (Amber ff99SB*-ILDN[36] with the TIP3P water model) and the ensemble (NPT; constant number of particles, pressure, and temperature) using the Desmond/GPU MD software;[37] further simulation details can be found in Part 3, Section 7.1 of the SM.

### 3.2.1 Windowing system

In this section we investigate the effect of different windowing systems on free energy estimation. As explained in Section 2.1, we expect that large windows will allow the process $\mathcal{K}^t$ to sample the rungs more effectively than will small windows, but that very large windows can lead to instabilities.

We simulated the system with three different window arrangements, labeled (A.1) through (A.3):

- (A.1): a single, very large window covering all 400 rungs
- (A.2): 392 windows of size 2, 8 windows of size 1 (each containing the last rung of an edge), and a single window of size 8 (containing the first rung of each of the eight edges)
- (A.3): our default setup: 58 windows, 48 of which have size 10, 8 of which have size 5, and the last 2 of which have sizes 160 (= 8 × 20) and 120 (= 8 × 15) that encompass the first 20 and 15 rungs, respectively, of every edge (see Figure 5)

First, simulation (A.1) failed after just a few moves of $\mathcal{K}^t$. The underlying cause was the instability of the free energy estimates, which led to $\mathcal{K}^t$ visiting a rung that is physically infeasible given $\mathcal{X}^{t-1}$, with atoms crashing into each other with unphysically large energies, ultimately leading to failure of the molecular dynamics integration scheme. We mention this result because it highlights the necessity of the windowing mechanism in practice.

At the other extreme, (A.2) used windows containing only two rungs. In this case there was little risk of visiting a physically infeasible rung, but the motion of $\mathcal{K}^t(r)$ for each sampler $r = 1,\ldots,32$ was severely restricted, as seen in Figure 6. Such behavior is unsurprising, because windows restrict the movement of the rung process $\mathcal{K}^t(r)$ to adjacent rungs, which promotes diffusive behavior of the rung process. Nevertheless, the simulation (A.2) successfully completed, and reported reasonable estimates of the free energy differences for each edge, with the caveat that the error bars for (A.2) tended to be larger for most edges than for (A.3), as shown in Figure 7. This comparison shows that the TSS free energy estimators are generally robust to the choice of windows, but that their selection affects the statistical efficiency of the simulation. As recommended in Section 2.4, the rung process can in practice be monitored, and if poor mixing properties are observed, the window size can be increased.

### 3.2.2 Visit control

In this section we explore the role of the visit control parameter $\eta \geq 0$ introduced in Section 1.3. Like the windowing system, visit control affects the dynamics of the rung process: Whenever $\eta > 0$, there is a bias toward under-visited rungs, which helps to prevent the rung process from getting trapped in a region of $\Lambda$. Unlike the windowing system, however, the choice of $\eta \geq 0$ does not affect the asymptotic variance of the estimates; it only modulates the speed at which the estimates approach a neighborhood of their equilibrium value. We present results from simulations with a range of $\eta$ values:

- (B.$\eta$): a visit-control parameter of $\eta$ ($\eta = 0, 1, 2, 4, 16,$ and $64$), with our default windowing system.

(Note that (B.2), because it uses our default value $\eta = 2$, is equivalent to (A.3); when we refer to (B.2) in the discussion, it refers to the same simulation.) We first compare simulation (B.0) to (B.2). In Figure 8, we see that the $\eta = 0$ simulation provides a significantly worse estimate of the free energies for each of the edges: The $\eta = 0$ trajectories all suffer from the exponential slowdown discussed in Section 1.4, and their error bars are noisier and an order of magnitude larger. The cause (and consequence) was the poor sampling of $\Lambda$: Without a mechanism to force them to explore the space $\Lambda$ more broadly, each of the 32 TSS samplers was stuck in a narrow region of rungs for a significant period early in the simulation. A comparison of the rung processes $\mathcal{K}^t(r)$ is shown in Figure 9 for two of the 32 replicas. This example illustrates how the lack of visit control can make it infeasible to quickly approach convergence.

Next we compare (B.2) to (B.$\eta$) for $\eta \neq 0, 2$. As indicated in Section 1.3, if $\eta$ is too large, the variance of $o_k/\gamma_k^\eta$ is large, and its fluctuations dominate the visit control rung distribution (33), creating an essentially random bias in the rung distribution. With $\eta = 64$, for instance, we see in Figure S1.a that, early in the simulation, the estimates fluctuated strongly compared to $\eta = 2$. If $\eta$ is too small, the bias towards under-visited rungs is not strong enough, and in Figure S1.a we see that the exponential slowdown is present to some extent for $\eta = 1$. The same simulations are shown on longer timescales in Figure S1.b, where the estimates for $\eta = 1$ and $\eta = 64$ are visibly worse than those for $\eta = 2$ for approximately the first third of the simulation; neither $\eta = 1$ nor $\eta = 64$ was worse than $\eta = 0$ (Figure 8), however. It is harder to decide between $\eta = 2, 4,$ and 16 based on visual inspection (see Figure S2.a and S2.b), but based on experience of other test systems, we generally find that a broad range of $\eta$ works well. The choice of $\eta$ affects only how quickly we reach a neighborhood of the asymptotic value, but does not affect the asymptotic variance; see Table 1 for the error estimates for each edge and for all simulated values of $\eta$.

### 3.2.3 Self-adjustment

The variance of any free energy estimator decreases as one decreases the interval at which samples are provided to the estimator, with diminishing returns as the interval shrinks to 0, owing to correlation between the samples. Because TSS performs both the sampling and estimation, it is able to reduce the variance of the free energy estimates in ways that are not possible when the estimator is decoupled from the method by which the samples are generated. We showed in Section 1.4 that when independent samples are available from $\rho_\lambda(x)dx$ for each $\lambda \in \Lambda$, the sampling of the rung process $\mathcal{K}^t$ strongly influenced the variance of the free energy estimates, a property called self-adjustment. We take the same approach with the current example. Let $\bar{P}_{\pi,Z}^{x,n}$ denote the transition kernel acting on $\mathcal{S} \times [K] \times [J]$ that applies $n$ steps of the MD transition kernel $\mathcal{T}_k(x'|x)dx$, similarly to $P_{\gamma,Z}^x$ in equation (5), but augmented with the windowing variable $j$, as described below equation (48) in Section 2.1. Instead of using, at time $t \geq 0$, $\bar{P}_{\pi^t,Z^t} = \bar{P}_{\pi^t,Z^t}^j \circ \bar{P}_{\pi^t,Z^t}^k \circ \bar{P}_{\pi^t,Z^t}^{x,n}$, which involves a single application of $\bar{P}_{\pi^t,Z^t}^k$ ($\nu = 1$) and $n$ applications of $\mathcal{T}_k$ per estimator update (as was the case for all the (A) and (B) simulations), we use the "high-frequency" kernel $\bar{P}_{\pi^t,Z^t}^\nu$, defined for $\nu > 1$ and any $\pi \in \mathcal{P}^\varepsilon(K)$ and $Z \in \mathbb{R}_{>0}^K$ by

$$\bar{P}_{\pi,Z}^\nu = \bar{P}_{\pi,Z}^j \circ \bar{P}_{\pi,Z}^k \circ \bar{P}_{\pi,Z}^{x,n/\nu} \circ \ldots \circ \bar{P}_{\pi,Z}^k \circ \bar{P}_{\pi,Z}^{x,n/\nu} = \bar{P}_{\pi,Z}^j \circ \left(\bar{P}_{\pi,Z}^k \circ \bar{P}_{\pi,Z}^{x,n/\nu}\right)^{(\nu)} \qquad (81)$$

Here it is assumed that $\nu$ divides $n$. In other words, rather than performing a single long MD trajectory followed by one sampling event of the rung $k$, we break the MD trajectory into $\nu$ slices and follow each slice with a sampling event of the rung $k$. The total number of MD steps between updates of the estimates $\pi^t, Z^t$ stays constant ($\nu$ applications of $n/\nu$ steps), but there are $\nu$ times as many intermediate updates to $\mathcal{K}^t$: $\mathcal{K}^{t,1}, \ldots, \mathcal{K}^{t,\nu}$ with $\mathcal{K}^{t+1} = \mathcal{K}^{t,\nu}$. All the intermediate moves at time $t \geq 0$ are performed within the same window using the estimates $\pi^t, Z^t$, and are generally less computationally demanding than the MD steps.

To study the effect of more rapid rung moves, as well as more rapid rung moves in combination with more rapid free energy updates, we performed simulations using the following setups:

- (C.1): we set $\nu = 25$ and the interval for free energy estimate updates to 9.6 ps; with our MD time step of 2 fs, this results in free energy updates every 4800 MD steps and rung moves every 192 steps.
- (C.2): we set $\nu = 1$ and the interval for free energy updates to 9.6 ps / 25 = 0.384 ps; this results in free energy updates and rung moves every 192 MD steps.

We compare (A.3), in which both free energy updates and rung moves occurred at intervals of 9.6 ps (4800 MD steps), to (C.1) and (C.2); these three simulations differ only in the intervals between free energy updates and rung moves. In Figure 10 we show the free energy estimates of simulations (A.3), (C.1), (C.2) along with their error bars. The simulations (C.1) and (C.2) have visibly lower variance than (A.3), with (C.2) having the lowest variance of all three (see Figure 10 and Table 1); indeed, (C.2) is statistically the best setup out of all our simulations. Notably, (C.1) has up to 1.8x lower estimated error than (A.3) depending on the edge, and (C.2) has up to 3.2x lower estimated error than (A.3) depending on the edge (see Table 1). These two comparisons suggest that a noteworthy portion of the error reduction from (A.3) to (C.2) is due to improved sampling of $\mathcal{K}^t$.

The significant variance reduction due to the improved rung sampling makes it relevant to consider the computational cost of an estimator update relative to the cost of a rung move. Estimator updates require information for all windows from all replicas, and thus incur a global communication cost that cannot be circumvented, on the top of the cost of solving the optimization problems described in 2.1 and 2.2. Rung moves, by contrast, require no communication between samplers, and can re-use estimates (all the moves in (62) use the same $\gamma, Z$), which is particularly important on architectures constrained by communication latencies. The effect of self-adjustment can be so dramatic that one can construct examples for which the asymptotic variance of the TSS free energy estimator is arbitrarily small relative to the asymptotic variance of a traditional estimator such as MBAR, for $\nu$ large enough (see Part 3, Section 8.4 of the SM for an example).

### 3.2.4 Concluding remarks on the numerical simulations

Each of the three mechanisms investigated in this section—windowing, visit control, and self-adjustment—affects the sampling of the space $\Lambda$ in a different way. The common theme is that, through their effect on $\mathcal{K}^t$, they can improve (or worsen) the free energy estimates. Our best setup (C.2) performs both sampling and estimation steps frequently, and yields useful estimates after just 1 ns of sampling (see Figure S3). We expect that within the framework provided by TSS there are other, as yet unknown parameters that could be estimated on the fly and used to modulate the sampling of $\mathcal{K}^t$ through the choice of rung distribution $\pi$. We also expect that one can find transition kernels with better mixing properties than the simulated tempering kernel $P_{\pi,Z}^k$. Finally, we expect that efficient implementations of the rung move will lead to sampling schemes in which $\nu > 1$ would be preferred.

**Acknowledgments**

The authors thank Justin Gullingsrud for help preparing the chemical system, Michael Bergdorf for insight into algorithm performance, Michael Eastwood for insightful discussions and a critical reading of the manuscript, and Berkman Frank for editorial assistance.

**Figures and Table**

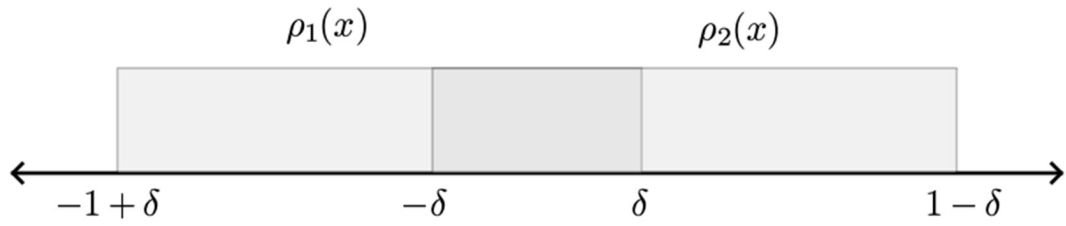

**Figure 1.** Densities of two uniform distributions, one on $[-1 + \delta, \delta]$, and the other on $[-\delta, 1 - \delta]$.

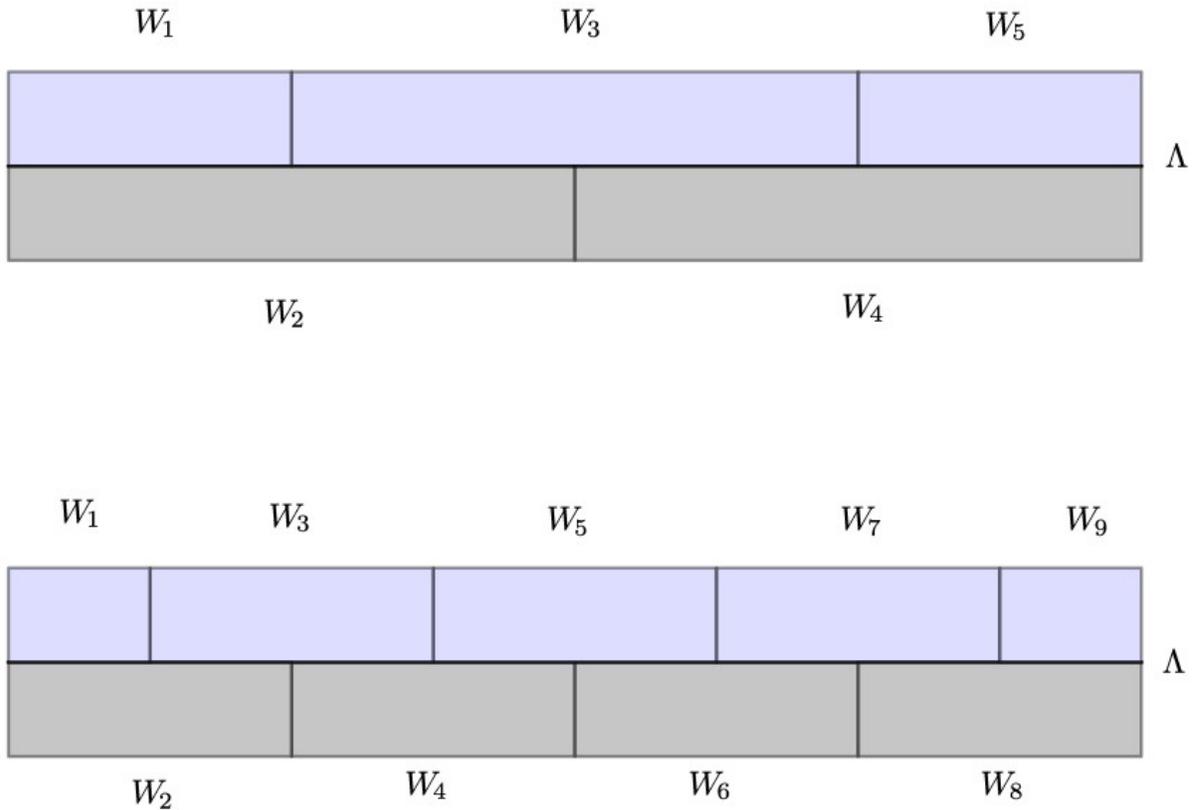

**Figure 2.** (Top): The eight edges, illustrated with the vertex being the common substructure X. The two large central windows occupy the first 25 and 30 rungs of each edge, containing 200 and 240 rungs, respectively. Each edge contains 60 rungs. The smaller central window thus contains rungs with indices $60i + k, i = 0, \ldots, 7, k = 0, \ldots, 24$. Aside from the two central windows, the remainder of the window structure is identical for each edge. (Bottom): Window structure for one of the eight edges. The windows are offset to satisfy the irreducibility requirements of Section 2.1. All non-central windows contain either 5 or 10 rungs.

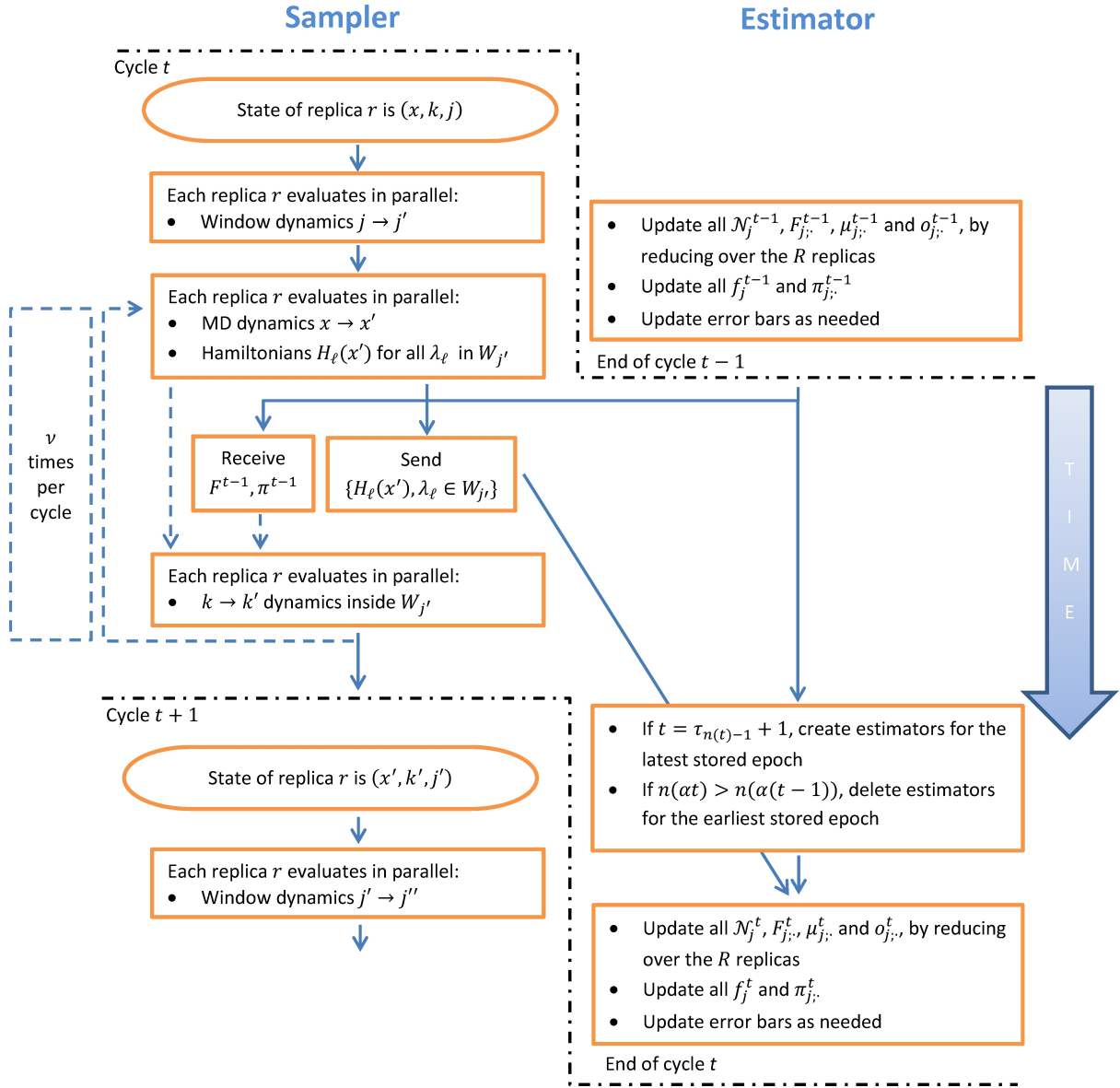

**Figure 3.** A full TSS cycle. The sampling moves (processes in the left column) involve energy and force evaluations, and the estimator updates (processes in the right column) perform the per-window stochastic approximation updates and produce global free energy estimates. The sampling and estimation can be executed on separate hardware devices, and both occur during a TSS cycle, the start and end of which is delineated between black dash-dotted lines. When the result of an operation is transmitted exactly once per TSS cycle, this is shown with a solid blue line; when the result may be transmitted more frequently, this is shown with a dashed blue line. The sampler device receives updated estimates $(F^{t-1}, \pi^{t-1})$ from the estimator device once per cycle, and these received estimates are used for the first time in a cycle only after the completion of the first MD dynamics and Hamiltonian evaluations. These Hamiltonian evaluations (evaluated at the end of this first run of MD dynamics) are sent to the estimator device for further updates. The MD dynamics, Hamiltonian evaluations, and rung move ($k \to k'$ dynamics) are iterated an additional $\nu - 1$ times within a cycle, with the same estimate $(F^{t-1}, \pi^{t-1})$ used for the rung moves in all $\nu - 1$ additional steps; these additional Hamiltonian evaluations are not sent to the estimator device. At the end of a cycle, the $x$ and $k$ dynamics have been updated in an interleaving fashion a total of $\nu$ times. While the additional $\nu - 1$ iterations of $(x, k)$ dynamics are performed for each replica by the sampler device, the estimator device updates the epoch counts, receives the Hamiltonian evaluations, uses them to compute the updated global estimates, and sends the updated estimates to the sampler device in time for the first $k' \to k''$ dynamics of the following cycle.

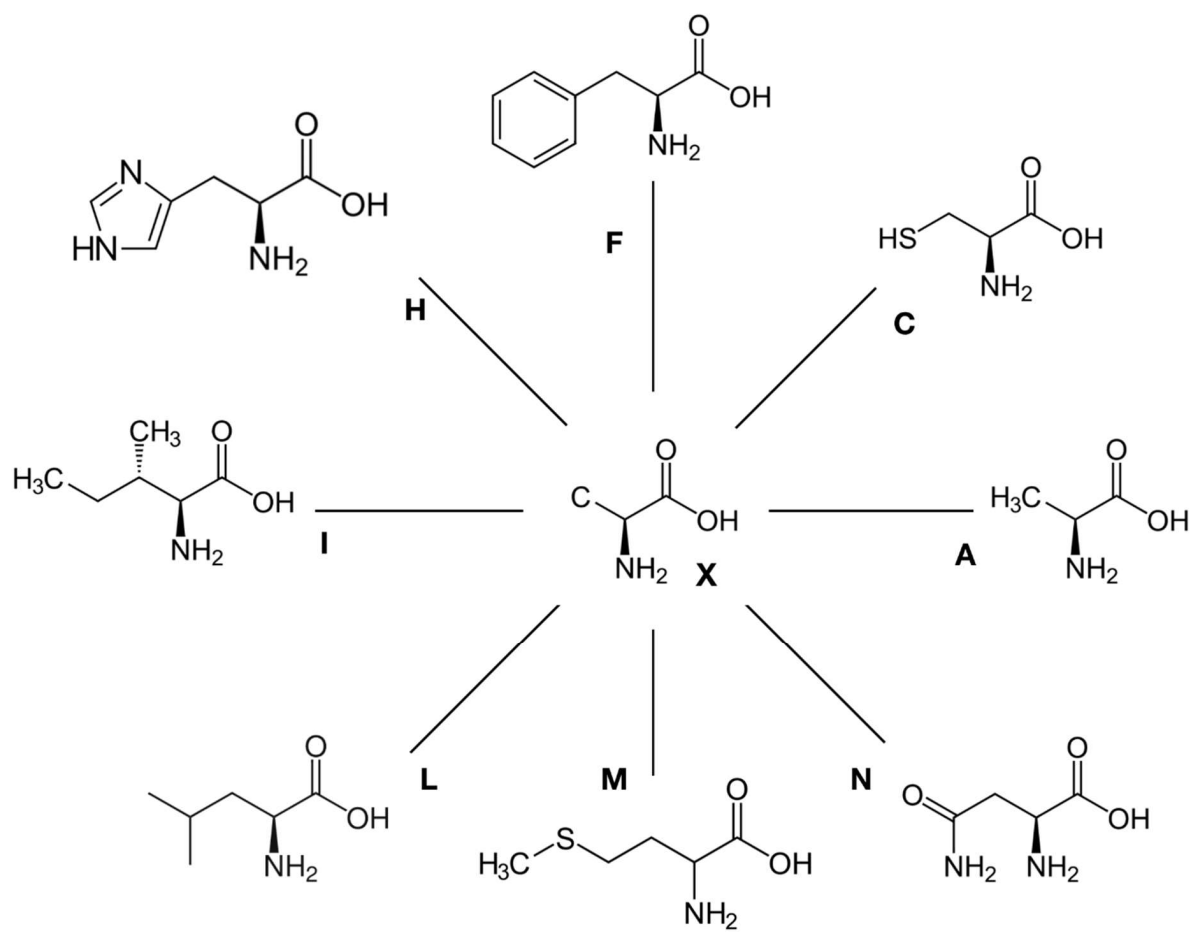

**Figure 4.** Collection of amino acids used in the simulations reported in Section 3. The central substructure, X, is common to the eight amino acid structures.

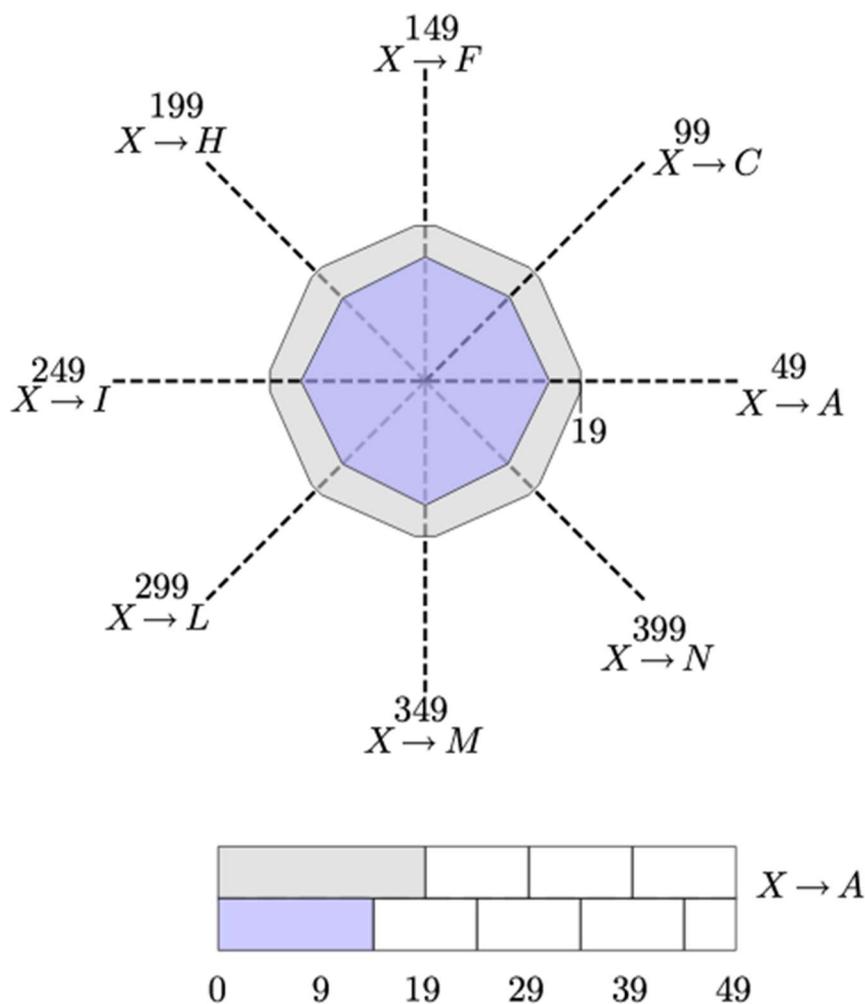

**Figure 5.** Default windowing system used in the amino acid simulations. (Top): The eight edges, illustrated with the vertex being the common substructure X. The two large central windows occupy the first 15 and 20 rungs of each edge, containing 120 and 160 rungs, respectively. Each edge contains 50 rungs. The smaller central window thus contains rungs with indices $50i + k, \ i = 0, \ldots, 7, \ k = 0, \ldots, 14$. Aside from the two central windows, the remainder of the window structure is identical for each edge. (Bottom): Window structure for one of the eight edges. The windows are offset to satisfy the irreducibility requirements of Section 2.1. All non-central windows contain either 5 or 10 rungs.

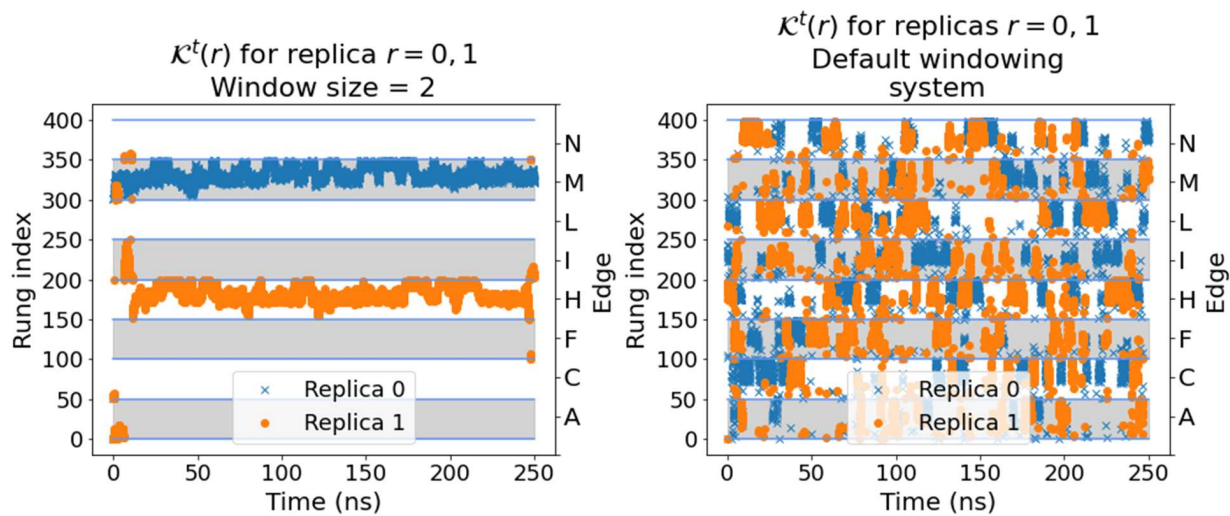

**Figure 6.** Small windows can cause diffusive rung dynamics. Plotted are simulations using windows of size 2 (A.2) and the default setup (A.3). As discussed in Section 2.4, the rung process can be monitored throughout the simulation, and the left-hand figure indicates that the window size of 2 is too small. (Left): The full rung trajectory for two of the replicas (with indices 0 and 1). The left $y$-axis labels indicate the rung index, and the right $y$-axis labels show which edge each range of rung indices corresponds to. The bands delineate the eight different edges. The rung trajectories show that the two samplers are confined to narrow regions of $\Lambda$, resulting in inefficient sampling of $\Lambda$; the behavior of the other replicas (not shown) is similar. (Right): The full rung trajectory for two replicas with our default windowing system. Each of the two replicas visits each edge within 40 ns, which indicates much more effective sampling than with windows of size 2.

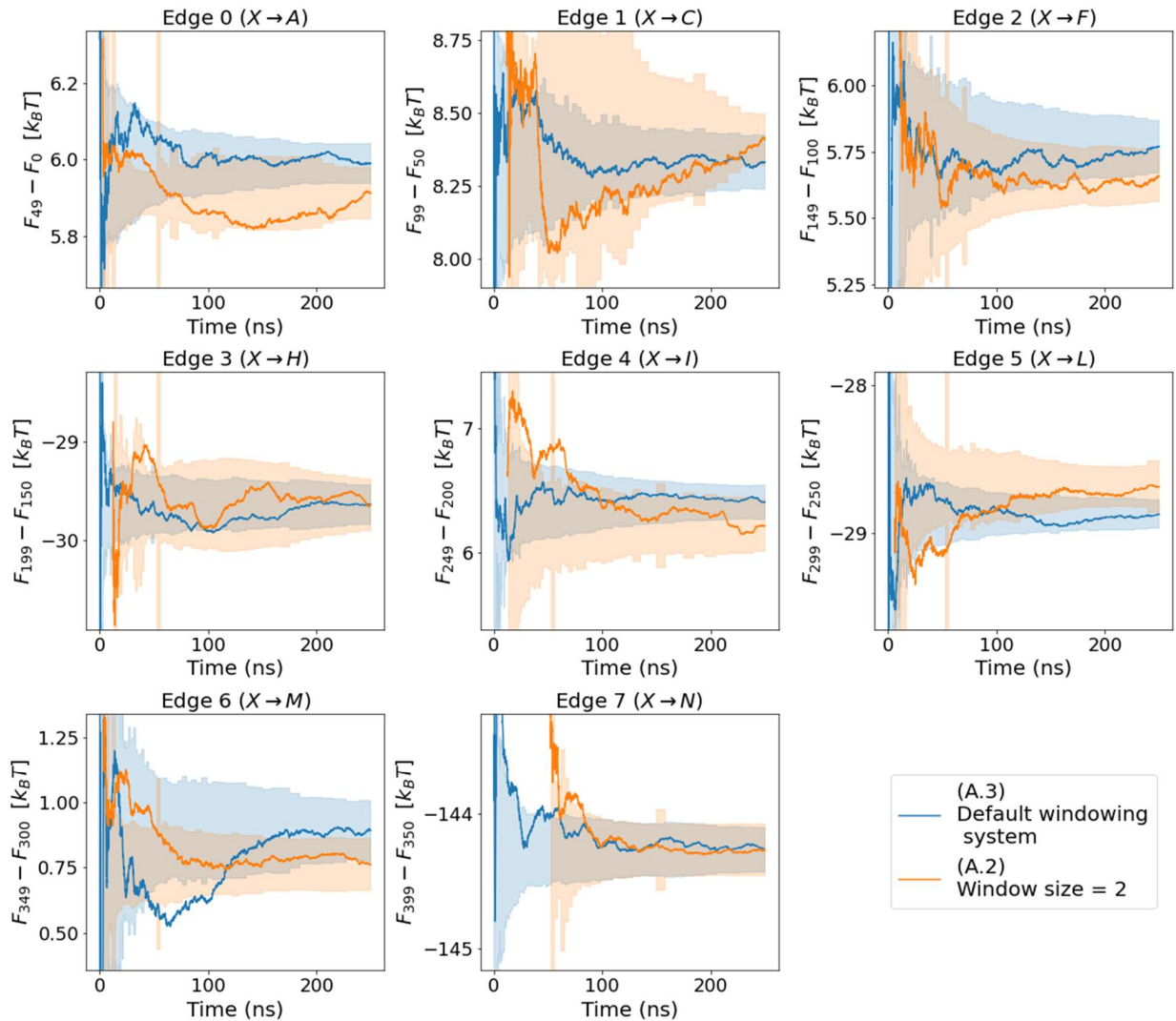

**Figure 7.** Window size influences the error bars of the free energy estimates. Plotted are simulations using windows of size 2 (A.2) and the default setup (A.3). The free energy difference estimates are shown for each of the eight edges along with their error bars. The error bars are larger and noisier for the simulation with windows of size 2 than when the default setup is used. It also takes longer to produce estimates for each edge with the size 2 windows than with our default windowing system; the delay is particularly notable in Edge 7, where it takes almost 50 ns before an estimate is available. This is due to the diffusive sampling of $\Lambda$, as explained in Figure 6; replicas do not sample the rungs associated to Edge 7 sufficiently to produce an estimate until after 50 ns.

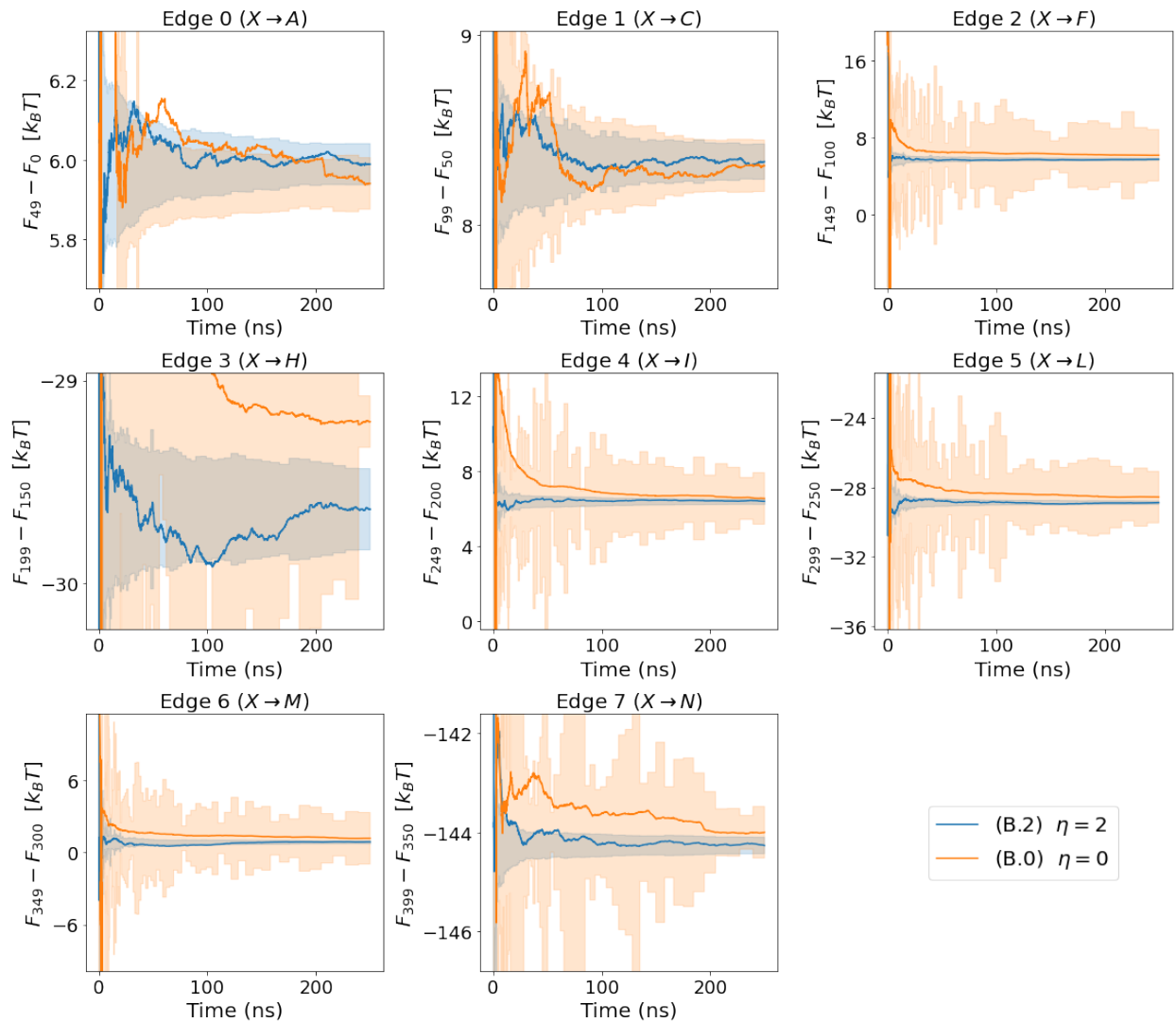

**Figure 8.** Visit control accelerates convergence of the estimates to a neighborhood of their equilibrium values. Simulations using $\eta = 0$ (B.0) and $\eta = 2$ (B.2) are shown. The free energy difference estimates are shown for each of the eight edges, along with their error bars. Use of visit control ($\eta = 2$) removes the exponential slowdown of the estimates apparent in its absence ($\eta = 0$), wherein the curves have a steady, small negative slope (see especially edges 2 through 7). Edges 0 and 1 rapidly converge to a neighborhood of the equilibrium even without visit control, because the substructure $X$ is very similar to both $A$ and $L$, making the estimation easier.

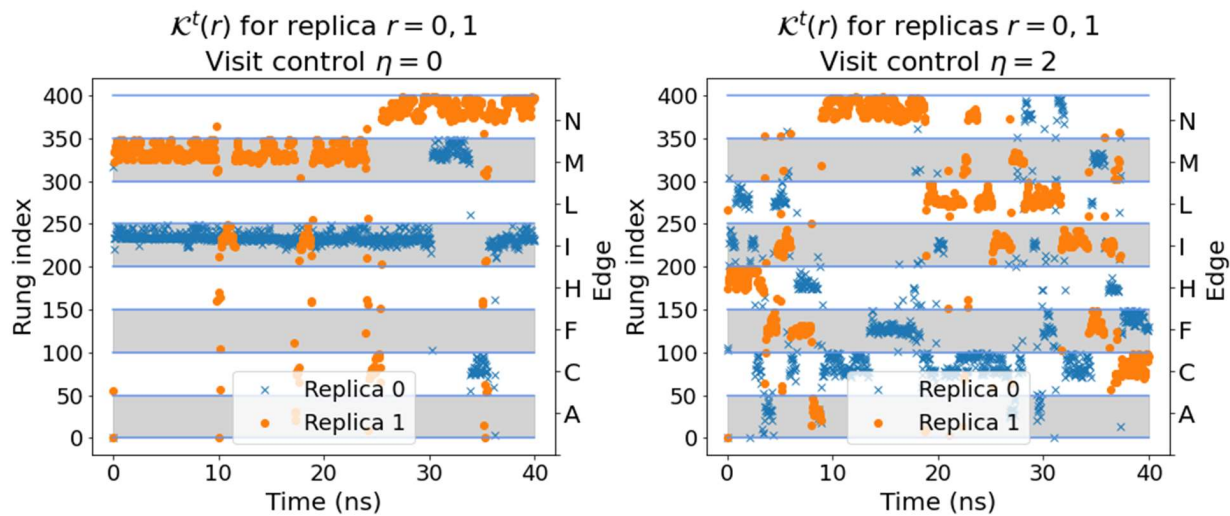

**Figure 9.** The visit control mechanism accelerates convergence in part by promoting efficient sampling of the $\Lambda$ space. Plotted are simulations using $\eta = 0$ (B.0) and $\eta = 2$ (B.2). The rung trajectories for two of the replicas (with indices 0 and 1), which are representative of the behavior of the remaining 30 replicas, are shown for the first 40 ns of simulation. Early on, the rung process in the $\eta = 0$ simulation struggled to explore $\Lambda$. The $\eta = 2$ simulation, on the other hand, immediately started exploring across all eight edges, because the visit control mechanism pushes the replicas toward under-visited rungs. Past the 50 ns mark (not shown), the replicas for $\eta = 0$ began effectively exploring the rung space, and the visitation patterns became indistinguishable between $\eta = 0$ and $\eta = 2$ as the simulation progressed. Efficient exploration of the rung space is a necessary (but insufficient) condition for efficient estimation of the free energy differences (see Figure 8).

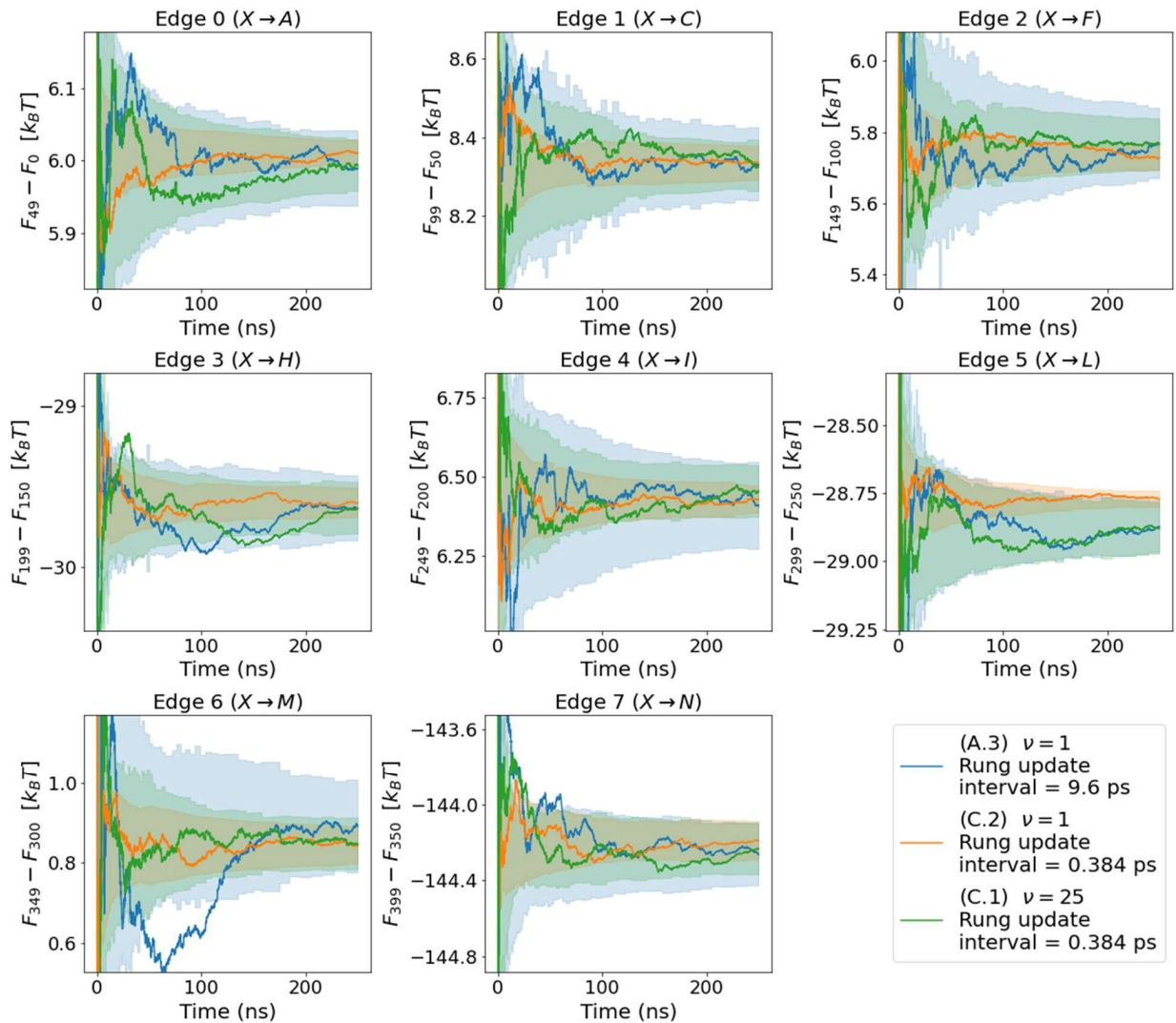

**Figure 10.** Increasing the frequency of rung updates can reduce the variance of the estimates, even if the number of estimator updates remains constant. The figure shows simulations that all use visit control with $\eta = 2$. Two simulations have $\nu = 1$ (i.e., one rung update per estimator update), one with the default interval between estimator updates of 9.6 ps (A.3) and one with a shorter interval of 0.384 ps (C.2); a third simulation (C.1) also has a short rung update interval of 0.384 ps but with $\nu = 25$, so estimator updates occur at the default interval. (C.2), which had high frequency rung and estimator updates performed the best out of all the simulations in terms of the estimated variance. Notably, (C.1) had smaller error bars than (A.3) for most of the edges, despite having the same estimator update frequency and only differing in rung update frequency. (Numerical values for the error bars are given in Table 1).

|  |  | Edge | | | | | | | |
|---|---|---|---|---|---|---|---|---|---|
|  |  | A | C | F | H | I | L | M | N |
| Simulation | (A.2) | 5.911 (0.034) | 8.409 (0.044) | 5.658 (0.048) | -29.638 (0.130) | 6.220 (0.104) | -28.686 (0.087) | 0.761 (0.049) | -144.265 (0.096) |
|  | (A.3) | 5.989 (0.026) | 8.331 (0.046) | 5.768 (0.049) | -29.636 (0.100) | 6.407 (0.069) | -28.873 (0.048) | 0.890 (0.058) | -143.66 (0.082) |
|  | (B.0) | 5.940 (0.032) | 8.310 (0.068) | 6.194 (1.336) | -29.203 (0.064) | 6.560 (0.686) | -28.552 (0.743) | 1.188 (1.072) | -143.999 (0.259) |
|  | (B.1) | 6.006 (0.029) | 8.247 (0.051) | 5.698 (0.061) | -29.663 (0.089) | 6.305 (0.103) | -28.765 (0.057) | 0.751 (0.063) | -144.291 (0.062) |
|  | (B.4) | 6.016 (0.030) | 8.380 (0.049) | 5.702 (0.046) | -29.558 (0.077) | 6.478 (0.087) | -28.756 (0.054) | 0.943 (0.050) | -144.250 (0.077) |
|  | (B.16) | 5.969 (0.030) | 8.324 (0.061) | 5.874 (0.053) | -29.513 (0.087) | 6.370 (0.061) | -28.738 (0.058) | 0.788 (0.058) | -144.115 (0.101) |
|  | (B.64) | 6.036 (0.030) | 8.391 (0.055) | 5.730 (0.074) | -29.472 (0.121) | 6.493 (0.071) | -28.755 (0.069) | 0.827 (0.055) | -144.169 (0.103) |
|  | (C.1) | 5.994 (0.018) | 8.327 (0.033) | 5.765 (0.036) | -29.636 (0.080) | 6.453 (0.041) | -28.881 (0.047) | 0.847 (0.032) | -144.235 (0.068) |
|  | (C.2) | 6.010 (0.009) | 8.334 (0.022) | 5.728 (0.018) | -29.602 (0.047) | 6.426 (0.022) | -28.773 (0.015) | 0.843 (0.023) | -144.191 (0.051) |

**Table 1.** Free energy difference estimate and root-mean-square error (in parentheses) at 250 ns for each edge, for each simulation. The expression used to compute the mean square error is provided in Section 6 of the SM.

# Times Square Sampling: Supplementary Materials

December 9, 2021

# Contents

# Part I
# Theory

## 1 Optimality for stochastic approximation

### 1.1 Changing variables from $(Z, \xi)$ to $(F, \mu)$

We recall the notation from the main text: the variables $(\mathcal{X}, \mathcal{K}) \in \mathcal{S} \times [K]$ are the pair of state and rung variables; the parameter $\theta = (F, \mu, o)$ is estimated throughout the simulation and determines the conditional distribution of $\mathcal{K}$ given $\mathcal{X}$. The parameter $\theta$ is equivalently used in the main text as $\theta = (Z, \xi, o)$, and the two are related by the changes of variables $Z = e^{-F}$ and $\xi = e^{-F}\mu$, or more explicitly as

$$Z_k = e^{-F_k}, \quad \xi_{km} = e^{-F_k}\mu_{km}, \qquad k \in [K], m \in [M]. \tag{1.1}$$

We have the original recursions

$$Z_k^{t+1} = Z_k^t + \frac{1}{t+1}\left(\frac{e^{-H_k(\mathcal{X}^{t+1})}}{\sum_{\ell\in[K]}\pi_\ell(\theta^t)e^{-H_\ell(\mathcal{X}^{t+1})}/Z_\ell^t} - Z_k^t\right) \qquad (1.2)$$

$$\xi_{km}^{t+1} = \xi_{km}^t + \frac{1}{t+1}\left(\frac{e^{-H_k(\mathcal{X}^{t+1})}}{\sum_{\ell\in[K]}\pi_\ell(\theta^t)e^{-H_\ell(\mathcal{X}^{t+1})}/Z_\ell^t}\psi_m(\mathcal{X}^{t+1}) - \xi_{km}^t\right) \qquad (1.3)$$

$$o_k^{t+1} = o_k^t + \frac{1}{t+1}\left(\frac{1}{\gamma_k(\theta^t)}\mathbf{1}_{\{k\}}(\mathcal{K}^{t+1}) - o_k^t\right) \qquad (1.4)$$

with mean field

$$\dot{Z}_k = \frac{Z_k^\star}{\sum_{\ell\in[K]}\pi_\ell(\theta(t))Z_\ell^\star/Z_\ell(t)} - Z_k(t) \qquad (1.5)$$

$$\dot{\xi}_{km} = \frac{Z_k^\star}{\sum_{\ell\in[K]}\pi_\ell(\theta(t))Z_\ell^\star/Z_\ell(t)}\mathbb{E}_k[\psi_m(\mathcal{X})] - \xi_{km}(t) \qquad (1.6)$$

$$\dot{o}_k = \frac{\pi_k(\theta(t))}{\gamma_k(\mu(t))}\frac{Z_k^\star/Z_k(t)}{\sum_{\ell\in[K]}\pi_\ell(\theta(t))Z_\ell^\star/Z_\ell(t)} - o_k(t) \qquad (1.7)$$

where $\theta(t) = (Z(t),\mu(t),o(t))$. In the variables $(F,\mu)$, the first two recursions are

$$F_k^{t+1} = F_k^t - \log\left(1 + \frac{1}{t+1}\left[\frac{e^{F_k^t - H_k(\mathcal{X}^{t+1})}}{\sum_{\ell\in[K]}\pi_\ell e^{F_\ell^t - H_\ell(\mathcal{X}^{t+1})}} - 1\right]\right) \qquad (1.8)$$

$$\mu_{km}^{t+1} = e^{F_k^{t+1} - F_k^t}\left(\mu_{km}^t + \frac{1}{t+1}\left[\frac{e^{F_k^t - H_k(\mathcal{X}^{t+1})}}{\sum_{\ell\in[K]}\pi_\ell e^{F_\ell^t - H_\ell(\mathcal{X}^{t+1})}}\psi_m(\mathcal{X}^{t+1}) - \mu_{km}^t\right]\right). \qquad (1.9)$$

Observe that

$$\mu_{km}^{t+1} = e^{F_k^{t+1} - F_k^t}\left(\mu_{km}^t + \frac{1}{t+1}\left[\frac{e^{F_k^t - H_k(\mathcal{X}^{t+1})}}{\sum_{\ell\in[K]}\pi_\ell e^{F_\ell^t - H_\ell(\mathcal{X}^{t+1})}}\psi_m(\mathcal{X}^{t+1}) - \mu_{km}^t\right]\right) \qquad (1.10)$$

$$= \frac{\left(\mu_{km}^t + \frac{1}{t+1}\left[\frac{e^{F_k^t - H_k(\mathcal{X}^{t+1})}}{\sum_{\ell\in[K]}\pi_\ell e^{F_\ell^t - H_\ell(\mathcal{X}^{t+1})}}\psi_m(\mathcal{X}^{t+1}) - \mu_{km}^t\right]\right)}{\left(1 + \frac{1}{t+1}\left[\frac{e^{F_k^t - H_k(\mathcal{X}^{t+1})}}{\sum_{\ell\in[K]}\pi_\ell e^{F_\ell^t - H_\ell(\mathcal{X}^{t+1})}} - 1\right]\right)} \qquad (1.11)$$

$$= \mu_{km}^t + \frac{1}{t+1}\frac{e^{F_k^t - H_k(\mathcal{X}^{t+1})}}{\sum_{\ell\in[K]}\pi_\ell e^{F_\ell^t - H_\ell(\mathcal{X}^{t+1})}}[\psi_m(\mathcal{X}^{t+1}) - \mu_{km}^t] + O(t^{-2}). \qquad (1.12)$$

In the variables $(F,\mu)$, we therefore have the mean fields

$$\dot{F}_k = 1 - \frac{e^{F_k - F_k^\star}}{\sum_{\ell\in[K]}\pi_\ell(\theta(t))e^{F_\ell - F_\ell^\star}} \qquad (1.13)$$

$$\dot{\mu}_{km} = \frac{e^{F_k - F_k^\star}}{\sum_{\ell\in[K]}\pi_\ell(\theta(t))e^{F_\ell - F_\ell^\star}}(\mu_{km}^\star - \mu_{km}). \qquad (1.14)$$

Using the change of variables (1.1) and the chain rule, we find that the two pairs of equations (1.13), (1.14) and (1.5), (1.6) are equivalent (although the corresponding discrete recursions are not equivalent). We will now identify the inverse of the Jacobian in the $(F, \mu, o)$ coordinates.

## 1.2 Identifying the inverse Jacobian

As discussed in Section 1.2 of the main text, in a stochastic approximation procedure of the form

$$\theta^{t+1} = \theta^t + \frac{1}{t+1}\Gamma G(Y^{t+1}; \theta^t), \qquad (1.15)$$

the matrix $\Gamma$ is called the gain matrix, and the noisy observable $G$ has an associated mean field $g(\theta) = \mathbb{E}_\theta[G(Y; \theta)]$. The matrix that minimizes the asymptotic variance of $\theta^t$ (when $\theta^t$ converges) is $\Gamma = -J_g(\theta^\star)^{-1}$, where $J_g(\theta)$ is the Jacobian of $g$ evaluated at $\theta$ (see Section 3.2.3, Proposition 4 of [3]).

Consequently, in order to demonstrate the optimality of the stochastic approximation procedure for $(F^t, \mu^t)$, it suffices to identify the inverse of the Jacobian of the expected gain $(g^F(\theta), g^\mu(\theta))$ in the main text. Tan accomplished this for the free energy estimates in [22], Theorem 2. It turns out this property can be extended to the case of estimating the averages

$$\mu_{km}^\star = \mathbb{E}_k[\psi_m(X)] = \int_S \psi_m(x)\rho_k(x)dx. \qquad (1.16)$$

The key to obtaining the optimal gain matrix is to use the assumption that $\sum_{\ell \in [K]} \pi_\ell(F, \mu, o) = 1$ for all $(F, \mu, o) \in \mathbb{R}^K \times \mathbb{R}^{KM} \times \mathcal{P}^\varepsilon(K)$, which allows us to avoid computing the derivatives of $\pi$ with respect to $F, \mu, o$ at the fixed point $F^\star, \mu^\star$; the Jacobian will turn out to be independent of $o \in \mathbb{R}_{\geq 0}^K$.

**Proof** [Of Proposition 1]

The optimal gain matrix for a stochastic approximation procedure is given by the inverse of the Jacobian evaluated at the unique equilibrium point. We will show that the Jacobian for the mean field $(\dot{F}, \dot{\mu}) = (g^F(F, \mu, o), g^\mu(F, \mu, o))$ is equal to the identity at the equilibrium $(F^\star, \mu^\star)$. We then show the same property extends to $\dot{o} = g^o(F, \mu, o)$ when

$$\pi_k(F, \mu, o) = \frac{\gamma_k(\mu)e^{F_k^\star - F_k}}{\sum_{\ell \in [K]} \gamma_\ell(\mu)e^{F_\ell^\star - F_\ell}} \quad (\equiv \pi_k^\circ(F, \mu, o)) \qquad (1.17)$$

The stochastic approximation procedure for estimating the averages $\mu_{km}^\star$ is defined by the system of equations

$$F_k^{t+1} = F_k^t + \frac{1}{t+1}(1 - R_k(X^{t+1}; F^t, \mu^t, o^t)) \qquad (1.18)$$

$$\mu_{km}^{t+1} = \mu^t + \frac{1}{t+1}R_k(X^{t+1}; F^t, \mu^t, o^t)(\psi_m(X^{t+1}) - \mu^t), \qquad (1.19)$$

where

$$R_k(x; F, \mu, o) = \frac{e^{F_k - H_k(x)}}{\sum_{\ell=1}^K \pi_\ell(F, \mu, o)e^{F_\ell - H_\ell(x)}}. \qquad (1.20)$$

The mean fields are

$$g_k^F(F, \mu, o) = 1 - \frac{e^{F_k - F_k^\star}}{\sum_{\ell \in [K]} \pi_\ell(F, \mu, o)e^{F_\ell - F_\ell^\star}} \qquad (1.21)$$

$$g_{km}^\mu(F, \mu, o) = \frac{e^{F_k - F_k^\star}}{\sum_{\ell \in [K]} \pi_\ell(F, \mu, o)e^{F_\ell - F_\ell^\star}}(\mu_{km}^\star - \mu_{km}). \qquad (1.22)$$

Let $g = (g_1^F, \ldots, g_K^F, g_{11}^\mu, \ldots, g_{KM}^\mu) \in \mathbb{R}^{K+KM}$, and let

$$r_k(F, \mu, o) = \frac{e^{F_k - F_k^\star}}{\sum_{\ell \in [K]} \pi_\ell(F, \mu, o) e^{F_\ell - F_\ell^\star}}. \qquad (1.23)$$

To save space, we omit writing the explicit dependence of $\pi_k$ on $\theta = (F, \mu, o)$. Note that

$$\frac{\partial r_k}{\partial F_j}(F, \mu, o) = \delta_{kj} \frac{e^{F_k - F_k^\star}}{\sum_{\ell \in [K]} \pi_\ell e^{F_\ell - F_\ell^\star}} \qquad (1.24)$$

$$- \frac{e^{F_k - F_k^\star}}{\left(\sum_{\ell \in [K]} \pi_\ell e^{F_\ell - F_\ell^\star}\right)^2} \left( \pi_j e^{F_j - F_j^\star} + \sum_{\ell \in [K]} \frac{\partial \pi_\ell}{\partial F_j} e^{F_\ell - F_\ell^\star} \right)$$

$$\frac{\partial r_k}{\partial \mu_{jn}}(F, \mu, o) = - \frac{e^{F_k - F_k^\star}}{\left(\sum_{\ell \in [K]} \pi_\ell e^{F_\ell - F_\ell^\star}\right)^2} \left( \sum_{\ell \in [K]} \frac{\partial \pi_\ell}{\partial \mu_{jn}} e^{F_\ell - F_\ell^\star} \right). \qquad (1.25)$$

We can write the derivatives of the mean field in terms of the partial derivatives of $r_k$, this time omitting writting explicitly the dependence of $r_k$ and its derivatives on $F, \mu$:

$$\frac{\partial g_k^F}{\partial F_j}(F, \mu, o) = -\frac{\partial r_k}{\partial F_j} \qquad (1.26)$$

$$\frac{\partial g_{km}^\mu}{\partial F_j}(F, \mu, o) = \frac{\partial r_k}{\partial F_j}(\mu_{km}^\star - \mu_{km}) \qquad (1.27)$$

$$\frac{\partial g_k^F}{\partial \mu_{jn}}(F, \mu) = -\frac{\partial r_k}{\partial \mu_{jn}} \qquad (1.28)$$

$$\frac{\partial g_{km}^\mu}{\partial \mu_{jn}}(F, \mu, o) = \frac{\partial r_k}{\partial \mu_{jn}}(\mu_{km}^\star - \mu_{km}) - \delta_{kj}\delta_{mn} r_k. \qquad (1.29)$$

To determine the optimal gain matrix, we evaluate the Jacobian at the point

$$(F^\star, \mu^\star, o) = (F_1^\star, \ldots, F_K^\star, \mu_{11}^\star, \ldots, \mu_{KM}^\star, o_1, \ldots, o_K),$$

where $o \in \mathbb{R}_{\geq 0}^K$ is arbitrary. To this end, we first evaluate the Jacobian at $(F^\star, \mu, o)$ for some arbitrary $(\mu, o) \in \mathbb{R}^{KM} \times \mathbb{R}_{\geq 0}^K$, and start by noting that

$$\frac{\partial r_k}{\partial \mu_{jn}}(F^\star, \mu) = -\sum_{\ell \in [K]} \frac{\partial \pi_\ell}{\partial \mu_{jn}}(F^\star, \mu) = -\frac{\partial}{\partial \mu_{jn}} \sum_{\ell \in [K]} \pi_\ell(F^\star, \mu) = 0, \qquad (1.30)$$

$$\frac{\partial r_k}{\partial F_j}(F^\star, \mu) = \delta_{kj} - \left( \pi_j + \sum_{\ell \in [K]} \frac{\partial \pi_\ell}{\partial F_j}(F^\star, \mu) \right) = \delta_{kj} - \pi_j, \qquad (1.31)$$

since $\sum_{\ell \in [K]} \pi_\ell(F,\mu,o) = 1$ for all $(F,\mu,o) \in \mathbb{R}^K \times \mathbb{R}^{KM} \times \mathbb{Z}_+^K$. It follows that

$$\frac{\partial g_k^F}{\partial F_j}(F^\star, \mu^\star, o) = -(\delta_{kj} - \pi_j(F^\star, \mu^\star, o)) \tag{1.32}$$

$$\frac{\partial g_{km}^\mu}{\partial F_j}(F^\star, \mu^\star, o) = 0 \tag{1.33}$$

$$\frac{\partial g_k^F}{\partial \mu_{jm}}(F^\star, \mu^\star, o) = 0 \tag{1.34}$$

$$\frac{\partial g_{km}^\mu}{\partial \mu_{jn}}(F^\star, \mu^\star, o) = -\delta_{kj}\delta_{mn}. \tag{1.35}$$

Using the notation $D_F g$ and $D_\mu g$ to denote the Jacobian of a function $g$ with respect to $F$ and $\mu$ respectively, we can write the Jacobian $J$ in block matrix form as

$$J(F^\star, \mu^\star, o) = \begin{pmatrix} D_F g^F & D_\mu g^F \\ D_F g^\mu & D_\mu g^\mu \end{pmatrix} \tag{1.36}$$

$$= \begin{pmatrix} -1+\pi_1 & \pi_2 & \cdots & \pi_K & 0 & \\ \pi_1 & -1+\pi_2 & \cdots & \pi_K & & \ddots \\ \vdots & \vdots & \ddots & \pi_K & & 0 \\ \pi_1 & \cdots & & -1+\pi_K & 0 & 0 \\ 0 & 0 & 0 & & -1 & 0 \\ & \ddots & & & & \ddots \\ 0 & 0 & 0 & 0 & & -1 \end{pmatrix} \tag{1.37}$$

$$= \begin{pmatrix} \mathbf{1}\pi^\top - I_{K \times K} & 0 \\ 0 & -I_{KM \times KM} \end{pmatrix}. \tag{1.38}$$

In the above evaluation of the Jacobian, the dependence on $o$ appears only through $\pi$ in the matrix $\mathbf{1}\pi^\top - I_{K \times K}$. This matrix admits a one-dimensional left null space spanned by the row eigenvector $\pi^\top = (\pi_1, \ldots, \pi_K)$ with eigenvalue 0. Thus $(\mathbf{1}\pi^\top - I_{K \times K})\mathbf{x}$ lies in the orthogonal complement of $\pi$ for all $\mathbf{x} \in \mathbb{R}^K$. The matrix $\mathbf{1}\pi^\top - I_{K \times K}$ is invertible though if its action is restricted to this complement, which is possible because 1.20 implies $\sum_k \pi_k(1 - R_k) = 0$. The inverse is obtained by simply ignoring the term $\mathbf{1}\pi^\top$ in 1.38. We derive

$$\Gamma^\star = -Dg(F^\star, \mu^\star)^{-1} = \begin{pmatrix} I_{K \times K} & 0_{K \times MK} \\ 0_{MK \times K} & I_{MK \times MK} \end{pmatrix}, \tag{1.39}$$

immediately leading to the recursions (1.18) and (1.19).

The second part of Proposition 1 asserts that if $\pi(\theta) = \pi^\circ(\theta)$, the recursion for $o^t$ also has minimal asymptotic variance. With $\pi^\circ$ defined by (1.17),

$$\pi_k^\circ(F, \mu) = \gamma_k(\mu) \frac{e^{F_k^\star - F_k}}{\sum_{\ell \in [K]} \gamma_\ell(\mu) e^{F_\ell^\star - F_\ell}}. \tag{1.40}$$

From equation (1.7), the mean field for $o$ is

$$\dot{o}_k(t) = \frac{\pi_k(\theta(t))}{\gamma_k(\mu(t))} \frac{e^{F_k(t) - F_k^\star}}{\sum_{\ell \in [K]} \pi_\ell(\theta(t)) e^{F_\ell(t) - F_\ell^\star}} - o_k(t), \tag{1.41}$$

in $(F,\mu)$ coordinates. Inserting $\pi(\theta(t)) = \pi^\circ(\theta(t))$ into (1.41), we find

$$\dot{o}_k(t) = 1 - o_k(t), \tag{1.42}$$

for which the Jacobian is trivially given by $-I_{K \times K}$. It follows that, when $\pi = \pi^\circ$, the recursion (1.4) is optimal. □

## 2 Visit control in the steady state approximation

In this section, we demonstrate the effect of visit control on the mean field of the free energy estimates when a steady-state approximation is used. In other words, we will study $g_k^F(F,\mu,o)$ (defined in (1.13)) under the assumption that $o = o^\star(F,\mu)$, which is equivalent to saying that $o = (o_1, \ldots, o_K)$ satisfies $g_k^o(F,\mu,o) = 0$ for all $k \in [K]$. In Section 2.1 below we explicitly write out the calculations to identify $g_k^F(F,\mu,o^\star(F,\mu))$, these calculations being used in Section 1.4 of the main text. Then, in Section 2.2 we study a general mean field and establish that the Lyapunov function for the system decreases faster in the steady state approximation.

### 2.1 Uniform distributions computations

Here, we explicitly work out the modified mean field equation for the free energy estimates, as presented in Section 1.4 of the main text. The model in question uses two Hamiltonians defined, for some $0 < \delta < 1/2$, as

$$H_1(x) = \infty \cdot 1_{[-1+\delta,\delta]^c}(x)$$
$$H_2(x) = \infty \cdot 1_{[-\delta,1-\delta]^c}(x),$$

where $A^c$ denotes the complement of a set $A \subset \mathbb{R}$. For the purpose of illustrating visit control, we let $\pi$ denote a general function of the tilts $o$ whose target range is contained in $\mathcal{P}^\varepsilon(K)$. We define $\Delta^t = F_2^t - F_1^t$ where $F_k^t$ is defined using the linearized estimator presented in Section 1.2 of the main text,

$$F_k^{t+1} = F_k^t + \frac{1}{t+1}\left(1 - \frac{e^{F_k^t - H_k(X^{t+1})}}{\sum_{\ell \in [K]} \pi_\ell e^{F_\ell^t - H_\ell(X^{t+1})}}\right), \tag{2.1}$$

which gives

$$\Delta^{t+1} = \Delta^t + \frac{1}{t+1}\left(\frac{e^{F_1^t - H_1(X^{t+1})}}{\sum_{\ell \in [K]} \pi_\ell e^{F_\ell^t - H_\ell(X^{t+1})}} - \frac{e^{F_2^t - H_2(X^{t+1})}}{\sum_{\ell \in [K]} \pi_\ell e^{F_\ell^t - H_\ell(X^{t+1})}}\right) \tag{2.2}$$

$$= \Delta^t + \frac{1}{t+1} \cdot \frac{e^{F_1^t - H_1(X^{t+1})} - e^{F_2^t - H_2(X^{t+1})}}{\pi_1 e^{F_1 - H_1(X^{t+1})} + \pi_2 e^{F_2 - H_2(X^{t+1})}} \tag{2.3}$$

$$= \Delta^t + \frac{1}{t+1} \cdot \frac{e^{-H_1(X^{t+1})} - e^{\Delta^t - H_2(X^{t+1})}}{\pi_1 e^{-H_1(X^{t+1})} + \pi_2 e^{\Delta^t - H_2(X^{t+1})}} \tag{2.4}$$

$$= \Delta^t + \frac{1}{t+1} G(X^{t+1};\Delta) \tag{2.5}$$

The mean field for $\Delta$ can be identified either by computing $\mathbb{E}[G(\mathcal{X}; \Delta)]$, or simply by using the expression for the mean field $g_k(F)$ (equation 11 in the main text),

$$g_k(F) = 1 - \frac{e^{F_k - F_k^\star}}{\sum_{\ell \in [K]} \pi_\ell e^{F_\ell - F_\ell^\star}}. \tag{2.6}$$

We then have, with $\Delta(t) = F_2(t) - F_1(t)$ (i.e., $F(t) = (0, \Delta(t))$ up to an additive shift by any multiple of the vector $(1,1)$) and $F^\star = 0$,

$$\dot{\Delta}(t) = g_2(F(t)) - g_1(F(t)) \tag{2.7}$$

$$= \frac{e^{F_1(t)} - e^{F_2(t)}}{\pi_1 e^{F_1(t)} + \pi_2 e^{F_2(t)}} \tag{2.8}$$

$$= \frac{1 - e^{\Delta(t)}}{\pi_1 + \pi_2 e^{\Delta(t)}}. \tag{2.9}$$

Now we use TSS's choice of $\pi$ for this example, which is defined for any $\eta > 0$ by

$$\pi_k^{\text{TSS}}(o) = \frac{\gamma_k o_k^{-\eta}}{\sum_{\ell \in [K]} \gamma_\ell o_\ell^{-\eta}}. \tag{2.10}$$

The dependence of $\gamma = (\gamma_1, \gamma_2)$ on $\mu$ has been ignored in the above expression, since in the example we take $\gamma = (1/2, 1/2)$. The steady state $o^\star$ as a function of $F$ is obtained by solving $g_k^o(F, \mu, o) = 0$, where

$$g_k^o(F, \mu, o) = \frac{1}{\gamma_k} \cdot \frac{\pi_k(o) e^{F_k - F_k^\star}}{\pi_1(o) e^{F_1 - F_1^\star} + \pi_2(o) e^{F_2 - F_2^\star}} - o_k. \tag{2.11}$$

(The dependence on $\mu$ of $g_k^o$ is present only to mimic the expressions presented in Section 1.3.) Using (2.10) and $F_k^\star = 0$ in (2.11), we find that $(o_k^\star)^{\eta+1} \propto e^{F_k}$, and so

$$\pi_k = \gamma_k (o_k^\star)^{-\eta} \propto e^{-\beta F_k}, \tag{2.12}$$

where $\beta = \eta/(\eta + 1)$. Using $\pi_1 + \pi_2 = 1$, we get

$$\pi(o^\star(F)) = \left(\gamma_1 o_1^{-\eta}, \gamma_2 o_2^{-\eta}\right) = \left(\frac{e^{-\beta F_1}}{e^{-\beta F_1} + e^{-\beta F_2}}, \frac{e^{-\beta F_1}}{e^{-\beta F_2} + e^{-\beta F_2}}\right) \tag{2.13}$$

$$= \left(\frac{1}{1 + e^{-\beta \Delta}}, \frac{e^{-\beta \Delta}}{1 + e^{-\beta \Delta}}\right) \tag{2.14}$$

$$= \left(\frac{e^{\beta \Delta/2}}{e^{\beta \Delta/2} + e^{-\beta \Delta/2}}, \frac{e^{-\beta \Delta/2}}{e^{\beta \Delta/2} + e^{-\beta \Delta/2}}\right) \tag{2.15}$$

$$= \frac{1}{2\cosh(\beta \Delta/2)} \left(e^{\beta \Delta/2}, e^{-\beta \Delta/2}\right) \tag{2.16}$$

Finally, using $o^\star$ in (2.7) yields

$$\frac{1-e^\Delta}{\pi_1+\pi_2 e^\Delta} = \frac{e^{-\Delta/2}-e^{\Delta/2}}{\pi_1 e^{-\Delta/2}+\pi_2 e^{\Delta/2}} \tag{2.17}$$

$$= -2\sinh(\Delta/2)\frac{1}{\pi_1 e^{-\Delta/2}+\pi_2 e^{\Delta/2}} \tag{2.18}$$

$$= -2\sinh(\Delta/2)\frac{2\cosh(\beta\Delta/2)}{e^{(\beta-1)\Delta/2}+e^{(1-\beta)\Delta/2}} \tag{2.19}$$

$$= -2\sinh(\Delta/2)\frac{\cosh(\beta\Delta/2)}{\cosh((1-\beta)\Delta/2)}, \tag{2.20}$$

which is the result reported in Section 1.4.

## 2.2 Improved decrescence of Lyapunov function

In Section 1.3 of the main text we introduced the visit control mechanism, which proposes the class of controls

$$\pi_k^{\text{TSS}}(\theta) = \frac{\gamma_k(Z_k/Z_k^\circ)^{\eta/(\eta+1)}}{\sum_{\ell\in[K]}\gamma_\ell(Z_\ell/Z_\ell^\circ)^{\eta/(\eta+1)}}. \tag{2.21}$$

for some paramater $\eta \geqslant 0$. The mean field for $(Z,o)$ is

$$\dot{Z}_k = \frac{Z_k^\star}{\sum_{\ell\in[K]}\pi_\ell^{\text{TSS}}(\theta(t))Z_\ell^\star/Z_\ell(t)} - Z_k(t) \tag{2.22}$$

$$\dot{o}_k = \frac{\pi_k^{\text{TSS}}(\theta(t))}{\gamma_k(\mu(t))}\frac{Z_k^\star/Z_k(t)}{\sum_{\ell\in[K]}\pi_\ell^{\text{TSS}}(\theta(t))Z_\ell^\star/Z_\ell(t)} - o_k(t). \tag{2.23}$$

The law 2.21 has been adopted because if the process $o(t)$ can be kept close to its stationary regime, then the rate at which $Z(t)$ approaches $Z^\star$ increases in magnitude monotonically with increasing $\eta \geq 0$. More precisely, we measure the rate of decay of the relative entropy (or Kullback-Leibler divergence)

$$V_\gamma(Z) = D_{\text{KL}}(\gamma || r) = -\sum_{k\in[K]}\gamma_k \log(r_k/\gamma_k), \tag{2.24}$$

where

$$r_k = \gamma_k x_k \bigg/ \left(\sum_{\ell\in[K]}\gamma_\ell x_\ell\right) \tag{2.25}$$

with $x_k = Z_k^\star/Z_k$. The relative entropy serves as a Lyapunov function for the dynamics of $Z$ (see Lemma 2). We will show that the rate of decay increases monotonically with increasing $\eta \geq 0$ if the process $o(t)$ is assumed to be stationary, thus proving Proposition 2.

To ease notation clutter, we will drop the dependence on $\mu$ in $\gamma_k(\mu)$ and write only $\gamma_k$, and we will also stop denoting explicit time dependence. At stationarity of the tilts $o$ (that is, $\dot{o}_k = 0$ for all $k \in [K]$), we have $Z_k^\circ = Z_k^\star$ for all $k \in [K]$ and

$$\pi_k^{\text{TSS}}(\theta) = \frac{\gamma_k(Z_k/Z_k^\star)^{\eta/(\eta+1)}}{\sum_{\ell\in[K]}\gamma_\ell(Z_\ell/Z_\ell^\star)^{\eta/(\eta+1)}} = \frac{\gamma_k x_k^{-\eta/(\eta+1)}}{\sum_{\ell\in[K]}\gamma_\ell x_\ell^{-\eta/(\eta+1)}}. \tag{2.26}$$

We compute

$$\dot{Z}_k = \frac{Z_k^\star}{\sum_{\ell \in [K]} \pi_\ell Z_\ell^\star / Z_\ell} - Z_k = \frac{\sum_\ell \gamma_\ell x_\ell^{-\eta/(\eta+1)}}{\sum_{\ell \in [K]} \gamma_\ell x_\ell^{1/(\eta+1)}} Z_k^\star - Z_k \tag{2.27}$$

and

$$\frac{d}{dt} V_\gamma(Z) = -\sum_{k \in [K]} \gamma_k \left( \frac{Z_k^\star / Z_k}{\sum_{\ell \in [K]} \gamma_\ell Z_\ell^\star / Z_\ell} - 1 \right) \frac{\dot{Z}_k}{Z_k} \tag{2.28}$$

$$= -\sum_{k \in [K]} \gamma_k \left( \frac{x_k}{\sum_{\ell \in [K]} \gamma_\ell x_\ell} - 1 \right) \left( \frac{\sum_\ell \gamma_\ell x_\ell^{-\eta/(\eta+1)}}{\sum_{\ell \in [K]} \gamma_\ell x_\ell^{1/(\eta+1)}} x_k - 1 \right). \tag{2.29}$$

If $\eta = 0$ (no visit control), then

$$\frac{d}{dt} V_\gamma(Z) = -\sum_{k \in [K]} \gamma_k \left( \frac{x_k}{\sum_{\ell \in [K]} \gamma_\ell x_\ell} - 1 \right)^2 \tag{2.30}$$

is negative. We wish to show that the rate of decay becomes progressively more negative as $\eta$ is increased to infinity. It suffices to show that

$$\frac{d}{d\eta} \frac{d}{dt} V_\gamma(Z) = \frac{1}{(\eta+1)^2} \left( \sum_{k \in [K]} \gamma_k x_k \left( \frac{x_k}{\sum_{\ell \in [K]} \gamma_\ell x_\ell} - 1 \right) \right) \tag{2.31}$$

$$\times \left( \frac{\sum_\ell \gamma_\ell x_\ell^{1/(\eta+1)} x_\ell^{-1} \log(x_\ell)}{\sum_{\ell \in [K]} \gamma_\ell x_\ell^{1/(\eta+1)}} - \frac{\sum_\ell \gamma_\ell x_\ell^{1/(\eta+1)} x_\ell^{-1}}{\sum_{\ell \in [K]} \gamma_\ell x_\ell^{1/(\eta+1)}} \frac{\sum_\ell \gamma_\ell x_\ell^{1/(\eta+1)} \log(x_\ell)}{\sum_{\ell \in [K]} \gamma_\ell x_\ell^{1/(\eta+1)}} \right) \tag{2.32}$$

is negative. The overall factor on the first row of the right-hand side is positive, as follows from the Cauchy-Schwarz inequality. The factor on the second row of the right-hand side is negative. This follows from the correlation inequality

$$\mathbb{E}[f(X)g(X)] \geqslant \mathbb{E}[f(X)]\mathbb{E}[g(X)] \tag{2.33}$$

applied to the strictly increasing functions $f(x) = -x^{-1}$ and $g(x) = \log(x)$, where $f, g : (0, \infty) \to \mathbb{R}$. The correlation inequality can be deduced from

$$\sum_{k \in [K]} \sum_{\ell \in [K]} w_k w_\ell (f(x_k) - f(x_\ell))(g(x_k) - g(x_\ell)) \geqslant 0, \tag{2.34}$$

which holds by monotonicity of $f$ and $g$, because either both $f(x_k) - f(x_\ell)$ and $g(x_k) - g(x_\ell)$ are positive or both are negative. Here $w_k = \gamma_k x_k^{1/(\eta+1)} / \left( \sum_{\ell \in [K]} \gamma_\ell x_\ell^{1/(\eta+1)} \right)$. Owing to $\gamma \in \mathcal{P}^\varepsilon(K)$, the correlation inequality is strict unless $x_k$ is a constant independent of $k \in [K]$.

## 3 Variance calculations for TSS and MBAR

Free energy difference estimators have been studied for at least the past 70 years. Following the introduction of the Zwanzig relation [26], which is a form of importance sampling, Bennett [2] developed an improved estimator known as the Bennett Acceptance Ratio (BAR), which is the two-rung version of MBAR.

These estimators are known by different names in different fields; [25] is perhaps the earliest mathematical treatment of the maximum likelihood estimator (which we call MBAR); other relevant works include [5, 21, 10, 19, 6, 20, 17], and undoubtedly many others.

The variance of any estimator for the free energy differences is related to the amount of overlap between the distributions. In the examples where we compute variance, the rung distribution $\pi$ is always taken to be fixed. For such fixed $\pi \in \mathcal{P}^\varepsilon(K)$, we define the overlap matrix $O$ by

$$O_{ij} = \mathbb{E}\left[\frac{e^{F_i^\star - H_i(X)} e^{F_j^\star - H_j(X)}}{\left(\sum_{k \in [K]} \pi_k e^{F_k^\star - H_k(X)}\right)^2}\right], \tag{3.1}$$

where $X$ is distributed according to the mixture distribution $\sum_{k \in [K]} \pi_k e^{F_k^\star - H_k(x)} dx$.

Throughout, we make the assumption that $O$ is irreducible. This is an assumption that goes back to Vardi [25]. The motivation for this assumption is that the differences $F_i - F_j$ can only be resolved if there is overlap between the densities $\rho_i(x)$ and $\rho_j(x)$, or alternatively, if for any $i, j \in [K]$ there is a sequence $(i_k)_{k=1}^n$ such that $i_1 = i$, $i_n = j$, and

$$\int_S \mathbf{1}\{\rho_{i_k}(x) > 0\} \mathbf{1}\{\rho_{i_{k+1}}(x) > 0\} dx > 0, \quad k = 1, \ldots, n-1. \tag{3.2}$$

See, for instance, [6] (bottom of p. 1071).

## 3.1 Variance of stochastic approximation

Consider the stochastic approximation recursion

$$\theta^{t+1} = \theta^t + \frac{1}{t+1} \Gamma G(X^{t+1}; \theta^t), \tag{3.3}$$

and let $g(\theta) = \mathbb{E}_\theta[G(X, \theta)]$ denote the mean field, where the expectation is with respect to the stationary distribution $p_\theta(x) dx$. When the function $G$ and the gain matrix $\Gamma$ satisfy the appropriate conditions to guarantee almost sure c convergence of the stochastic approximation algorithm, the asymptotic covariance of the iterates $\theta^t$ is the solution $\Sigma$ to the Lyapunov equation (see eqn. 3.2.7, Part I in Benveniste et al., 1990)

$$\left(\frac{I}{2} + \Gamma Dg(\theta^\star)\right) \Sigma + \Sigma \left(\frac{I}{2} + \Gamma Dg(\theta^\star)\right)^\top + \Gamma R(\theta^\star) \Gamma^\top = 0. \tag{3.4}$$

Here, $\theta^\star$ is the unique solution to $g(\theta) = 0$, $Dg(\theta)$ is the Jacobian of $g$ with respect to $\theta$, $I$ is the identity matrix and $R(\theta^\star)$ is the integrated autocovariance, defined by

$$R(\theta) = \sum_{n=-\infty}^{\infty} \text{Cov}_\theta(G(X_0; \theta), G(X_n; \theta)), \tag{3.5}$$

where $X_n$ is generated according to the Markov chain which has stationary distribution $p_\theta(x) dx$.

The optimal choice of gain matrix, in the sense of minimizing $\Sigma$ in the rank-ordering on positive semidefinite matrices, is the gain

$$\Gamma = -Dg(\theta^\star)^{-1}, \tag{3.6}$$

which leads to the asymptotic covariance matrix (eqn. 3.2.8 in Benveniste et al. 1990)

$$\Sigma = Dg(\theta^\star)^{-1} R(\theta^\star) Dg(\theta^\star)^{-T}. \tag{3.7}$$

The optimal gain matrix for the $F$ component was identified in (1.39) to be $\Gamma = -I_{K \times K}$, so the asymptotic covariance matrix of the TSS free energy estimates TSS is given by

$$\Sigma_{\text{TSS}} = R(\theta^\star), \tag{3.8}$$

where the observable function $G$ which determines the integrated covariance $R$ in (3.5) is the function $G = G^F$ from the main text, which we recall is defined componentwise by

$$G_k^F(x; \theta) = 1 - \frac{e^{F_k - H_k(x)}}{\sum_{\ell \in [K]} \pi_\ell(\theta) e^{F_\ell - H_\ell(x)}}. \tag{3.9}$$

Suppose $\pi(\theta)$ is fixed and equal to some $\pi^\star \in \mathcal{P}^\epsilon(K)$, and that the TSS samples $(\mathcal{X}^{t+1}, \mathcal{K}^{t+1})$ are drawn independently from the density $p_{F^t, \pi^\star}$, the expression for the asymptotic covariance matrix simplifies to

$$(\Sigma_{\text{TSS}})_{ij} = \text{Cov}_{\theta^\star}(G^{F_i^\star}(\mathcal{X}; \theta^\star), G^{F_j^\star}(\mathcal{X}; \theta^\star)) \tag{3.10}$$

$$= \mathbb{E}_{\theta^\star}\left[\left(1 - \frac{e^{F_i^\star - H_i(\mathcal{X})}}{\sum_{\ell \in [K]} \pi_\ell^\star e^{F_\ell^\star - H_\ell(\mathcal{X})}}\right)\left(1 - \frac{e^{F_j^\star - H_j(\mathcal{X})}}{\sum_{\ell \in [K]} \pi_\ell^\star e^{F_\ell^\star - H_\ell(\mathcal{X})}}\right)\right] \tag{3.11}$$

$$= -1 + \mathbb{E}_{\theta^\star}\left[\frac{e^{F_i^\star - H_i(\mathcal{X})}}{\sum_{\ell \in [K]} \pi_\ell^\star e^{F_\ell^\star - H_\ell(\mathcal{X})}} \frac{e^{F_j^\star - H_j(\mathcal{X})}}{\sum_{\ell \in [K]} \pi_\ell^\star e^{F_\ell^\star - H_\ell(\mathcal{X})}}\right]. \tag{3.12}$$

Written in matrix form, this is

$$\Sigma_{\text{TSS}} = O - \mathbf{1}\mathbf{1}^\top, \tag{3.13}$$

where $\mathbf{1} = (1, \ldots, 1) \in \mathbb{R}^K$. Since the vector $\mathbf{1}$ is orthogonal to any vector $v_{ij} = e_j - e_i$, where $e_i, e_j$ are the $i^{th}$ and $j^{th}$ unit basis vectors respectively, for the purposes of determining the variance of the free energy difference estimates we can take

$$\Sigma_{\text{TSS}} = O. \tag{3.14}$$

Thus, $O$ is the asymptotic covariance matrix for TSS when samples $(\mathcal{X}^{t+1}, \mathcal{K}^{t+1})$ are drawn independently from $p_{\theta^t}(x,k)dx\mathcal{C}(dk)$. The expression for the asymptotic variance of the maximum likelihood estimator is presented in the next Section (Section 3.2); the proof that $\Sigma_{\text{TSS}} = O$ is smaller than the MLE covariance matrix is in Section 3.4.

## 3.2 Variance of maximum likelihood estimation

Given $n$ independent samples $\mathcal{X}_1, \ldots, \mathcal{X}_n$, $n_k$ of which are from distribution $\rho_k(x)dx$, $n = n_1 + \cdots + n_K$, the maximum likelihood estimator for the free energy differences (which we call MBAR) is a solution $F$ of the system of equations

$$e^{-F_k} = \frac{1}{n} \sum_{i=1}^n \frac{e^{-H_k(\mathcal{X}_i)}}{\sum_{\ell=1}^K \frac{n_\ell}{n} e^{F_\ell - H_\ell(\mathcal{X}_i)}}. \tag{3.15}$$

The solution is unique up to a common additive constant that disappears when computing a free energy difference. The existence and uniqueness of solutions to this system of equations has been investigated in depth by Vardi [25]. The large sample variance has also been studied in [6], [5], and subtle aspects of the theory have been in studied in [10]. The MBAR estimator, which achieves the Cramer-Rao lower bound, has a covariance matrix which can be expressed as (see eqn. 6.4 in [10]):

$$\Sigma_{\text{MBAR}} = (\Pi - \Pi O \Pi)^+ - \Pi^{-1}, \tag{3.16}$$

where $\Pi = \text{diag}(\pi_1, \ldots, \pi_K)$ and the superscript $+$ denotes the Moore-Penrose pseudoinverse. In the specific case $K = 2$, one can simplify the expression for the asymptotic variance of $F_2 - F_1$, which works out to be:

$$\lim_{t \to \infty} t \text{Var}(F_2^t - F_1^t) = \frac{1}{\pi_1 \pi_2} [O_{12}^{-1} - 1]. \tag{3.17}$$

This is the same expression as the one for two-state MBAR on the last page of [20]. For further properties of the matrix $\Sigma_{\text{MBAR}}$, see Section 3.4.

## 3.3 Uniform distributions computations

In this section we derive the expressions for the variance used in the example from Section 1.4 of the main text. Let $(\mathcal{X}^t, \mathcal{K}^t)$ be the Markov process generated with the simulated tempering process while updating $F^t$ with stochastic approximation. More precisely, the density of the transition kernel for the process is $P_{F^t}$, where for each fixed $F \in \mathbb{R}^2$, $P_F = P_F^2 \circ P_F^1$, and

$$P_F^1(x', k'|x, k) = \delta(x' - x) \frac{\pi_{k'} e^{F_{k'} - H_{k'}(x)}}{\sum_{\ell=1}^2 \pi_\ell e^{F_\ell - H_\ell(x)}}, \quad P_F^2(x'', k''|x', k') = \delta(k'' - k') \rho_{k'}(x''). \tag{3.18}$$

The updates to $F^t$ are given by

$$F_k^{t+1} = F_k^t + \frac{1}{t+1} G_k(\mathcal{X}^{t+1}; F^t) = F_k^t + \frac{1}{t+1} \left( 1 - \frac{e^{F_k^t - H_k(\mathcal{X}^{t+1})}}{\sum_{\ell=1}^2 \pi_\ell e^{F_\ell - H_\ell(\mathcal{X}^{t+1})}} \right). \tag{3.19}$$

Since only free energy differences are of interest, we can define $\Delta^t = F_2^t - F_1^t$. We can also allow for a general gain matrix $\Gamma$, which in the present case will just be a scalar. The general update for $\Delta^t$ is

$$\Delta^{t+1} = \Delta^t + \frac{\Gamma}{t+1} \left( \frac{e^{-H_1(\mathcal{X}^{t+1})} - e^{\Delta^t - H_2(\mathcal{X}^{t+1})}}{\pi_1 e^{-H_1(\mathcal{X}^{t+1})} + \pi_2 e^{\Delta^t - H_2(\mathcal{X}^{t+1})}} \right). \tag{3.20}$$

By an abuse of notation, we will denote the term in parentheses in the above expression by $G(x; \Delta)$ (where $x = \mathcal{X}^{t+1}$ and $\Delta = \Delta^t$ above).

For fixed $\pi \in \mathcal{P}^\varepsilon(K)$, we will compute the asymptotic variance of the TSS free energy estimator $\Delta^t$ when $\mathcal{X}^t$ is sampled independently from $\rho_k(x)dx$ given $\mathcal{K}^t = k$, and where $\mathcal{K}^t$ is sampled with $\nu \geqslant 1$ iterations of the simulated tempering transition kernel $P_{F^t}^k$. To evaluate the expression for the asymptotic variance given in (3.7), we will require $R(\theta^\star)$, the integrated covariance, and $Dg(\theta^\star)$, the derivative of the mean field evaluated at $\theta^\star$. This last computation is particularly simple, because the Jacobian of the mean field for the free energies evaluated at $\theta^\star$ was found in Section 1 to be $-I_{K \times K}$.

We will use the expression in (3.16) to compute the asymptotic variance of MBAR, with the assumpton that a proportion $\pi_k$ of samples is drawn independently from $\rho_k(x)dx$.

### 3.3.1 Computing the variance for TSS and MBAR

To compute the variance of the TSS free energy estimates, we need to compute the integrated covariance (3.5). For certain sampling schemes, it is possible to compute the integrated covariance explicitly by solving an associated Poisson equation. Tan has derived some of these expressions for the case of free energy estimation when independent samples are available from each distribution $\rho_k(x)dx$ (see the Supplementary Material in [22]), while a more general result on the covariance structure of Gibbs samplers with Rao-Blackwellized estimators is available in [13]. In this section we calculate the covariance matrices by explicit summation.

Recall the exact form of the Hamiltonians in the uniform distributions example,

$$H_1(x) = \infty \cdot 1_{[-1+\delta,\delta]^c}(x), \quad H_2(x) = \infty \cdot 1_{[-\delta,1-\delta]^c}(x). \tag{3.21}$$

Observe that $G(x;\Delta)$ takes on three distinct values: $1/\pi_1$ when $x < -\delta$, $-1/\pi_2$ when $x > \delta$, and $(1-e^\Delta)/(\pi_1 + \pi_2 e^\Delta)$ when $x \in [-\delta,\delta]$. Moreover, the likelihood of observing any one of these values at the next iteration is determined by the value of $G$ at the present iteration. In other words, $Y^t = G(X^t;\Delta)$ is a three-state Markov chain for any $\Delta \in \mathbb{R}$, with transition matrix

$$P_\Delta = \begin{pmatrix} 1-2\delta & 2\delta & 0 \\ (1-2\delta) \cdot \frac{\pi_1}{\pi_1 + \pi_2 e^\Delta} & 2\delta & (1-2\delta) \cdot \frac{\pi_2 e^\Delta}{\pi_1 + \pi_2 e^\Delta} \\ 0 & 2\delta & 1-2\delta \end{pmatrix}, \tag{3.22}$$

for which the unique stationary distribution is

$$p_\Delta = \left((1-2\delta) \cdot \frac{\pi_1}{\pi_1 + \pi_2 e^\Delta}, 2\delta, (1-2\delta) \cdot \frac{\pi_2 e^\Delta}{\pi_1 + \pi_2 e^\Delta}\right). \tag{3.23}$$

The advantage of this perspective is that we can compute the integrated covariance (3.5) explicitly for this chain at $\Delta = \Delta^\star$ with $Y^0$ distributed according to $p_{\Delta^\star}$ (defined in (3.23)). To simplify the calculations we take $\pi_1 = \pi_2 = 1/2$, so that

$$P_{\Delta^\star} = \begin{pmatrix} 1-2\delta & 2\delta & 0 \\ (1-2\delta)/2 & 2\delta & (1-2\delta)/2 \\ 0 & 2\delta & 1-2\delta \end{pmatrix}, \tag{3.24}$$

To compute the integrated covariance we need to compute the powers of $P_{\Delta^\star}$, so we diagonalize the matrix. We have $P_{\Delta^\star} = ADA^{-1}$, where

$$A = \begin{pmatrix} 1 & 1 & 1 \\ 1 & 0 & 1-\frac{1}{2\delta} \\ 1 & -1 & 1 \end{pmatrix}, A^{-1} = \frac{1}{2}\begin{pmatrix} 1-2\delta & 4\delta & 1-2\delta \\ 1 & 0 & -1 \\ 2\delta & -4\delta & 2\delta \end{pmatrix}, D = \text{diag}(1, 1-2\delta, 0). \tag{3.25}$$

Thus, for integers $t \geqslant 0$ we have $P_{\Delta^\star}^t = AD^t A^{-1}$, which is explicitly

$$P_{\Delta^\star}^t = \frac{1}{2}\begin{pmatrix} 1-2\delta+p_\delta^t & 4\delta & 1-2\delta-p_\delta^t \\ 1-2\delta & 4\delta & 1-2\delta \\ 1-2\delta-p_\delta^t & 4\delta & 1-2\delta+p_\delta^t \end{pmatrix} = \frac{1}{2}\begin{pmatrix} p_\delta+p_\delta^t & 4\delta & p_\delta-p_\delta^t \\ p_\delta & 4\delta & p_\delta \\ p_\delta-p_\delta^t & 4\delta & p_\delta+p_\delta^t \end{pmatrix}, \tag{3.26}$$

where $p_\delta = (1-2\delta)$.

Now we assume $Y^t$ is the chain generated with (3.22) and $\Delta = \Delta^\star$. In this case $Y^t$ takes on only the values $\pm 2$ and $0$, and has $\mathbb{E}[Y^t] = \mathbb{E}_{\Delta^\star}[G(X; \Delta^\star)] = 0$. We use the values $\{-2, 0, 2\}$ to index the entries of the transition matrix, so that, for instance, $(P^1_{\Delta^\star})_{-2,-2} = (1 - 2\delta)$ and $(P^1_{\Delta^\star})_{0,2} = (1 - 2\delta)/2$. Noting that $(P^t_{\Delta^\star})_{2,2} = (P^t_{\Delta^\star})_{-2,-2}$ and $(P^t_{\Delta^\star})_{2,-2} = (P^t_{\Delta^\star})_{-2,2}$, we have

$$\begin{aligned}
\mathbb{E}[Y^0 Y^t] &= \sum_{y^0 \in \{-2,0,2\}} \sum_{y^t \in \{-2,0,2\}} y^0 y^t \mathbb{P}(Y^t = y^t | Y^0 = y^0) p_{\Delta^\star}(y^0) \\
&= \sum_{y^0 \in \{\pm 2\}} \sum_{y^t \in \{\pm 2\}} y^0 y^t \mathbb{P}(Y^t = y^t | Y^0 = y^0) p_{\Delta^\star}(y^0) \\
&= 4 \left( \frac{1-2\delta}{2} [(P^t_{\Delta^\star})_{-2,-2} - (P^t_{\Delta^\star})_{-2,2}] + \frac{1-2\delta}{2} [(P^t_{\Delta^\star})_{2,2} - (P^t_{\Delta^\star})_{2,-2}] \right) \\
&= 4(1-2\delta)[(P^t_{\Delta^\star})_{-2,-2} - (P^t_{\Delta^\star})_{-2,2}], \\
&= 4(1-2\delta) p_\delta^t. \tag{3.27}
\end{aligned}$$

Thus,

$$\sum_{t=1}^\infty \mathbb{E}[Y^0 Y^{vt}] = 4(1-2\delta) \sum_{t=1}^\infty p_\delta^{vt} = 4(1-2\delta) \frac{p_\delta^v}{1 - p_\delta^v}. \tag{3.28}$$

The variance of TSS when using $v$ moves per estimator update is

$$\Sigma^v_{\text{TSS}} = \text{Var}(Y^0) + 2 \sum_{t=1}^\infty \text{Cov}(Y^0, Y^{vt}) = 4(1-2\delta) + 8(1-2\delta) \frac{p_\delta^v}{1 - p_\delta^v}. \tag{3.29}$$

On the other hand, from (3.17) the variance of MBAR is

$$\frac{1}{\pi_1 \pi_2} [O_{12}^{-1} - 1] = 4 \left[ \left( \int_{-1+\delta}^{1-\delta} \frac{e^{-H_1(x) - H_2(x)}}{\pi_1 e^{-H_1(x)} + \pi_2 e^{-H_2(x)}} dx \right)^{-1} - 1 \right] \tag{3.30}$$

$$= 4 \left[ \frac{1}{2\delta} - 1 \right] = \frac{2 p_\delta}{\delta}. \tag{3.31}$$

Recalling that $p_\delta = 1 - 2\delta$, we find $\Sigma^v_{\text{TSS}} \leq \Sigma_{\text{MBAR}}$ whenever $v \geq 1 - \log(1+2\delta)/\log(1-2\delta)$, and we note that the latter inequality is satisfied by all $v \geq 2$ if $\delta \in (0, 1/2)$.

## 3.4 General proof of variance reduction for TSS

In this section we provide the details of the proof of Proposition 3 in the main text. Let $\Pi = \text{diag}(\pi_1, \ldots, \pi_K)$ and let $O$ denote the overlap matrix defined in (3.1). Since we are considering asymptotic covariance matrices, we assume the rung distribution $\pi$ has converged to some fixed vector with strictly positive entries. In particular, $\Pi$ is positive definite.

The covariance matrix of the TSS free energy difference estimators was derived in Section 3.1; the MBAR estimator and its covariance matrix were defined in 3.2. The two covariance matrices are denoted $\Sigma_{\text{TSS}}$ and $\Sigma_{\text{MBAR}}$ respectively, and both are defined below for convenience. For either covariance matrix $\Sigma$, the variance of the free energy difference estimates $F_j - F_i$ is obtained through the quadratic form

$$(\mathbf{e}_j - \mathbf{e}_i)^\top \Sigma (\mathbf{e}_j - \mathbf{e}_i), \tag{3.32}$$

where $\mathbf{e}_i, \mathbf{e}_j \in \mathbb{R}^K$ are unit basis vectors with the $i^{\text{th}}$ and $j^{\text{th}}$ entries respectively equal to 1. The claim of Theorem 2 in the main text is proved if we can establish that

$$(\mathbf{e}_j - \mathbf{e}_i)^\top (\Sigma_{\text{MBAR}} - \Sigma_{\text{TSS}})(\mathbf{e}_j - \mathbf{e}_i) \geqslant 0, \quad \forall i, j \in [K]. \tag{3.33}$$

**Theorem 1** *Suppose that the distributions $\rho_i, i \in [K]$ are such that the overlap matrix $O$ is irreducible. Let $\pi \in \mathcal{P}^\varepsilon(K)$ and suppose MBAR uses a proportion $\pi_k > 0$ of independent samples from $\rho_k(x)dx$, whereas for each fixed $F \in \mathbb{R}^K$ and the same $\pi = (\pi_1, \ldots, \pi_K)$, TSS uses independent samples from $p_{\pi,F}(x,k)dxC(dk)$. Then, for any $i, j \in [K]$, the asymptotic variance of the estimate $F_i - F_j$ for TSS is smaller than that of MBAR.*

**Proof** Recall that the asymptotic covariance matrix of the TSS free energy difference estimates was determined in Section 3.1 to be $\Sigma_{\text{TSS}} = O$ (see (3.8)). The asymptotic covariance matrix of MBAR is (see e.g. equation 6.4, [10], or [5, 21, 6])

$$\Sigma_{\text{MBAR}} = (\Pi - \Pi O \Pi)^+ - \Pi^{-1}, \tag{3.34}$$

where $\Pi = \text{diag}(\pi_1, \ldots, \pi_K)$ and the superscript $+$ denotes the Moore-Penrose pseudoinverse.

The positive semidefinite symmetric matrix $\Pi^{1/2} O \Pi^{1/2}$ admits 1 as a positive eigenvalue with positive normalized eigenvector $\mathbf{e} = \Pi^{1/2} \mathbf{1}$. Indeed,

$$(\Pi^{1/2} O \Pi^{1/2} \mathbf{e})_i = \pi_i^{1/2} \sum_{j \in [K]} O_{ij} \pi_j = \pi_i^{1/2} \sum_{j \in [K]} \pi_j \int_S \frac{\rho_i(x) \rho_j(x)}{\sum_{k \in [K]} \pi_k \rho_k(x)} dx = \pi_i^{1/2} \int_S \rho_i(x) dx = \pi_i^{1/2}. \tag{3.35}$$

Since $\Pi^{1/2} O \Pi^{1/2}$ is also irreducible, the Perron-Frobenius theorem implies that $\Pi^{1/2} O \Pi^{1/2}$ admits 1 as a simple eigenvalue and all other eigenvalues (which are necessarily non-negative real numbers) have magnitude less than 1. Therefore, $B = I - \Pi^{1/2} O \Pi^{1/2}$ and $A = \Pi - \Pi O \Pi = \Pi^{1/2} B \Pi^{1/2}$ are positive semidefinite matrices with one-dimensional kernels. The kernels of $B$ and $A$ are spanned by the normalized eigenvectors $\mathbf{e} = \Pi^{1/2} \mathbf{1}$ and $\mathbf{e}' = \mathbf{1}/K^{1/2}$, respectively. The Moore-Penrose pseudoinverse $A^+$ has the property $AA^+ = I - \mathbf{e}' \mathbf{e}'^\top$. We compute

$$A(\Sigma_{\text{MBAR}} - \Sigma_{\text{TSS}}) = I - \mathbf{e}' \mathbf{e}'^\top + \Pi(I - O\Pi)\Pi^{-1} - \Pi(I - O\Pi)O = (\Pi O)^2 - \mathbf{e}' \mathbf{e}'^\top. \tag{3.36}$$

Using that $\mathbf{e}'^\top A = (A\mathbf{e}')^\top = \mathbf{0}^\top$, we further compute

$$A(\Sigma_{\text{MBAR}} - \Sigma_{\text{TSS}})A = (\Pi O)^2 A - \mathbf{e}' \mathbf{e}'^\top A = (\Pi O)^2 A = (\Pi O)A(O\Pi), \tag{3.37}$$

which is positive semidefinite. If $\mathbf{v} \in \langle \mathbf{1} \rangle^\perp$, then $\mathbf{u} = A^+ \mathbf{v} \in \langle \mathbf{1} \rangle^\perp$ has the property $\mathbf{v} = A\mathbf{u}$. Then

$$\mathbf{v}^\top (\Sigma_{\text{MBAR}} - \Sigma_{\text{TSS}}) \mathbf{v} = \mathbf{u}^\top A (\Sigma_{\text{MBAR}} - \Sigma_{\text{TSS}}) A \mathbf{u} = \mathbf{u}^\top (\Pi O) A (O\Pi) \mathbf{u} = \mathbf{w}^\top A \mathbf{w} \geq 0, \tag{3.38}$$

where $\mathbf{w} = (O\Pi)\mathbf{u}$. The theorem follows because $\mathbf{v} = \mathbf{e}_j - \mathbf{e}_i$ lie in the orthogonal complement $\langle \mathbf{1} \rangle^\perp$ for all $1 \leq i, j \leq K$. □

# 4 Proof of convergence

The proof of convergence relies on the well-known technique of relating the stochastic system to the (deterministic) mean-field system of coupled ordinary differential equations (ODEs), and establishing some form of convergence of the deterministic system. Let us sketch the proof. We rely on existing results that phrase the convergence result in terms of regularity properties of the Markov process and stability properties of the system of ODEs. We observe that the quantity

$$\frac{e^{F_k - F_k^\star}}{\sum_{j \in [K]} \pi_j e^{F_j - F_j^\star}} \tag{4.1}$$

can be controlled as a function of $F$ and $\varepsilon > 0$, uniformly in $\pi \in \mathcal{P}^\varepsilon(K)$, and consequently we can say a lot about the system without knowing anything about $\pi$ specifically. The key step is to find a Lyapunov function that is uniform in $\pi$; see Lemma 2 and Proposition 1. We then observe that $F(t) \to F^\star$ regardless of $\pi \in \mathcal{P}^\varepsilon(K)$, and consequently the system of mean field equations can be approximated as

$$1 - \frac{e^{F_k(t) - F_k^\star}}{\sum_{\ell \in [K]} \pi_\ell(\theta(t)) e^{F_\ell(t) - F_\ell^\star}} \approx 0 \tag{4.2}$$

$$\frac{e^{F_k(t) - F_k^\star}}{\sum_{\ell \in [K]} \pi_\ell(\theta(t)) e^{F_\ell(t) - F_\ell^\star}} (\mu_{km}^\star - \mu_{km}(t)) \approx \mu_{km}^\star - \mu_{km}(t) \tag{4.3}$$

$$\frac{\pi_k(\theta(t))}{\gamma_k(\mu(t))} \cdot \frac{e^{F_k(t) - F_k^\star}}{\sum_{\ell \in [K]} \pi_\ell(\theta(t)) e^{F_\ell(t) - F_\ell^\star}} - o_k \approx \frac{\pi_k(\theta(t))}{\gamma_k(\mu(t))} - o_k. \tag{4.4}$$

The rest of the work is establishing that under a stability assumption on the dynamics for $o$, the whole system is globally asymptotically stable. We then use a converse Lyapunov theorem to produce a Lyapunov function necessary for the proof of convergence of the stochastic approximation procedure.

The general results we use are Theorem 5.5 and Proposition 6.1 in [1]. To apply this theorem, we need to check three conditions: a stability condition (Condition (A1) in [1]), a drift condition (DRI), and a condition on the gain sequences (A4). The stability condition requires us to find a global Lyapunov function, which is done in Section 4.2, while the drift condition requires us to check various regularity properties, which is done in Section 4.4.

We will make further assumptions on the compactness of $\Theta$ (with regard to the domain of $o = (o_1, \ldots, o_K)$; see (4.17)), but we note that none of these assumptions are used in establishing the existence of a global Lyapunov function for the free energy estimates $F$. When there exists a global Lyapunov function, [1] outlines a reprojection technique that offers stronger convergence results (see also [11] for discussion on reprojection techniques). We have not implemented such a technique, because in general we observe boundedness and eventual convergence of the estimates.

## 4.1 Setup

We recall the definitions from the main text. The parameter of interest is $\theta = (F, \mu, o) = (Z, \xi, o)$, and the rung distributions $\pi(\theta)$ and $\gamma(\mu)$ are user-defined deterministic functions $\pi : \Theta \to \mathcal{P}^\varepsilon(K)$ and $\gamma : \mathbb{R}^{KM} \to \mathcal{P}^\varepsilon(K)$, where we recall that $\mathcal{P}^\varepsilon(K)$ is the space of probability distributions $\omega$ on $[K] = \{1, \ldots, K\}$ with $\omega_k \geqslant \varepsilon$ for all $k \in [K]$. The functions $\psi_1(x), \ldots, \psi_m(x)$ can be evaluated for every $x \in \mathcal{S}$, and their averages with respect to $\rho_1(x)dx, \ldots, \rho_K(x)dx$ define the quantities $\mu_{km}^\star$, which are assumed finite. In a standard TSS simulation, the $\mu_{km}^\star$ are used to determine the coordinate invariant distribution defined in the main text and discussed in Section 5.1.3 of the present supplementary material.

We will demonstrate the convergence of $F^t$ and $\mu^t$ in the stochastic algorithm

$$F_k^{t+1} = F_k^t + \frac{1}{t+1}(1 - R_k^t) \tag{4.5}$$

$$\mu_{km}^{t+1} = \mu_{km}^t + \frac{1}{t+1} R_k^t (\psi_m(X^{t+1}) - \mu_{km}^t) \tag{4.6}$$

$$o_k^{t+1} = o_k^t + \frac{1}{t+1}\left(\frac{1_{\{k\}}(\mathcal{K}^{t+1})}{\gamma_k(\mu^t)} - o_k^t\right), \tag{4.7}$$

where

$$R_k^t = R_k(X^{t+1}; \theta^t), \quad R(x, \theta) = \frac{e^{F_k - H_k(x)}}{\sum_{\ell \in [K]} \pi_\ell(\theta) e^{F_\ell - H_\ell(x)}}. \tag{4.8}$$

We will only assume that $\mathbf{1} = (1, \ldots, 1) \in \mathbb{R}^K$ is a globally asymptotically stable equilibrium for the system

$$\dot{o}_k(t) = \frac{\pi_k(F^\star, \mu^\star, o(t))}{\gamma_k(\mu^\star)} - o_k(t), \quad k = 1, \ldots, K. \tag{4.9}$$

In particular, this implies $\pi_k(\theta(t)) \to \gamma_k(\mu^\star)$ for general functions $\pi$ (not just $\pi(\theta) = \gamma(\mu)$), which will consequently imply the convergence in distribution of the rung process $\mathcal{K}^t$ to $\gamma(\mu^\star)$.

## 4.2 Stability assumption (A1)

The goal of this section is to find a Lyapunov function which satisfies Assumption (A1) of [1] for the mean field equations (1.13), (1.14) and (1.7). This is done in Proposition 3; the conditions (A1) are listed immediately preceding the proposition, while the relevant definitions from ODE stability theory are provided in Section 4.2.2.

### 4.2.1 Change of coordinates

The stochastic approximation recursions (4.5) depend on $F$ only up to a constant shift $F_k \to F_k + \alpha$ common to all $k \in [K]$. It will be simpler to apply results from stability theory if we identify a unique equilibrium point for the mean field of $F$. To this end, we define the vector subspace $\mathbb{L} = \text{span}(\vec{v}) = \{\alpha \vec{1} : \alpha \in \mathbb{R}\}$ and the norm

$$\|F\|_{\mathbb{L}}^2 = \inf_{\alpha \in \mathbb{R}} \|F - \alpha \vec{1}\|^2. \tag{4.10}$$

The space $(\mathbb{R}^K/\mathbb{L}, \|\cdot\|_{\mathbb{L}})$ is isometric to $(\mathbb{R}^{K-1}, \|\cdot\|)$, where $\|\cdot\|$ denotes the standard Euclidean 2-norm. Observe that

$$\|F\|_{\mathbb{L}}^2 = \|F - \text{proj}_{\mathbb{L}}(F)\|^2 = \sum_{k \in [K]} (F_k - \langle F \rangle)^2 = \|F - \langle F \rangle\|^2, \tag{4.11}$$

where

$$\langle F \rangle = \frac{1}{K} \sum_{k \in [K]} F_k. \tag{4.12}$$

In this norm, we wish to establish $\|F(t) - F^\star\|_{\mathbb{L}} \to 0$ as $t \to \infty$ in a precise sense defined below as global asymptotic stability.

Now, let $\bar{F}_k(t) = F_k(t) - \langle F_k(t)\rangle$, and note that $\|\bar{F}(t)\| = \|\bar{F}(t)\|_\mathbb{L}$. To further simplify notation, we also assume that $F_k^\star = 0$ for all $k \in [K]$, which amounts to an additional change of variables $\bar{F}_k \to \bar{F}_k - F_k^\star$. The dynamics for $\bar{\theta} = (\bar{F}, \mu, o)$ are

$$\dot{\bar{F}}_k(t) = \left[\frac{1}{K}\sum_{j\in[K]}\frac{e^{\bar{F}_j(t)}}{\sum_{\ell\in[K]}\pi_\ell(\bar{\theta}(t))e^{\bar{F}_\ell(t)}}\right] - \frac{e^{\bar{F}_k(t)}}{\sum_{\ell\in[K]}\pi_\ell(\bar{\theta}(t))e^{\bar{F}_\ell(t)}} \quad (4.13)$$

$$\dot{\mu}_{km}(t) = \frac{e^{\bar{F}_k(t)}}{\sum_{\ell\in[K]}\pi_\ell(\bar{\theta}(t))e^{\bar{F}_\ell(t)}}(\mu_{km}^\star - \mu_{km}(t)) \quad (4.14)$$

$$\dot{o}_k(t) = \frac{\pi_k(\bar{\theta}(t))}{\gamma_k(\mu(t))} \cdot \frac{e^{\bar{F}_k(t)}}{\sum_{\ell\in[K]}\pi_\ell(\bar{\theta}(t))e^{\bar{F}_\ell(t)}} - o_k(t) \quad (4.15)$$

We concisely write this system as $\dot{\bar{\theta}} = \bar{g}(\bar{\theta})$, where $\bar{g} = (\bar{g}^F, g^\mu, g^o)$ is defined by (4.13), (4.14) and (4.15). To summarize, we seek to show that $\|\bar{F}(t)\|_\mathbb{L} \to 0$, $\|\mu - \mu^\star\| \to 0$ and $\|o_k(t) - \mathbf{1}\| \to 0$. Note that the stability properties of $(\bar{F}, \mu, o)$ are exactly the same as those of $(F, \mu, o)$, owing to the equality

$$\frac{e^{F_k(t)}}{\sum_{\ell\in[K]}\pi_\ell(\theta(t))e^{F_\ell(t)}} = \frac{e^{F_k(t) - \langle F(t)\rangle}}{\sum_{\ell\in[K]}\pi_\ell(\theta(t))e^{F_\ell(t) - \langle F(t)\rangle}}. \quad (4.16)$$

Finally, we observe that the state space for $o(t)$ can be assumed compact. Since $\gamma \in \mathcal{P}^\varepsilon(K)$, from the integrated form of the recurrence (4.7),

$$o_k^T = \frac{1}{T}\sum_{t=1}^{T}\frac{\mathbf{1}_{\{k\}}(\mathcal{K}^t)}{\gamma_k(\mu^{t-1})}, \quad (4.17)$$

we infer $o^t$, hence $o(t)$, lie in the convex hull $O^\varepsilon \subset \mathbb{R}^K_{\geq 0}$ defined by $o_k^t \geq 0$ for all $k \in [K]$ and $\sum_{k\in[K]} o_k^t \leq \varepsilon^{-1}$. We thus redefine $\Theta$ as $\Theta = \mathbb{R}^K \times \mathbb{R}^{KM} \times O^\varepsilon$.

### 4.2.2 Producing a Lyapunov function

We provide some definitions from stability theory. These will be used to establish convergence of $F(t)$ and $\mu(t)$, and if needed, convergence of $o(t)$. We use an $\varepsilon$-$\delta$ formulation of asymptotic stability, which are equivalent to modern formulations in terms of classes functions (see, for instance, [15, 9, 23]). Consider a general non-autonomous system

$$\dot{z}(t) = h(z(t), t). \quad (4.18)$$

It is always assumed that $h$ is everywhere continuously differentiable and therefore satisfies the mild regularity properties required for existence and uniqueness of solutions to ordinary differential equations (see for instance [7]). The variable $z$ takes values in $Z \subset \mathbb{R}^{n_z}$. We adopt the convention that $z(t; z_0, t_0)$ denotes a solution of (4.18) with initial condition $z(t_0; z_0, t_0) = z_0$. If $t_0$ is not specified and the system is autonomous, we sometimes use $z(t; z_0)$ to denote the solution with $z(0; z_0) = z_0$.

The following definition is taken from properties (ii), (v) in Section 3 of [15].

**Definition 1** *The equilibrium $z^\star$ of the non-autonomous system is called uniformly asymptotically stable in the large if the following two properties are satisfied:*

1. *(uniform stability) for any $\eta > 0$, there is $\delta = \delta(\eta) > 0$, independent of $t_0$, such that for any $|z_0 - z^\star| < \delta$ and $t_0 \geq 0$, we have $|z(t; z_0, t_0) - z^\star| < \eta$ for all $t \geq t_0$.*

2. (*uniform asymptotic stability*) for any $\eta > 0$ and $\rho > 0$ there is $T = T(\rho, \eta) < \infty$, independent of $t_0$, such that $|z(t; z_0, t_0) - z^\star| < \eta$ whenever $t \geq t_0 + T$, $t_0 \geq 0$, and $|z_0 - z^\star| < \rho$.

We will also refer to uniform asymptotic stability in the large by the name **global asymptotic stability**.

A standard way of establishing global asymptotic stability is through the use of Lyapunov functions.

**Definition 2** *A function $\alpha : \mathbb{R}_{\geq 0} \to \mathbb{R}_{\geq 0}$ is said to be of class $\mathcal{K}_\infty$ if it is continuous, zero at zero, strictly increasing, and satisfies $\lim_{c \to \infty} \alpha(c) = \infty$. A function $W : D \subset \mathbb{R}^n \to \mathbb{R}_{\geq 0}$ is said to be positive definite if it is continuous, $W(0) = 0$, and $W(x) > 0$ for all $x \neq 0$.*

*Let $V(t,z) : [0,\infty) \times Z \to [0,\infty)$ be a function defined on some subset $Z \subset \mathbb{R}^n$. Suppose that there exist functions $\alpha_1, \alpha_2 \in \mathcal{K}_\infty$ such that*

$$\alpha_1(\|z\|) \leq V(t,z) \leq \alpha_2(\|z\|), \quad t \geq 0, z \in Z. \tag{4.19}$$

*Suppose also that there is a positive definite function $W : \mathbb{R}_{\geq 0} \to \mathbb{R}_{\geq 0}$ such that*

$$V(t,z) \geq W(\|z\|), \quad t \geq 0, z \in Z. \tag{4.20}$$

*We call such a function $V$ a **candidate Lyapunov function**. A function $V$ is said to be a **Lyapunov function** for the system (4.18) if, in addition, there is a positive definite function $\sigma : \mathbb{R}_{\geq 0} \to \mathbb{R}_{\geq 0}$ such that, for all $t \geq 0$*

$$\frac{d}{dt} V(t, z(t)) \leq -\sigma(\|z(t)\|), \quad t \geq 0. \tag{4.21}$$

We remark that the functions $\alpha_1, \alpha_2$ are said to be of class $\mathcal{K}_\infty$ (see e.g. [9]).

Let $z = (x,y) \in X \times Y$ for some $X, Y \subset \mathbb{R}^{n_x} \times \mathbb{R}^{n_y}$ and consider the system

$$\begin{pmatrix} \dot{x}(t) \\ \dot{y}(t) \end{pmatrix} = \begin{pmatrix} f(x(t), y(t)) \\ g(x(t), y(t)) \end{pmatrix}, \tag{4.22}$$

which we write as $\dot{z} = h(z)$ with $h(z(x,y)) = (f(x,y), g(x,y))$. We use the norm $|z| = |x| + |y|$, so in particular $|z - z^\star| \geq |x - x^\star| + |y - y^\star|$. Let $z_0 = (x_0, y_0)$ and let $y(t; z_0)$ denote the solution of (4.22) with initial condition $(x_0, y_0)$. We can view $\dot{x} = f(x, y(t; z_0))$ as a time-dependent system for $x$ with initial condition $x_0$. Below, we study the stability properties of this time-dependent system for any trajectory $y(t)$, $y : [0, \infty) \to Y$.

**Lemma 1** *Let $U(x)$ be a candidate Lyapunov function in the sense of Definition 2, and suppose that there is a positive definite function $\sigma$ such that*

$$\sigma(|x - x^\star|) \leq \inf_{y \in Y} -\langle \nabla U(x), f(x,y) \rangle. \tag{4.23}$$

*Then the state $x^\star$ is a globally asymptotically stable equilibrium for the time-dependent system*

$$\dot{x} = f(x, y(t)) \tag{4.24}$$

*uniformly in the smooth map $y(t)$, $y : [0, \infty) \to Y$.*

**Proof** Let $y(t)$ be any smooth trajectory $y: [0, \infty) \to Y$, let $x_0 \in X$, and denote by $x(t; x_0)$ the solution of (4.24). Using the chain rule and the inequality (4.23), we find

$$\frac{d}{dt} U(x(t;x_0)) \leq \sup_y \langle \nabla U(x(t;x_0)), f(x(t;x_0), y) \rangle \leq -\sigma(|x(t;x_0) - x^\star|).$$

By Theorem 22 in [15], the equilibrium $x^\star$ is uniform-asymptotically stable "in the large" (i.e., globally asymptotically stable). It follows from the uniformity of the bound (4.23) that the stability is uniform in the trajectory $y(\cdot)$, in the sense that in the second condition of Definition 1, $T(\rho, \eta)$ can be chosen independently of $y(\cdot)$ (see the proof of Theorem 13, or the proof of Theorem 3.1 in section X.3 of [7]). □

We first identify a candidate Lyapunov function for the system $\bar{F}$ which will be suitably uniform with respect to $o$. We then use Lemma 1 to establish the global asymptotic stability of the equilibrium $(F^\star, \mu^\star)$.

**Lemma 2** *For any $\gamma \in \mathcal{P}^\varepsilon(K)$, the function*

$$V_\gamma(F) = \log\left(\sum_{k \in [K]} \gamma_k e^{F_k}\right) - \sum_{k \in [K]} \gamma_k F_k$$

*is a candidate Lyapunov function for the equilibrium $F^\star = 0$, in the sense that $V_\gamma$ is positive definite, and there exist $\alpha_1, \alpha_2$ as in Definition 2 such that $\alpha_1(\|F\|_\mathbb{L}) \leq V_\gamma(F) \leq \alpha_2(\|F\|_\mathbb{L})$.*

**Proof** Let $r(F) = (r_1(F), \ldots, r_K(F))$ be defined by

$$r_k(F) = \frac{\gamma_k e^{F_k}}{\sum_{\ell \in [K]} \gamma_\ell e^{F_\ell}}$$

and note that $\gamma$ and $r(F)$ are mutually absolutely continuous for all $F \in \mathbb{R}^K$ and that $V_\gamma(F) = D_{\mathrm{KL}}(\gamma \| r(F))$, where $D_{\mathrm{KL}}(\cdot \| \cdot)$ is the relative entropy (also known as the Kullback-Leibler divergence):

$$V_\gamma(F) = \log\left(\sum_{k \in [K]} \gamma_k e^{F_k}\right) - \sum_{k \in [K]} \gamma_k F_k = -\sum_{k \in [K]} \gamma_k \log\left(\frac{r_k(F)}{\gamma_k}\right) = D_{\mathrm{KL}}(\gamma \| r(F)). \quad (4.25)$$

Relative entropy is known to be positive definite (apply Jensen's inequality). For any two probability distributions $\gamma', \gamma''$, $D_{\mathrm{KL}}(\gamma' \| \gamma'') = 0$ if and only if $\gamma' = \gamma''$; the condition $r_k(F) = \gamma_k$ for all $k \in [K]$ forces $F_k = \mathrm{const}$, which is equivalent to $\|F\|_\mathbb{L} = 0$. It is also easy to check that $V_\gamma(F) \to \infty$ linearly with $\|F\|_\mathbb{L} \to \infty$. Together with the fact that $\varepsilon \leq \gamma_k \leq 1 - \varepsilon$ for all $k \in [K]$, we find that there exist $\alpha_1, \alpha_2 \in \mathcal{K}_\infty$ such that $\alpha_1(\|F\|_\mathbb{L}) \leq V_\gamma(F) \leq \alpha_2(\|F\|_\mathbb{L})$ for all $F \in \mathbb{R}^K$. □

We record the following simple estimate for use in the subsequent proposition.

**Lemma 3** *For any $\gamma \in \mathcal{P}^\varepsilon(K)$,*

$$\frac{e^{F_k}}{\sum_{\ell \in [K]} \gamma_\ell e^{F_\ell}} \geq e^{-2\|F\|_\mathbb{L}}. \quad (4.26)$$

**Proof** Recall that $\|F\|_\mathbb{L}^2 = \sum_{k \in [K]} (F_k - \langle F \rangle)^2$, so for any $\ell, k \in [K]$ we have $|F_\ell - F_k| \leq |F_\ell - \langle F \rangle| + |\langle F \rangle - F_k| \leq 2\|F\|_\mathbb{L}$. Assume without loss of generality that $F_K = \max_{k \in [K]} F_k$. Then $e^{F_k - F_K} \leq 1$ for all $k \in [K]$, hence $\sum_{k \in [K]} \gamma_k e^{F_k - F_K} \leq 1$ and

$$\frac{e^{F_k}}{\sum_{\ell \in [K]} \gamma_\ell e^{F_\ell}} = \frac{e^{F_k - F_K}}{\sum_{\ell \in [K]} \gamma_\ell e^{F_\ell - F_K}} \geq e^{F_k - F_K} \geq e^{-2\|F\|_\mathbb{L}}. \quad (4.27)$$

□

**Proposition 1** *Let $\bar{\theta}^0 = (\bar{F}^0, \mu^0, o^0) \in \Theta$ be any initial condition, and let $\bar{\theta}(t) = \bar{\theta}(t; \bar{\theta}^0)$ denote the solution of the time-independent system $\dot{\bar{\theta}} = \bar{g}(\bar{\theta})$ defined by (4.13), (4.14), (4.15). The state $(F^\star, \mu^\star)$ is globally asymptotically stable for the time-dependent system*

$$\dot{\bar{F}}_k(t) = \left[\frac{1}{K}\sum_{j \in [K]} \frac{e^{\bar{F}_j(t)}}{\sum_{\ell \in [K]} \pi_\ell(\bar{\theta}(t)) e^{\bar{F}_\ell(t)}}\right] - \frac{e^{\bar{F}_k}}{\sum_{\ell \in [K]} \pi_\ell(\bar{\theta}(t)) e^{\bar{F}_\ell(t)}} \qquad (4.28)$$

$$\dot{\mu}_{km}(t) = \frac{e^{\bar{F}_k(t)}}{\sum_{\ell \in [K]} \pi_\ell(\bar{\theta}(t)) e^{\bar{F}_\ell(t)}} (\mu^\star_{km} - \mu_{km}) \qquad (4.29)$$

*uniformly in the trajectory $o(t) = o(t; \bar{\theta}^0)$.*

**Proof** For any fixed $\gamma \in \mathcal{P}^\varepsilon(K)$, define the function

$$U(\bar{F}, \mu) = V_\gamma(\bar{F}) + \frac{1}{2} \sum_{k \in [K]} \sum_{m=1}^{M} (\mu_{km} - \mu^\star_{km})^2. \qquad (4.30)$$

We have shown in Lemma 2 that $V_\gamma$ is a candidate Lyapunov function, and it is clear that any sum of quadratic terms is as well; $U$ is therefore a candidate Lyapunov function. If we differentiate $U(\bar{F}(t), \mu(t))$ with respect to $t$ and apply the chain rule, we get

$$\frac{d}{dt} U(\bar{F}, \mu) = \langle \nabla_F U(\bar{F}, \mu), \bar{g}^F(\bar{F}, \mu, o) \rangle + \langle \nabla_\mu U(\bar{F}, \mu), g^\mu(\bar{F}, \mu, o) \rangle, \qquad (4.31)$$

where in the above display $\nabla_F$ and $\nabla_\mu$ are derivatives with respect to $F$ and $\mu$ respectively, and the $t$ dependence of $\bar{F}, \mu$ and $o$ has been suppressed to avoid excessively many parentheses. To apply Lemma 1 with $x = (\bar{F}, \mu)$ and $y = o$, we need to bound the terms $\langle \nabla_F U, \bar{g}^F(\bar{F}, \mu, o) \rangle$ and $\langle \nabla_\mu U, g^\mu(\bar{F}, \mu, o) \rangle$ (defined in (4.13), (4.14)) in terms of $\|\bar{F}\| + \|\mu - \mu^\star\|$, uniformly in the trajectory $o(t)$.

We first handle the inner product involving $\bar{g}^F$. Writing $\bar{\theta} = (\bar{F}, \mu, o)$, we have

$$\langle \nabla V_\gamma(F), \bar{g}^F(\bar{\theta}) \rangle = \sum_{k \in [K]} \gamma_k \left(\frac{e^{\bar{F}_k}}{\sum_{\ell \in [K]} \gamma_\ell e^{\bar{F}_\ell}} - 1\right) \left(\left[\frac{1}{K}\sum_{j \in [K]} \frac{e^{\bar{F}_j}}{\sum_{\ell \in [K]} \pi_\ell(\bar{\theta}) e^{\bar{F}_\ell}}\right] - \frac{e^{\bar{F}_k}}{\sum_{\ell \in [K]} \pi_\ell(\bar{\theta}) e^{\bar{F}_\ell}}\right) \qquad (4.32)$$

$$= -\sum_{k \in [K]} \gamma_k \left(\frac{e^{\bar{F}_k}}{\sum_{\ell \in [K]} \gamma_\ell e^{\bar{F}_\ell}} - 1\right) \left(\frac{e^{\bar{F}_k}}{\sum_{\ell \in [K]} \pi_\ell(\bar{\theta}) e^{\bar{F}_\ell}} - 1\right) \qquad (4.33)$$

To establish the uniform descrescence property (4.23), it suffices to show that the sum above can be bounded below uniformly in $\pi \in \mathcal{P}^\varepsilon(K)$. For $c \geq 0$, define

$$\sigma^F(c) = \inf_{\gamma \in \mathcal{P}^\varepsilon(K)} \inf_{F': \|F'\|_\mathbb{L} = c} \sum_{k=1}^K \gamma_k \left(\frac{e^{F'_k}}{\sum_{\ell \in [K]} \gamma_\ell e^{F'_\ell}} - 1\right) \left(\frac{e^{F'_k}}{\sum_{\ell \in [K]} \gamma'_\ell e^{F'_\ell}} - 1\right). \qquad (4.34)$$

We claim that $\sigma^F$ is positive definite: $\sigma^F(c) = 0$ if and only if $c = 0$, and $\sigma^F > 0$ otherwise. Clearly $\sigma^F(0) = 0$, with the infimum being attained by $F' = (0, \ldots, 0)$. To see that $\sigma^F(c) > 0$ for $c > 0$, let $F' \in \mathbb{R}^K, \gamma' \in \mathcal{P}^\varepsilon(K)$ be arbitrary and define the two functions $f, g : \mathbb{R} \to \mathbb{R}$ by

$$f(w) = \frac{e^w}{\sum_{\ell \in [K]} \gamma_\ell e^{F'_\ell}} - 1, \quad g(w) = \frac{e^w}{\sum_{\ell \in [K]} \gamma'_\ell e^{F'_\ell}} - 1. \qquad (4.35)$$

Both $f$ and $g$ are monotonically increasing, so if $f(F'_k) - f(F'_\ell) \geqslant 0$ we must have $g(F'_k) - g(F'_\ell) \geqslant 0$, and if $f(F'_k) - f(F'_\ell) \leqslant 0$ we must have $g(F'_k) - g(F'_\ell) \leqslant 0$. Therefore,

$$\sum_{k \in [K]} \sum_{\ell \in [K]} \gamma_k \gamma_\ell (f(F'_k) - f(F'_\ell))(g(F'_k) - g(F'_\ell)) = \frac{\sum_{\ell \in [K]} \gamma_\ell e^{F'_\ell}}{\sum_{\ell \in [K]} \gamma'_\ell e^{F'_\ell}} \sum_{k \in [K]} \sum_{\ell \in [K]} \gamma_k \gamma_\ell \left| f(F'_k) - f(F'_\ell) \right|^2 \geqslant 0. \quad (4.36)$$

Since all terms in the above sum are positive, equality holds if and only if $f(F'_1) = \cdots = f(F'_K)$, which forces $F_1 = \cdots = F_K$ or, equivalently, $\|F\|_{\mathbb{L}} = 0$. Observe that

$$\sum_{k \in [K]} \gamma_k f(F'_k) g(F'_k) \geqslant \left( \sum_{k \in [K]} \gamma_k f(F'_k) \right) \left( \sum_{k \in [K]} \gamma_k g(F'_k) \right) = 0, \quad (4.37)$$

where the inequality follows by re-arranging (4.36) and the equality follows from the specific form of $f$, which satisfies $\sum_{k \in [K]} \gamma_k f(F'_k) = 0$ for any $F' \in \mathbb{R}^K$. Since $\mathcal{P}^\varepsilon(K)$ is compact and the set $\{F : \|F\|_{\mathbb{L}} = c\}$ is compact in $\mathbb{R}^{K-1}$, a standard argument shows that the infima in (4.34) are attained. It follows that $\sigma^F$ is positive definite and that

$$\sup_{\mu \in \mathbb{R}^{KM}, o \in O^\varepsilon} \langle \nabla_F V_\gamma(\bar{F}), \bar{g}^F(\bar{F}, \mu, o) \rangle \leqslant -\sigma^F(\|\bar{F}\|_{\mathbb{L}}). \quad (4.38)$$

Now we return to the second term in (4.31). Using the bound in Lemma 3, we have

$$\sum_{k \in [K]} \sum_{m=1}^M \frac{e^{\bar{F}_k}}{\sum_{j \in [K]} \pi(\bar{\theta}) e^{\bar{F}_j}} (\mu_{km} - \mu^\star_{km})^2 \geqslant e^{-2\|\bar{F}\|_{\mathbb{L}}} \sum_{k \in [K]} \sum_{m=1}^M (\mu_{km} - \mu^\star_{km})^2 = e^{-2\|\bar{F}\|_{\mathbb{L}}} \|\mu - \mu^\star\|^2. \quad (4.39)$$

Finally, let $\sigma^\mu(c, d) = e^{-2c} d$. It is straightforward to find a positive definite function $\sigma$ that $\sigma(c + d) \leqslant \sigma^F(c) + \sigma^\mu(c, d)$ (though $\sigma$ will satisfy $\lim_{c \to \infty} \sigma(c) = 0$; see [8] for this style of argument). From equations (4.38) and (4.39) we find

$$\sup_{o \in O^\varepsilon} \{ \langle \nabla_F U(\bar{F}, \mu), \bar{g}^F(\bar{F}, \mu, o) \rangle + \langle \nabla_\mu U(\bar{F}, \mu), g^\mu(\bar{F}, \mu, o) \rangle \} \leqslant -\sigma(\|\bar{F}\|_{\mathbb{L}} + \|\mu - \mu^\star\|). \quad (4.40)$$

Thus the conditions of Lemma 1 are satisfied, with $x^\star = (F^\star, \mu^\star)$ ($F^\star = 0$) and an arbitrary trajectory $y(t) = o(t)$. It follows that $(F^\star, \mu^\star)$ is a globally asymptotically stable equilibrium for the system $\dot{\bar{\theta}} = \bar{g}(\bar{\theta})$. $\square$

Proposition 1 establishes the global asymptotic stability of $(F^\star, \mu^\star)$ independently of the trajectory $o(t) = o(t; \bar{\theta}^0)$. We remark that this is enough to guarantee convergence of the stochastic iterates $(F^t, \mu^t)$ if $\pi(\bar{\theta}) = \gamma(\mu)$, i.e. we do not add any visit control mechanisms.

We now proceed to establish the global asymptotic stability of $(F^\star, \mu^\star, \mathbf{1})$ using only the assumptions that $\pi$ is continuously differentiable and that $\mathbf{1}$ is a globally asymptotically stable equilibrium for the system (4.9). We require a definition (see Definition 56.1 in [23]) and a Lemma.

**Definition 3** *Suppose that the system*

$$\dot{z}(t) = h(z(t), t) \quad (4.41)$$

*has a unique solution and unique equilibrium $z^\star$. Consider any function $s$ for which the system*

$$\dot{z}(t) = h(z(t), t) + s(z(t), t) \quad (4.42)$$

*has a unique solution, denoted $\tilde{z}(t;z_0,t_0)$. The equilibrium $z^\star$ for the system* (4.41) *is called* **totally stable** *if for any $\eta > 0$, there are $\delta_1(\eta) > 0$ and $\delta_2(\eta) > 0$ such that*

$$|z_0 - z^\star| < \delta_1, \quad |s(z(t),t)| < \delta_2 \qquad (4.43)$$

*together imply $|\tilde{z}(t;z_0,t_0)| < \eta$.*

**Lemma 4** *Suppose that h is continuously differentiable. If $z^\star$ is a globally asymptotically stable equilibrium for the system* (4.41), *it is totally stable.*

*Moreover, there is an $r_1 > 0$ such that, for any $0 < \delta < r_1$ and any $r_2 \in (\delta, r_1]$, there is a $\Delta > 0$ and a $T > 0$ such that, for any $t_0 \geq 0$ and any $z_0 \in \mathbb{R}^{n_z}$, $\delta \leq |z_0 - z^\star| \leq r_2$, the solution $z(t;t_0,z_0)$ of* (4.41) *satisfies $|z(t;t_0,z_0) - z^\star| < \delta$ for all $t \geq t_0 + T$, provided that $|s(t,z)| < \Delta$ for all $t \geq t_0, |z - z^\star| \leq r_1$.*

The proof of the first statement in Lemma 4 can be found in [23], Theorem 56.4; the proof of the second statement can be found in [7], Theorem 5.2, Section X.5.

**Lemma 5** *Let Y be compact, and suppose that $y^\star$ is a globally asymptotically stable equilibrium for the system*

$$\dot{y}(t) = g(x^\star, y(t)) \qquad (4.44)$$

*and that $\dot{x}(t) = f(x(t), y(t))$ satisfies the conditions of Lemma 1. Then $(x^\star, y^\star)$ is a globally asymptotically stable equilibrium of the system $(\dot{x}, \dot{y}) = (f(x(t), y(t)))$.*

**Proof** First we prove uniform stability. Let $\eta > 0$ be given. Since $y^\star$ is a globally asymptotically stable equilibrium of $\dot{y} = g(x^\star, y)$, it is totally stable by Lemma 4. Therefore, for the given $\eta > 0$, there exist $\delta_1 = \delta_1(\eta) > 0$ and $\delta_2 = \delta_2(\eta) > 0$ such that the unique solution $y(t;z_0)$ of the system

$$\dot{y}(t) = g(x(t), y(t)) = g(x^\star, y(t)) + s(t, y(t)), \quad y(0) = y_0 \qquad (4.45)$$

satisfies

$$|y(t;z_0) - y^\star| < \eta/2 \qquad (4.46)$$

for all $t \geq 0$ whenever $|y_0 - y^\star| < \delta_1$ and $|s(t,y)| < \delta_2$ for all $t \geq 0$ and $y \in B(y^\star, \eta)$.

Let $r(x,y) = g(x,y) - g(x^\star, y)$, and observe from (4.45) that $s(t,y) = r(x(t;z_0),y)$. Since $g$ is continuous, $r$ is continuous. By continuity of $r$, for each fixed $y$ there is $\beta_y = \beta_y(\delta_2) > 0$ such that

$$|r(x,y)| < \delta_2 \qquad (4.47)$$

whenever $|x - x^\star| < \beta_y$. Without loss of generality, we can assume the map $y \to \beta_y$ is continuous, and hence it achieves its minimum on the compact set $\overline{B(y^\star, \eta)}$; we denote $\beta = \min_{y \in \overline{B(y^\star, \eta)}} \beta_y > 0$.

Next, we use the uniform stability of $x^\star$ for $\dot{x} = f(x, y(t))$, which holds uniformly in the trajectory $y(\cdot)$, $y : [0, \infty) \to Y$. In particular, it holds for $y(t) = y(t;z_0)$, where $y(t;z_0)$ denotes the solution of (4.45) with $y(0;z_0) = y_0$. For the given $\beta > 0$, there is $\delta_x > 0$ such that

$$|x(t;z_0) - x^\star| < \min(\beta, \eta/2) \qquad (4.48)$$

for all $t \geq 0$ whenever $x_0 \in B(x^\star, \delta_x)$.

We claim that $\delta = \min(\delta_1, \delta_x) > 0$ satisfies the conditions of uniform stability. Let $z_0$ satisfy $|z_0 - z^\star| < \delta$, so by the triangle inequality $|x_0 - x^\star| < \delta$ and $|y_0 - y^\star| < \delta$. By the choice of $\delta \leq \delta_x$, (4.48) holds for all $t \geq 0$, so $x(t;z_0)$ stays within $\eta/2$ of $x^\star$. Since $x(t;z_0)$ also stays within $\beta$ of $x^\star$ for all $t \geq 0$, (4.47) holds, and therefore $|s(t,y)| = |r(x(t;z_0),y)| < \delta_2$ for all $t \geq 0$ and $y \in B(y^\star, \eta)$. Note also that $|y_0 - y^\star| < \delta \leq \delta_1$ by the assumption $|z_0 - z^\star| < \delta$. By the total stability of (4.45), we have $|y(t;z_0) - y^\star| < \eta/2$ for all $t \geq 0$. Thus we have $|z(t;z_0) - z^\star| \leq \eta/2 + \eta/2 = \eta$ for all $t \geq 0$, which establishes the uniform stability of $z^\star$ for the system $\dot{z} = g(z)$.

Next we establish equiasymptotic stability of $z^\star$. Let $\eta > 0$ and $\rho > 0$ be given. Using the total stability of (4.44), there are once again $\delta_1 = \delta_1(\eta) > 0$, $\delta_2 = \delta_2(\eta) > 0$ such that

$$|y(t;t_0,y_0) - y^\star| < \eta/2 \tag{4.49}$$

for all $t \geq t_0$ whenever $|y_0 - y^\star| < \delta_1$ and $|s(t,y)| < \delta_2$. Here, we write $y(t;t_0,y_0)$ to distinguish the solution starting from $y_0$ at $t_0$ (rather than $t_0 = 0$), and we note that $y(t;z_0) = y(t;t_0,y(t_0;z_0))$. Notice also that we may assume without loss of generality that $\delta_1 < \eta/2$, as otherwise the statement is vacuous.

To show that the perturbed system (4.45) still brings trajectories close to the equilibrium $y^\star$ of the unperturbed system, we apply the second part of Lemma 4 to the system (4.45) with $r_1 = \sup_{y,y' \in Y} |y - y'|$, $r_2 = r_1$ and the $\delta_1$ we have just obtained above. The theorem yields $\Delta = \Delta(\delta_1, r_2) > 0$ and $T_Y = T_Y(\eta, r_2)$ such that, the solution $y(t;t_0,y_0)$ of (4.45) with $y(t_0) = y_0$ satisfies

$$|y(t;z_0) - y^\star| < \delta_1 \tag{4.50}$$

for $t \geq T_Y + t_0$, provided that $|s(t,y)| < \Delta$ for all $t \geq t_0$ and $r_2 \geq |y_0 - y^\star| \geq \delta_1$.

Now we use the fact that $\dot{x} = f(x, y(t))$ is globally asymptotically stable uniformly in the trajectory $y(\cdot)$ to satisfy the bound $|s(t,y)| < \min(\delta_2, \Delta)$ for sufficiently large $t$. To this end, we again use the continuity of $g$ to find a $\delta_x = \delta_x(\delta_2, \Delta) > 0$ such that

$$|s(t,y)| < \min(\delta_2, \Delta) \tag{4.51}$$

holds whenever $x \in B(x^\star, \delta_x)$. This time the $\delta_x$ must work for all $y \in Y$, which is possible by compactness of $Y$. By the global asymptotic stability, for the given $\rho > 0$ there is $T_X = T_X(\rho, \eta, \delta_x) < \infty$ such that

$$|x(t;z_0) - x^\star| < \min(\eta/2, \delta_x) \tag{4.52}$$

whenever $|x_0 - x^\star| < \rho$ and $t \geq T_X$.

We now compile the estimates. Set $T = T(\rho, \eta) = T_X + T_Y$. We claim that for any $z_0 \in B(z^\star, \rho)$, we have $|z(t;z_0) - z^\star| < \eta$ for all $t \geq T$. First, it is straightforward that $x(t;z_0)$ arrives within $\eta/2$ of $x^\star$ by time $T_X$, owing to (4.52). By the same estimate, $x(t;z_0)$ is within $\delta_x$ of $x^\star$ at all times $t \geq T_X$ for any $z_0 \in B(z^\star, \rho)$, so (4.51) holds for all $t \geq T_X$ and $y \in Y$.

Now, at time $t_0 = T_X$, either $y(t;z_0)$ is in $B(y^\star, \delta_1)$ or it is not. If it is, then together with (4.51) we have that (4.49) holds by the total stability estimate, and hence $|y(t;z_0) - y^\star| < \eta/2$ for all $t \geq T_X$. If not, then by estimates (4.50) and (4.51) we have $|y(t;z_0) - y^\star| < \delta_1 < \eta/2$ for all $t \geq T_X + T_Y$. Either way, $|y(t;z_0) - y^\star| < \eta/2$ for all $t \geq T$. Together with estimate (4.52) and the triangle inequality, we obtain the desired result. □

**Proposition 2** *Suppose that (4.9) is uniformly globally asymptotically stable. Then the time-independent system $\dot{\bar{\theta}}(t) = \bar{g}(\bar{\theta}(t))$ is (uniformly) globally asymptotically stable.*

**Proof** We know from Proposition 1 that the equilibrium $(0, \mu^\star)$ is globally asymptotically stable for the system $(\dot{\bar{F}}, \dot{\mu}) = (\bar{g}^F, g^\mu)$. The proof of the present proposition follows by applying Lemma 5 with $x = (\bar{F}, \mu)$ and $y = o$ to the system $f = (\bar{g}^F, g^\mu)$ and $g = g^o$. □

Next, we establish the global asymptotic stability of the system (4.15) when $\pi$ is taken to be TSS's visit control.

**Lemma 6** *For any $\varepsilon_\pi \in [0, 1]$, let*

$$\pi_k^{\varepsilon_\pi}(o) = (1-\varepsilon_\pi)\pi_k(o) + \varepsilon_\pi \gamma_k(\mu^\star), \tag{4.53}$$

*where $\pi_k$ is the unregularized TSS visit control,*

$$\pi_k(o) = \frac{\gamma_k(\mu^\star) o_k^{-\eta}}{\sum_{\ell \in [K]} \gamma_\ell(\mu^\star) o_\ell^{-\eta}}. \tag{4.54}$$

*Then, for any $\varepsilon_\pi \in [0, 1]$, the system (4.9) with $\pi = \pi^{\varepsilon_\pi}$ has $\mathbf{1}$ as its globally asymptotically stable equilibrium.*

**Proof** When $\eta = 0$, we have $\pi_\pi^\varepsilon = \gamma$ for all $\varepsilon_\pi \in [0,1]$ and the differential equation (4.9) for $o$ becomes $\dot{o}_k = 1 - o_k$, which is trivially globally asymptotically stable. For $\eta > 0$ it suffices to exhibit a global Lyapunov function. Let $U(o) = \frac{1}{2}\sum_{k \in [K]} \gamma_k(1-o_k)^2$, and note that $U$ is a candidate Lyapunov function (in the sense of Definition 2 in Section 4.2). We have, for any $\varepsilon_\pi \in [0, 1f]$,

$$\langle \nabla U(o), \dot{o} \rangle = -\sum_{k \in [K]} \gamma_k(1-o_k)\left(\frac{\pi_k^{\varepsilon_\pi}}{\gamma_k(\mu^\star)} - o_k\right) \tag{4.55}$$

$$= -\sum_{k \in [K]} \gamma_k(1-o_k)\left((1-\varepsilon_\pi)\frac{o_k^{-\eta}}{\sum_{\ell \in [K]} \gamma_\ell o_\ell^{-\eta}} - (1-\varepsilon_\pi)o_k + \varepsilon_\pi - \varepsilon_\pi o_k\right) \tag{4.56}$$

$$= -(1-\varepsilon_\pi)\left[1 - E_\gamma[O] - \frac{E_\gamma[O^{1-\eta}]}{E_\gamma[O^{-\eta}]} + E_\gamma[O^2]\right] - \varepsilon_\pi E_\gamma[(1-O)^2], \tag{4.57}$$

where $E_\gamma[f(O)] = \sum_{k \in [K]} \gamma_k f(o_k)$ for any function $f : [K] \to \mathbb{R}$. By the correlation inequality (see (2.33) and (2.34)) applied to the two increasing functions $f(o) = o$ and $g(o) = -o^{-\eta}$ (note $\eta > 0$ and $o_k > 0$), we have

$$-E_\gamma[O^{1-\eta}] = E_\gamma[f(O)g(O)] \geq E_\gamma[f(O)]E_\gamma[g(O)] = -E_\gamma[O]E_\gamma[O^{-\eta}], \tag{4.58}$$

and so $E_\gamma[O^{1-\eta}]/E_\gamma[O^{-\eta}] \leq E_\gamma[O]$. Therefore,

$$1 - E_\gamma[O] - \frac{E_\gamma[O^{1-\eta}]}{E_\gamma[O^{-\eta}]} + E_\gamma[O^2] \geq 1 - 2E_\gamma[O] + E_\gamma[O^2] = E_\gamma[(1-O)^2]. \tag{4.59}$$

It follows from (4.57) that $\langle \nabla U(o), \dot{o} \rangle \leq -2U(o)$, which concludes the proof. □

Finally, we conclude the existence of an appropriate Lyapunov function which satisfies the assumptions below. We will invoke a classical converse Lyapunov theorem; see [9, 7, 15, 23].

**Assumption (A1).** Assume $\Theta$ is an open subset of $\mathbb{R}^q$, $g : \Theta \to \mathbb{R}^q$ is continuous and there exists a continuously differentiable function $w : \Theta \to [0, \infty)$ such that

1. There exists $M_0 > 0$ such that

$$\mathcal{L} = \{\theta \in \Theta : \langle \nabla w(\theta), g(\theta) \rangle = 0\} \subset \{\theta \in \Theta : w(\theta) < M_0\}, \tag{4.60}$$

2. There exists $M_1 \in (M_0, \infty]$ such that $\{\theta \in \Theta : w(\theta) \leq M_1\}$ is a compact set,

3. For any $\theta \in \Theta \setminus \mathcal{L}$, $\langle \nabla w(\theta), g(\theta) \rangle < 0$,

4. The closure of $w(\mathcal{L})$ has an empty interior.

**Proposition 3** *There exists a global Lyapunov function for the system $\dot{\bar{\theta}} = \bar{g}(\bar{\theta}(t))$ which satisfies Assumption (A1).*

**Proof** By Theorem 23 in [15], there is a function $w(\bar{\theta})$ and $\alpha_1, \alpha_2 \in \mathcal{K}_\infty$ such that $\alpha_1(\|\bar{\theta} - \bar{\theta}^\star\|) \leqslant w(\bar{\theta}) \leqslant \alpha_2(\|\bar{\theta} - \bar{\theta}^\star\|)$. The function $w$ vanishes only at the point $\bar{\theta}^\star = (0, \mu^\star, \mathbf{1})$, is positive otherwise, and smooth ($w \in C^\infty$). Moreover, $\langle \nabla w(\bar{\theta}), \bar{g}(\bar{\theta}) \rangle \leq -\sigma(\|\bar{\theta} - \bar{\theta}^\star\|)$ for some positive definite function $\sigma$.

Now we verify that this function $w$ satisfies Assumption (A1). In our case $\mathcal{L} = \{\bar{\theta}^\star\} = \{(0, \mu^\star, \mathbf{1})\}$ and $w(\bar{\theta}) = 0$ if and only if $\bar{\theta} = \bar{\theta}^\star \in \mathcal{L}$. Thus 1. is satisfied with any $M_0 > 0$.

Since $\alpha_1 \in \mathcal{K}_\infty$ and $w(\bar{\theta}) \geqslant \alpha_1(\|\bar{\theta} - \theta^\star\|)$, the function $w(\bar{\theta})$ is coercive ($w(\bar{\theta}) \to \infty$ as $\|\bar{\theta}\| \to \infty$) and therefore has compact level sets, so there is $M_1 > 0$ such that $\{\bar{\theta} \in \Theta : w(\bar{\theta}) \leqslant M_1\}$ is compact. Thus condition 2. can be satisfied, and we can choose any $M_0 \in (0, M_1)$ to satisfy condition 1.

Condition 3. follows immediately from the positive definiteness of the function $\sigma$.

Finally, $w(\mathcal{L}) = \{0\}$, and so the closure has empty interior. $\square$

## 4.3 Assumption (A4)

Assumption (A4) in [1] places conditions on the step sizes. The stochastic approximation algorithm with step size $\beta_k = 1/k$ satisfies

$$\sum_{k=1}^\infty \gamma_k = \infty, \quad \sum_{k=1}^\infty \gamma_k^\delta < \infty \tag{4.61}$$

for any $\delta > 1$. In the notation of [1], we need $\delta \in (1, p(1+\alpha)/(p+\alpha))$ for some $p \geqslant 2$ and $\alpha \in (0, 1]$. In our case we take $\alpha = 1$, $p \geqslant 2$, and we set $\varepsilon_k = C\gamma_k^\eta$ for some constant $C$ and $\eta$ such that

$$\frac{\delta - 1}{\alpha} \leqslant \eta \leqslant \frac{p - \delta}{p}, \tag{4.62}$$

taking $p$ as large as necessary. The sequence $\{\varepsilon_k\}_{k \geqslant 1}$ then satisfies Assumption (A4) of [1].

## 4.4 Drift conditions (DRI)

The drift conditions require us to verify basic properties to ensure ergodicity of the transition kernels for each fixed $\theta \in \Theta$. We reproduce the drift condition (DRI) in Section 6. of [1], which guarantees their Assumptions (A2) and (A3). Let $\Omega = \mathcal{S} \times [K]$. For any function $U : \Omega \to [1, \infty)$, define the norm

$$\|g\|_U = \sup_y \frac{|g(y)|}{U(y)}, \tag{4.63}$$

and let $\mathcal{L}_U = \{g : \Omega \to \mathbb{R}^{n_\theta} \mid \|g\|_U < \infty\}$.

**(DRI).** For any $\theta \in \Theta$, $P_\theta$ is $\psi$-irreducible and aperiodic. In addition, there exist a function $U : \Omega \to [1, \infty)$ and constants $p \geqslant 2$, $\beta \in [0, 1]$ such that for any compact set $\mathbb{K} \subset \Theta$,

**(DRI1)** there exist an integer $m$, constants $\lambda \in (0, 1)$, $b, \kappa, \delta > 0$ and a set $C \subset \Theta$ with a probability measure probability measure $q$ such that

1. $\sup_{\theta \in \mathbb{K}} P_\theta^m U^p(y) \leq \lambda U^p(y) + b 1_C(y)$,
2. $\sup_{\theta \in \mathbb{K}} P_\theta U^p(y) \leq \kappa U^p(y)$, for all $y \in \Omega$,
3. $\inf_{\theta \in \mathbb{K}} P_\theta^m(y, A) \geq \delta q(A)$, for all $y \in C$ and Borel-measurable $A \subset \Omega$

**(DRI2)** there exists a constant $c > 0$ such that, for all $y \in \Omega$,

1. $\sup_{\theta \in \mathbb{K}} |G_\theta(y)| \leq cU(y)$
2. $\sup_{(\theta, \theta') \in \mathbb{K}} |\theta - \theta'|^{-\beta} |G_\theta(y) - G_{\theta'}(y)| \leq cU(y)$

**(DRI3)** there exists a constant $c > 0$ such that, for all $(\theta, \theta') \in \mathbb{K} \times \mathbb{K}$,

1. $\|P_\theta g - P_{\theta'} g\|_U \leq c \|g\|_U |\theta - \theta'|^\beta$, for all $g \in \mathcal{L}_U$
2. $\|P_\theta g - P_{\theta'} g\|_{U^p} \leq c \|g\|_{U^p} |\theta - \theta'|^\beta$, for all $g \in \mathcal{L}_{U^p}$

### 4.4.1 Verification of (DRI)

It is well know that the simulated tempering sampler $P_\theta = P_\theta^x \circ P_\theta^k$ produces a Markov process which is $\psi$-irreducible and aperiodic, so we only sketch a proof below.

**Lemma 7** *For each fixed $\theta \in \Theta$, the chain defined by $P_\theta$ is $\psi$-irreducible and aperiodic.*

**Proof** We note that $\varphi = p_{\pi, F}(x, k) dx \mathcal{C}(dk)$ defines a positive measure on $\mathcal{S} \times [K]$, and we assume that the distributions are chosen such that $\mathcal{S} \times [K]$ is $\varphi$-communicating, hence the chain is $\varphi$-irreducible in the sense of Meyn and Tweedie [16] (see Proposition 4.2.2). It is also aperiodic, since the probability of staying at any $k \in [K]$ is nonzero and the dynamics in $x$ are assumed aperiodic. Thus the chain is $\psi$-irreducible and aperiodic. □

Note that the transition kernel $P_\theta$ is continuous in $(x, k)$, hence lower semi-continuous. An application of Fatou's lemma shows that the chain is therefore a T-chain (see p. 131, [16]). We use the following two results:

**Lemma 8** *(Theorem 6.2.9 in [16]) If a $\psi$-irreducible Markov chain with transition kernel $P(x, \cdot), x \in X$ is weak Feller and if $\mathrm{supp}(\psi)$ has nonempty interior, then the chain is a T-chain.*

Our chain is therefore a T-chain.

**Proposition 4** *(Theorem 6.2.5, (ii) in [16]) If $\{X^n\}$ is a $\psi$-irreducible T-chain then every compact set is petite.*

Note that for aperiodic and irreducible chains, as is the case here, every petite set is small (see [16], Theorem 5.5.7). We take $\beta = 1$ and any $p \geq 2$ for which Assumption (A4) is satisfied. We also take $U(y) = 1$. Let the compact set $\mathbb{K} \subset \Theta$ be given. Below, we verify the drift conditions.

**(DRI1)** As in [4], we take any $0 < \lambda < 1$, set $b = 1 - \lambda < 1$. Since the chain is an aperiodic T-chain, every compact set is small. Since $C = \mathcal{S} \times [K]$ is compact, $C$ is a small set, and

$$\sup_{\theta \in \mathbb{K}} (P_\theta U^p)(y) \leq \lambda U^p(y) + b 1_C(y) \equiv 1. \quad (4.64)$$

With $\kappa > 1$, we also have
$$\sup_{\theta \in \mathbb{K}} (P_\theta U^p)(y) \leqslant \kappa U^p(y) = \kappa. \tag{4.65}$$

Finally, since the whole set is small, there are $\delta > 0$, a positive integer $s \geqslant 1$ and a probability measure $q$ such that for any Borel-measurable $A \subset \mathcal{S} \times [K]$,
$$\inf_{\theta \in \mathbb{K}} P_\theta^s(y, A) \geqslant \delta q(A), \tag{4.66}$$

so the minorization condition is satisfied.

**(DRI2)** From (1.2), (1.3), (1.4) with $\theta = (Z, \xi, o)$ and $y = (x, \ell)$, the functions $G(y; \theta) = (G^Z(y; \theta), G^\xi(y; \theta), G^o(y; \theta))$ are defined by

$$G_k^Z(y; \theta) = \frac{e^{-H_k(x)}}{\sum_{\ell \in [K]} \pi_\ell(\theta) e^{-H_\ell(x)}/Z_\ell} - Z_k \tag{4.67}$$

$$G_{km}^\xi(y; \theta) = \frac{e^{-H_k(x)}}{\sum_{\ell \in [K]} \pi_\ell(\theta) e^{-H_\ell(x)}/Z_\ell} \psi_m(x) - \xi_{km} \tag{4.68}$$

$$G_k^o(y; \theta) = \frac{1}{\gamma_k(\mu)} \mathbf{1}_{\{k\}}(\ell) - o_k. \tag{4.69}$$

Note that each component is continuous in $(y, \theta)$. Moreover, $\Omega = \mathcal{S} \times [K]$ is compact, and the drift conditions allow us to restrict to a supremum over $\theta \in \mathbb{K}$ where $\mathbb{K}$ is a compact set. Since $(y, \theta) \to G(y; \theta)$ is continuous, it achieves its supremum on the compact set $\Omega \times \mathbb{K}$. Thus part 1. of condition (DRI2) is satisfied with $U(y) = 1$.

To see that
$$\sup_{(\theta, \theta') \in \mathbb{K} \times \mathbb{K}} \frac{|G(y, \theta) - G(y; \theta')|}{|\theta - \theta'|} \leqslant C, \tag{4.70}$$

note that $G^Z(y; \theta), G^\xi(y, \theta)$ and $G^o(y, \theta)$ are all continuously differentiable in $\theta$ for each fixed $y = (x, k)$, so (4.70) follows from the mean value theorem. Thus part 2. of (DRI2) is satisfied.

**(DRI3)** To establish the Lipschitz-style property of $P_\theta = P_\theta^k \circ P_\theta^x$ with respect to $\theta$, it suffices to establish it separately for both kernels, i.e. with $P^1 = P^x$ and $P^2 = P^k$, we want
$$\|P_\theta^i h\|_U \leqslant C_i \|h\|_U, \quad \|P_\theta^i h - P_{\theta'}^i h\|_U \leqslant C_i \|h\|_U |\theta - \theta'|^\beta. \tag{4.71}$$

Indeed, by the triangle inequality,
$$\|P_\theta h - P_{\theta'} h\|_U \leqslant \|P_\theta^2 (P_\theta^1 h - P_{\theta'}^1 h)\|_U + \|P_\theta^2 P_{\theta'}^1 h - P_{\theta'}^2 P_{\theta'}^1 h\|_U \tag{4.72}$$
$$\leqslant C_2 \|P_\theta^1 h - P_{\theta'}^1 h\|_U + C_2 \|P_{\theta'}^1 h\|_U |\theta - \theta'|^\beta \tag{4.73}$$
$$\leqslant 2 C_1 C_2 |\theta - \theta'|^\beta \|h\|_U. \tag{4.74}$$

For the first kernel, recall from the main text that $P_\theta^x(x', k'|x, k) = \mathcal{T}_k^{(n)}(x'|x)\mathrm{d}x'$, which does not depend on $\theta$. Thus (4.71) holds trivially for $i = 1$ and both $U$ and $U^p$. For the second kernel, note that for any

$g \in \mathcal{L}_U$,

$$P_\theta^k g(x,k) - P_{\theta'}^k g(x,k) \tag{4.75}$$

$$\sum_{k' \in [K]} \int_S g(x',k')(P_\theta^k(x',k'|x,k) - P_{\theta'}^k(x',k'|x,k))dx' \tag{4.76}$$

$$\sum_{k' \in [K]} \int_S g(x',k')\delta(x-x')\left(\frac{\pi_{k'}(\theta)e^{-H_{k'}(x')}/Z_{k'}}{\sum_{\ell \in [K]} \pi_\ell(\theta)e^{-H_\ell(x')}/Z_\ell} - \frac{\pi_{k'}(\theta')e^{-H_{k'}(x')}/Z_{k'}}{\sum_{\ell \in [K]} \pi_\ell(\theta')e^{-H_\ell(x')}/Z_\ell}\right)dx' \tag{4.77}$$

$$= \sum_{k' \in [K]} g(x,k')\left(\frac{\pi_{k'}(\theta)e^{-H_{k'}(x)}/Z_{k'}}{\sum_{\ell \in [K]} \pi_\ell(\theta)e^{-H_\ell(x)}/Z_\ell} - \frac{\pi_{k'}(\theta')e^{-H_{k'}(x)}/Z_{k'}}{\sum_{\ell \in [K]} \pi_\ell(\theta')e^{-H_\ell(x)}/Z_\ell}\right) \tag{4.78}$$

The term in parentheses above is continuously differentiable w.r.t. $\theta$, and by the mean value theorem its absolute value can be bounded above $c\|\theta - \theta'\|$, where $c > 0$ is a constant. It follows that

$$|P_\theta^k g(x,k) - P_{\theta'}^k g(x,k)| \leq cK\|g\|_U\|\theta - \theta'\|. \tag{4.79}$$

Since $U \equiv 1$ this demonstrates conditions 1. and 2. of (DRI3) and concludes the verification of assumptions.

Having verified all the assumptions necessary to apply Theorem 5.5 in [1], we can now state our convergence theorem. As remarked at the beginning of Section 4, we have restricted $\Theta$ to be a (arbitrarily large) compact set, and in Section 4.2 we have shown that $\theta^\star$ is attracting on all of $\Theta$. We do not use the full strength of the result in [1], which allows the iterates $\theta^t$ to escape arbitrarily large compact sets. Instead, we state the result in a way that assumes boundedness of the iterates and deduces convergence (see, for instance, [3], Corollary 16 in Section 1.8, or [11]).

**Theorem 2** *Suppose $\mathcal{T}_k(x'|x)dx$ defines a Feller transition kernel with unique stationary distribution $\rho_k(x)dx$ for each $k \in [K]$. Suppose also that $\pi : \Theta \to \mathcal{P}^\varepsilon(K)$ and $\gamma : \mathbb{R}^{KM} \to \mathcal{P}^\varepsilon(K)$ are continuously differentiable, and that $\mathbf{1} = (1,\ldots,1) \in \mathbb{R}^K$ is a globally asymptotically stable equilibrium for the differential equation $\dot{o}_k = \pi(F^\star,\mu^\star,o)/\gamma_k(\mu^\star) - o_k$. If the iterates $\theta^t$ are almost surely bounded, then, with probability 1, the iterates $\theta^t = (F^t,\mu^t,o^t)$ defined by (1.8), (1.9) and (1.4) satisfy $\lim_{t\to\infty} d(\theta^t,\{\theta^\star\}) = 0$, and the corresponding process $(\mathcal{X}^t, \mathcal{K}^t)$ with transition kernel $P_{\theta^t} = P_{\theta^t}^k \circ P_{\theta^t}^x$ converges in distribution to $\gamma_k \rho_k(x)dx\mathcal{C}(dk)$.*

**Proof** Proposition 6.1 in [1] shows that the drift conditions (DRI) verified above ensure that assumptions (A2) and (A3) in the same paper are satisfied, along with assumptions (A1) and (A4), verified above in Sections 4.3 and 4.2. We therefore fulfill the conditions of Theorem 5.5 in [1] on the compact sets $\mathcal{X} \times [K]$ and $\Theta$, which establishes $\theta^t \to \theta^\star$ for a suitably re-projected algorithm which is distinct from our algorithm only in that it may modify the process once it exits $\Theta$. By assumption, $\theta^t \in \Theta$ almost surely, so $\theta \to \theta^\star$ w.p.1.

Since $\lim_{t\to\infty} o_k^t = 1$ with probability 1 for each $k \in [K]$, we have $\lim_{t\to\infty} \pi(\theta^t) = \lim_{t\to\infty} \gamma(\mu^t) = \gamma(\mu^\star)$. The transition kernel $P_{\theta^\star}$ is assumed to have $\pi_k \rho_k(x)dx\mathcal{C}(dk)$ as its unique invariant density, and so $(\mathcal{X}^t, \mathcal{K}^t) \Rightarrow \gamma_k \rho_k(x)dx\mathcal{C}(dk)$, as required. □

Theorem 2 can be applied to the TSS algorithm. Let $\pi^{\text{TSS}}$ be the $\varepsilon_\pi$-regularized visit controlled rung distribution,

$$\pi^{\text{TSS}}(F,\mu,o) = (1-\varepsilon_\pi)\pi_k(\mu,o) + \varepsilon_\pi \gamma_k(\mu), \tag{4.80}$$

where

$$\pi_k(F,\mu,o) = \frac{\gamma_k(\mu)o_k^{-\eta}}{\sum_{\ell \in [K]} \gamma_\ell(\mu)o_\ell^{-\eta}}. \tag{4.81}$$

From Lemma 6, when $(F,\mu) = (F^\star,\mu^\star)$, the system (4.9) with $\pi = \pi^{\text{TSS}}$ has $\mathbf{1}$ as its global asymptotically stable equilibrium, and $\pi^{\text{TSS}}$ is continuously differentiable so long as $\gamma$ is. Moreover, $\varepsilon_\pi$ can be chosen such that $\pi^{\text{TSS}}(F,\mu,o) \in \mathcal{P}^\varepsilon(K)$ for all $(F,\mu,o) \in \Theta$. Propositions 1, 2 and 3 establish the existence of a global Lyapunov function, while Sections 4.3 and 4.4 verify Assumption (A4) and the drift conditions. Under the conditions of Theorem 2, the TSS iterates $\theta^t = (F^t,\mu^t,o^t)$ defined through (4.5), (4.6), (4.7) converge to $(F^\star,\mu^\star,\mathbf{1})$.

# Part II
# Implementation

## 5 Stochastic approximation with history forgetting and windows

This section gives the mathematical details for the stochastic approximation procedure modified for use with the windowing system (Section 2.1 in the main text) and history forgetting (Section 2.2). To assist the reader interested in the implementation of the algorithm, we relate the equations to pieces of the code, the source code of which is available at `https://github.com/DEShawResearch/tss`.

**Notation.** The position $x \in \mathcal{S}$, rung $k \in [K]$ and window $j \in [J]$ of replica $r = 1,\ldots,R$ at time $t$ will be written as $\mathcal{X}^t(r)$, $\mathcal{K}^t(r)$ and $\mathcal{J}^t(r)$, respectively.

### 5.1 Per-window estimates

The recursions written here are for a single window $W_j$ and multiple replicas $r = 1,\ldots,R$. Define also the times $\tau_l$ introduced in the main text, Section 2.2. Given $\phi > 1$, we consider the exponentially increasing sequence of integers defined recursively by $\tau_0 = 0$, $\tau_1 = 1$ and

$$\tau_{l+1} = \lceil \phi \tau_l \rceil, l \geqslant 1. \tag{5.1}$$

Guidance on choosing $\phi > 1$ is provided in Section 2.4 of the main text.

#### 5.1.1 Epoch-based counting processes

Let $\alpha \in [0,1)$ be the fraction of data that is dropped asymptotically, so that throughout the simulation, the fraction of stored samples is approximately $1 - \alpha$. The "recent history" is the data generated between time $\alpha t$ and $t$. Recall from the main text that at time $t$ we only update epoch $n(t)$, and potentially drop or start a new epoch; all the intermediate (i.e., not the first nor the last) epochs are untouched.

The number of visits made to the window $W_j$ by replica $r$ in epoch $l$ is

$$\mathcal{N}_j^{t,l}(r) = \sum_{s=\tau_{l-1}+1}^{\tau_l \wedge t} \mathbf{1}_{\{j\}}(\mathcal{J}^s(r)). \tag{5.2}$$

For each replica $r \in [R]$, the number of visits paid to the window $W_j$ in the recent history of the simulation is

$$\mathcal{N}_j^t(r) = \sum_{l=n(\alpha t)}^{n(t)} \mathcal{N}_j^{t,l}(r) \tag{5.3}$$

which can be further combined over all replicas

$$\mathcal{N}_j^t = \sum_{r=1}^{R} \mathcal{N}_j^t(r) \tag{5.4}$$

Only the epochs $l$ with $n(\alpha t) \leqslant l \leqslant n(t)$ are needed at the current time $t$. Notice that the current epoch is typically incomplete, so the population counters with $l = n(t)$ change with time. During any cycle $t$ that marks the beginning of a new epoch, new counters $\mathcal{N}_j^{t,n(t)}$ are instantiated on the store and initialized to 0. During any cycle $t$ for which $\alpha t$ crosses into a new epoch, the chronologically earliest counters $\mathcal{N}_j^{t,l}(r)$, the ones with $l = n(\alpha t) - 1$, are removed from the store. Within a cycle, both actions happen immediately after the execution of the Monte Carlo moves in the cycle. The number of counters committed to the store remains bounded during the simulation. We note that, by its very definition, $\mathcal{N}_j^{t,l}(r)$ needs to be incremented at the current time $t$ only for $l = n(t)$ and only for those $j \in \cup_{r=1}^{R} \{\mathcal{J}^t(r)\}$.

### 5.1.2 Epoch-based free energy estimates

For each $j \in [J]$ and $k \in W_j$ and each replica $r = 1, \ldots, R$, define the importance sampling ratio

$$\mathcal{R}_{j;k}^t(r) = 1_{\{j\}}(\mathcal{J}^t) \frac{e^{-H_k(X^t(r))}}{d_j^t(r)}, \tag{5.5}$$

where

$$d_j^t(r) = \begin{cases} \sum_{k \in W_j} \pi_{j;k}^{t-1} e^{F_{j;k}^{t-1} - H_k(X^t(r))}, & \text{if } \mathcal{N}_j^t \geqslant 1; \\ 1 & \text{else.} \end{cases} \tag{5.6}$$

The weights $\pi_{j;k}^{t-1}$ satisfy $\sum_{k \in W_j} \pi_{j;k}^{t-1} = 1$ and will be defined in the next section, whereas the overall estimates $F_{j;k}^{t-1}$ will be defined in Section 5.2.

For all $n(\alpha t) \leqslant l \leqslant n(t)$, $k \in [K]$ and $j \in \text{win}(k)$, we compute and store the epoch free energy estimates

$$e^{-\mathcal{F}_{j;k}^{t,l}} = \frac{1}{\mathcal{N}_j^{t,l}} \sum_{s=\tau_{l-1}+1}^{\tau_l \wedge t} \sum_{r=1}^{R} \mathcal{R}_{j;k}^s(r), \tag{5.7}$$

where we make the convenient but mathematically inconsequential convention that $\mathcal{F}_{j;k}^{t,l} = 0$ if $\mathcal{N}_j^{t,l} = 0$. At cycle $t$, all epoch free energies remain unchanged, except for those with $l = n(t)$. The **per-window, per-epoch** free energies for $l = n(t)$ are computed with the following recurrence: for $j = \text{win}(\mathcal{K}^t)$ and for all $k \in W_j$,

$$\mathcal{F}_{j;k}^{t,l} = \mathcal{F}_{j;k}^{t-1,l} - \log\left(1 + \frac{1}{\mathcal{N}_j^{t,l}} \left[\mathcal{N}_j^{t-1,l} - \mathcal{N}_n^{t,l} + e^{\mathcal{F}_{j;k}^{t-1,l}} \sum_{r=1}^{R} \mathcal{R}_{j;k}^t(r)\right]\right). \tag{5.8}$$

Equation (5.8) corresponds to `window_t::update_epoch_fe_estimates` in `window.cxx` in the source code in the GitHub repository.

We note that the convention $\mathcal{F}_{j;k}^{t,l} = 0$ if $\mathcal{N}_j^{t,l} = 0$ is indeed inconsequential, because the otherwise undefined $\mathcal{F}_{j;k}^{t-1,l}$ would have canceled out anyway on the right-hand side of equation (5.8). The **per-window, all-epoch** free energy estimates $F_{j;k}^t$ are defined by

$$e^{-F_{j;k}^t} = \sum_{l=n(\alpha t)}^{n(t)} \frac{\mathcal{N}_j^{t,l}}{\sum_{l'=n(\alpha t)}^{n(t)} \mathcal{N}_j^{t,l'}} e^{-\mathcal{F}_{j;k}^{t,l}} \tag{5.9}$$

for all $k \in [K]$ and $j \in \text{win}(k)$. Equation (5.9) corresponds to `window_t::estimate_all_epochs` in `window.cxx` in the source code.

### 5.1.3 Epoch-based rung distribution estimates

Much of the computation of the asymptotic rung distribution is implemented in `covariance_estimator.{cxx,hxx}`. The default TSS rung distribution is the coordinate-invariant one, which we recall is written in terms of the metric

$$g_{ij}(\lambda) = \mathbb{E}\left[\frac{\partial H_\lambda}{\partial \lambda_i}(X) \cdot \frac{\partial H_\lambda}{\partial \lambda_j}(X)\right] - \mathbb{E}\left[\frac{\partial H_\lambda}{\partial \lambda_i}\right] \mathbb{E}\left[\frac{\partial H_\lambda}{\partial \lambda_j}\right], \tag{5.10}$$

where the expectation is taken with respect to $p_\lambda(x)\mathrm{d}x$. An unregularized form of the asymptotic TSS metric is given by

$$\frac{\sqrt{\det(g(\lambda_k))}\mathrm{vol}(\lambda_k)}{\sum_{j \in [K]} \sqrt{\det(g(\lambda_j))}\mathrm{vol}(\lambda_j)}, \tag{5.11}$$

where $\mathrm{vol}(\lambda_k)$ is a volume element arising from the discretization of the space $\Lambda$ into the set $\{\lambda_1, \ldots, \lambda_K\}$. For instance, if $\Lambda = [0,1]$, then we might choose the discretization $\{[k/K, (k+1)/K), k = 0, \ldots, K-2\} \cup \{[1 - 1/K, 1]\}$, and then $\mathrm{vol}(\lambda_k) = 1/K$ for all $k \in [K]$.

We will let $G_\lambda$ denote the matrix whose $i,j$ entry is $g_{ij}(\lambda)$, and $G_k$ will denote $G_{\lambda_k}$. We can evaluate the partial derivatives using a second-order central difference,

$$\frac{\partial H_k}{\partial \lambda_i}(x) \approx \frac{H_{k+1}(x) - H_{k-1}(x)}{2}. \tag{5.12}$$

The central difference may be undefined for some of the rungs $k \in [K]$ such that $\lambda_k$ lies on or near the boundary of $\Lambda$. In that case, we can use a first-order forward or backward difference, as needed. It is important that the approximation is tied to the original graph edges and is independent of the choice of windows. To avoid clutter, we set

$$\psi_k(x) = \nabla_\lambda H(x)|_{\lambda = \lambda_k}, \tag{5.13}$$

and note for the sake of generality that the estimation method described below is applicable to an arbitrary function $\psi(x)$.

As we have done for free energies, for all $l$ such that $n(\alpha t) \leq l \leq n(t)$ and for all $k \in [K]$ and $j \in \text{win}(k)$, we introduce the epoch estimators $\Psi$ by way of the formula

$$\psi_{j;k}^{t,l} = \frac{1}{\mathcal{N}_j^{t,l}} \sum_{s=\tau_{l-1}+1}^{\tau_l \wedge t} \sum_{r=1}^{R} e^{F_{j;k}^{t,l}} \mathcal{R}_{j;k}^s(r) \psi_k(\mathcal{X}^s), \tag{5.14}$$

where the equality in (5.14) is meant componentwise whenever $\psi_{j;k}^{t,l}$ is a vector. As before, we set $\psi_{j;k}^{t,l} = 0$ if $\mathcal{N}_j^{t,l} = 0$. The definition of the epoch estimators is such that the epoch estimators $\psi_{j;k}^{t,l}$ of a constant function $\psi_k(x)$ is the constant itself. The overall estimates are obtained by weighting the epoch estimates according to

$$\psi_{j;k}^{t} = \sum_{l=n(\alpha t)}^{n(t)} \frac{\mathcal{N}_j^{t,l} e^{-F_{j;k}^{t,l}}}{\sum_{l'=n(\alpha t)}^{n(t)} \mathcal{N}_j^{t,l'} e^{-F_{j;k}^{t,l'}}} \psi_{j;k}^{t,l}. \tag{5.15}$$

The estimates are left undefined if the denominator of (5.15) vanishes. As for the free energies, at the current time $t$ the epoch estimates need to be updated only for the last epoch, the one with $l = n(t)$, as follows: for all $j \in \mathrm{win}(\mathcal{K}^t)$ and $k \in W_m$,

$$\psi_{j;k}^{t,l} = e^{F_{j;k}^{t,l} - F_{j;k}^{t-1,l}}$$
$$\times \left( \psi_{j;k}^{t-1,l} + \frac{1}{\mathcal{N}_j^{t,l}} \left[ (\mathcal{N}_j^{t,l} - \mathcal{N}_j^{t-1,l}) \psi_{j;k}^{t-1,l} + e^{F_{j;k}^{t-1,l}} \sum_{r=1}^R \mathcal{R}_{j;k}^t(r) \psi_k(X^t(r)) \right] \right). \quad (5.16)$$

For each $j \in [J]$ and $k : \lambda_k \in W_j$, the estimates for the covariance matrices $G_{j;k}^t$ at step $t$ are given by the formula

$$G_{j;k}^t = (\psi_1)_{j;k}^t - (\psi_2)_{j;k}^t ((\psi_2)_{j;k}^t)^\top \quad (5.17)$$

where $(\psi_1)_{j;k}^t$, $(\psi_2)_{j;k}^t$ are the overall averages of $\psi_2(x) = \nabla_\lambda H_\lambda(x)|_{\lambda=k}$ and $\psi_1(x) = \psi_2(x)\psi_2(x)^\top$. Equation (5.17) is implemented in `covariance_estimator_t::update_epoch_estimates` in `covariance_estimator.cxx`.

The coordinate invariant distribution is regularized for any fixed $\varepsilon_\gamma \in (0, 1]$ as

$$\gamma_{j;k}^t = \frac{[(1-\varepsilon_\gamma)\det(G_{j;k}^t)^{1/2} + \varepsilon_\gamma B^t]\mathrm{vol}(\lambda_k)}{\sum_{k' \in W_j} [(1-\varepsilon_\gamma)\det(G_{j;k'}^t)^{1/2} + \varepsilon_\gamma B^t]\mathrm{vol}(\lambda_{k'})}, \quad (5.18)$$

where

$$B^t = \max_{m=1,\ldots,M} \max_{k=1,\ldots,K} \det(G_{j;k}^t)^{1/2}. \quad (5.19)$$

The regularized per-window rung distribution (5.18) is implemented in `window_t::normalized_weights` in `window.cxx`.

### 5.1.4 Epoch-based tilts

For all $n(\alpha t) \le l \le n(t)$, $k \in [K]$ and $j \in [J]$, we compute and store

$$o_{j;k}^{t,l} = \frac{1}{\mathcal{N}_j^{t,l}} \sum_{s=\tau_{l-1}+1}^{\tau_l \wedge t} \sum_{r=1}^R \mathbf{1}_{\{j\}}(\mathcal{J}^t) \frac{\mathbf{1}_{\{k\}}(\mathcal{K}^t)}{\gamma_{j;k}^{-1}}. \quad (5.20)$$

We compute these estimates through the **per-window, per-epoch** recursion

$$o_{j;k}^{t,l} = o_{j;k}^{t-1,l} + \frac{1}{\mathcal{N}_j^{t,l}} \left[ (\mathcal{N}_j^{t-1,l} - \mathcal{N}_j^{t,l}) o_{j;k}^{t-1,l} + \sum_{r=1}^R \mathbf{1}_{\{j\}}(\mathcal{J}^t) \frac{\mathbf{1}_{\{k\}}(\mathcal{K}^t)}{\gamma_{j;k}^t} \right]. \quad (5.21)$$

Equation (5.21) corresponds to `window_t::update_tilts` in `window.cxx` in the code. We then provide combine the epoch estimates to obtain **per-window, all-epoch** tilt estimates $o_{j;k}^t$ through

$$o_{j;k}^t = \sum_{l=n(\alpha t)}^{n(t)} \frac{\mathcal{N}_j^{t,l}}{\mathcal{N}_j^t} o_{j;k}^{t,l}. \quad (5.22)$$

Equation (5.22) corresponds to `window_t::overall_tilt` in `window.cxx` in the code.

## 5.2 Global estimates

### 5.2.1 Global rung weights

The C++ functions referred to in this section rely heavily on ideas and notation introduced in Section 2.1 and 2.2 of the main text.

As a first step, we must estimate the window marginal probability $\bar{p}_j = \bar{P}_{\pi,Z}(\mathcal{J} = j)$. Our estimate $p_j$ of $\bar{p}_j$ is obtained by solving the eigenvector problem

$$Qp = p, \tag{5.23}$$

where $Q = (q_{ij})$, $1 \leqslant i,j \leqslant J$,

$$q_{ij} = \frac{1}{2} \frac{\sum_{k \in W_i \cap W_j} \gamma_{j;k} o_{j;k}}{\sum_{k:\lambda_k \in W_j} \gamma_{j;k} o_{j;k}}. \tag{5.24}$$

We note that early in the simulation, not all entries of the matrix are available, since only a limited number of windows $W_j$ have been visited and thus have $\mathcal{N}_j^t \geqslant 1$. We restrict the eigenvector problem only to the visited windows, and set $p_j^t = 0$ (i.e., the estimate $p_j$ at time $t$) if $\mathcal{N}_j^t = 0$. A column $j$ of $Q^t$ (the matrix $Q$ at time $t$) such that $W_j$ has been visited at least once may have non-zero components $q_{ij}^t$ corresponding to non-visited windows $W_i$. Those components are summed up, with the sum added to the $(j,j)$ diagonal entry. The non-visited components are then zeroed. We obtain a lower-dimensional left-stochastic matrix acting on the space of probability vectors defined over the space of visited states. We can use QR to find the invariant stochastic eigenvector of minimum $L^2$ norm, in the rare cases where the invariant stochastic vector isn't already unique. This function is implemented in `free_energy_estimator_t::window_weights_qr`.

Once the marginal window probabilities at time $t$, $p_j^t$, have been obtained through (5.23), we can compute the probability $q_k^t$ of visiting rung $k \in [K]$ irrespective of the window as follows:

$$q_k^t = \sum_{j \in \mathrm{win}(k)} p_j^t \frac{\gamma_{j;k} o_{j;k}^t}{\sum_{\ell \in W_j} \gamma_{j;\ell} o_{j;\ell}^t}. \tag{5.25}$$

Equation (5.25) corresponds to `free_energy_estimator_t::global_weights` in `free_energy_estimator.cxx` in the code.

### 5.2.2 Global visit control free energies

To determine the global visit control free energies, we must identify the offsets $f_j^t$ introduced at the end of Section 2.2.1 of the main text. These offsets solve the following system of equations, for $j \in [J]$ and a visit control parameter $\eta > 0$

$$f_j^t = (\eta+1) \log \left( \sum_{k:\lambda_k \in W_j} q_k^t \frac{\gamma_{j;k} e^{\frac{1}{\eta+1} F_{j;k}^t}}{\sum_{i \in \mathrm{win}(k)} p_i^t \gamma_{i;k} e^{\frac{1}{\eta+1}(F_{i;k}^t - f_i^t)}} \right). \tag{5.26}$$

Equivalently, these unknown offsets $f_j^t$ minimize the objective function

$$\frac{1}{\eta+1} \sum_{j \in [J]} p_j^t f_j^t + \sum_{k \in [K]} q_k^t \log \left( \sum_{j \in \mathrm{win}(k)} p_j^t \gamma_{j;k} e^{\frac{1}{\eta+1}(F_{j;k}^t - f_j^t)} \right), \tag{5.27}$$

which is a convex function in $f_j^t$.

There are multiple avenues for solving (5.26) or, equivalently, minimizing (5.27). We employ perhaps the simplest approach, which is to perform fixed point iteration, represented succinctly as

$$f = h(f), \tag{5.28}$$

where $f = (f_1, \ldots, f_J)$ and $h(f) = (h_1(f), \ldots, h_J(f))$ is defined componentwise by

$$h_j(f) = (\eta + 1) \log(g_j(f)), \quad j \in [J], \tag{5.29}$$

with

$$g_j(f) = \sum_{k \in W_j} q_k^t \frac{\gamma_{j;k}^t e^{\frac{1}{\eta+1} F_{j;k}^t}}{\sum_{i \in \text{win}(k)} p_i^t \gamma_{j;k}^t e^{\frac{1}{\eta+1}(F_{i;k}^t - f_i^t)}}. \tag{5.30}$$

This fixed point iteration is performed in `free_energy_estimator_t::improve_fes`.

### 5.2.3 Global reported free energies

The visit control free energies already define global estimators for the free energies, but these are $\eta$-dependent. Because they are functions of the empirical frequencies of visiting windows and of the tilts $o^t = \{o_{j;k}^t; j \in [J], k : \lambda_k \in W_j\}$, these estimators exhibit large error bars, which is typical of event counting. As described in Section 2.2.2 of the main text, we can provide global estimators with lower error bars by letting $\eta \to \infty$; we call these the reported TSS free energies, and which we denote by $F_k^{\text{TSS}}$

The limit $\eta \to \infty$ forces us to take $o_{j;k}^t = 1$ for all $j \in [J], k \in [K]$. Consequently, we can identify the global rung probabilities as in (5.25), this time with $o_{j;k}^t$ replaced by 1. We define the reported rung density to be

$$\gamma_k^{\text{TSS}} = \sum_{j \in [J]} p_j \gamma_{j;k}, \tag{5.31}$$

for each $k \in [K]$.

As discussed in the main text, in the limit $\eta \to \infty$, we can compute the gradient of the objective (5.27) and set the limiting gradients to zero. The system of equations so-obtained is

$$\sum_{j \in [J]} (\delta_{ij} - t_{ij}) f_j^{\text{TSS}} = \sum_{k \in W_i} \gamma_{i;k}^t \left( F_{i;k}^t - \sum_{j \in \text{win}(k)} \frac{p_j^t \gamma_{j;k}^t}{\gamma_k^{\text{TSS}}} F_{j;k}^t \right), \tag{5.32}$$

where

$$t_{ij} = \sum_{k \in W_i \cap W_j} \gamma_{i;k}^t \frac{p_j^t \gamma_{j;k}^t}{\gamma_k^{\text{TSS}}}. \tag{5.33}$$

Once the offsets $f_j^{\text{TSS}}$ are obtained, we can compute the reported TSS free energies for each $k \in [K]$ through

$$F_k^{\text{TSS}} = \frac{1}{\gamma_k^{\text{TSS}}} \sum_{j \in \text{win}(k)} p_j \gamma_{j;k} (F_{j;k} - f_j^{\text{TSS}}). \tag{5.34}$$

This procedure is carried out in `free_energy_estimator_t::solve_fes_infinite_eta`.

# 6 Jackknife error estimation

The task of computing error bars is made particularly easy by the fact that we compute and store epoch estimates. Thus, estimates for the errors are obtainable, for example, by way of the jackknife method. The $m$-th jackknife replicate $F_{j;k}^{t,(m)}$ is defined by

$$e^{-F_{j;k}^{t,(m)}} = \sum_{l=n(\alpha t), l \neq m}^{n(t)} \frac{\mathcal{N}_{m;k}^{t,l}}{\sum_{l'=n(\alpha t), l' \neq m}^{n(t)} \mathcal{N}_{j;k}^{t,l'}} e^{-F_{j;k}^{t,l}}, \qquad (6.1)$$

which is equation (5.9) with the $m$-th epoch deleted. The use of replicates for the evaluation of error bars will be presented only for free energy differences, since it proceeds analogously for the other estimated quantities or functions of them. For instance, the $m$-th jackknife replicate $\langle \psi \rangle_{j;k}^{t,(m)}$ is obtained by deleting the $m$-th epoch in equation (5.15).

We assume that the TSS simulation is sufficiently advanced that all windows have been visited at least once, even when any one of the epochs is deleted. Otherwise, the analysis must be restricted to the set of all windows that have been so visited. Global reported free energies, as described in Section 2.2.2 of the main text, can then be computed for each replicate $m$. The global $m$-th jackknife replicates are denoted by $F_k^{t,(m)}$.

For $l \in n(\alpha t), \ldots, n(t)$, the numbers

$$a_l = \frac{\tau_l \wedge t - \tau_{l-1}}{t - \tau_{n(\alpha t)-1}} \qquad (6.2)$$

are positive and sum to 1. Asymptotically, the error for $F_{e',k'}^t - F_{e,k}^t$ is normally distributed about a small bias, with the MSE estimated by

$$\widehat{\mathrm{MSE}}(t; k, k') = \frac{1}{n(t) - n(\alpha t)} \sum_{l=n(\alpha t)}^{n(t)} \frac{(1-a_l)^2}{a_l} \left[ \left( F_{e',k'}^{t,(l)} - F_{e,k}^{t,(l)} \right) - \left( F_{e',k'}^t - F_{e,k}^t \right) \right]^2. \qquad (6.3)$$

This MSE estimator is implemented in `free_energy_estimator_t::update_errors` in `free_energy_estimator.cxx`.

In our analysis, we do not use the debiased jackknife estimator. We assume the simulation is sufficiently well converged that the main source of error is the noise term, the variance of which is approximated by (6.3). The coefficients $(1-a_l)^2/a_l$ have been chosen such that the estimator $\widehat{\mathrm{MSE}}(t)$ has a chi-squared distribution in the limit $t \to \infty$ (assuming the central limit theorem holds).

Standard precautions apply for the estimation of the MSE to be accurate. The epochs must be long enough that the jackknife pseudo-values

$$a_l^{-1} \left[ \left( F_{e',k'}^t - F_{e,k}^t \right) - (1-a_l) \left( F_{e',k'}^{t,(l)} - F_{e,k}^{t,(l)} \right) \right] \qquad (6.4)$$

can be assumed to be independent normal variables. At the same time, the number of epochs must be large enough that the MSE estimator defined by (6.3) has sufficiently small noise to approximate well the targeted MSE. For sufficiently large $t$, the instantaneous number of epochs $n(t) - n(\alpha t) + 1$ fluctuates between two values $n_{\mathrm{epochs}}$ and $n_{\mathrm{epochs}} + 1$, as new epochs are added and old epochs are eliminated. We recommend that the number of epochs be at least 32. The proper number of epochs can be dialed in by setting $\phi$ according to the formula

$$\phi = \alpha^{-1/n_{\mathrm{epochs}}}. \qquad (6.5)$$

As an example, for $\alpha = 0.19$ and $n_{\mathrm{epochs}} = 32$, we compute $\phi = 0.19^{-1/32} \approx 1.0533$.

# Part III
# Numerics

## 7 Numerical example: calculation of free energy differences in aqueous solution

### 7.1 Simulation details

The molecular dynamics time step was 2 fs. The system was run in the NPT ensemble (constant number of particles, constant pressure, and constant temperature). The barostat and thermostat were applied using the multigrator [12]. The MTK barostat [14] was used, and the Langevin thermostat was used.

For the integration, a reversible reference system propagator algorithm (RESPA) scheme [24] was employed, with an inner time step of 2 fs and an outer time step (for the far electrostatics) of 6 fs.

The volume of the cell used for the simulation was approximately 51.5 Å × 51.5 Å × 51.5 Å throughout the simulation. The reference temperature was 298.15 K, and the reference pressure was 1.01325 bar.

We used a softcore potential, described in Section 7.2. For the electrostatic calculations, we used the $u$-series [18] with a cutoff radius of 12 Å.

### 7.2 Softcore potential

The Lennard-Jones 6-12 potential is

$$V_{\text{LJ}}(r) = 4\varepsilon_{\text{LJ}} \phi((r/\sigma_{\text{LJ}})^6), \tag{7.1}$$

where $\phi(z) = z^{-1}(z^{-1} - 1)$. The function $\phi$ has derivative

$$\phi'(z) = z^{-2}(-2z^{-1} + 1). \tag{7.2}$$

Let $r_s \in (0, \sigma_{\text{LJ}})$ be a softening ratio and $z_s = (r_s/\sigma_{\text{LJ}})^6$ the corresponding softening $z$-point location in $(0, 1)$. We replace $\phi(z)$ with the softened version

$$\tilde{\phi}(z) = \begin{cases} \phi(z_s) + \phi'(z_s)(z - z_s) & \text{if } z \leqslant z_s, \\ \phi(z) & \text{if } z > z_s. \end{cases} \tag{7.3}$$

The softened *and* tapered interparticle potential is

$$V_{ip}(r) = \begin{cases} 4\varepsilon_{\text{LJ}} \tilde{\phi}((r/\sigma_{\text{LJ}})^6) - V_{\text{LJ}}(r_c) & \text{if } r < r_c, \\ 0 & \text{otherwise,} \end{cases} \tag{7.4}$$

where $r_c > r_s$ is the cutoff radius.

We can make this parametric potential into a softcore potential by multiplying it by $\lambda \in [0, 1]$ and setting $z_s$ depending on $\lambda$. While the dependence of $z_s$ on $\lambda$ can be optimized, for the purposes of our simulation in Section 3 of the main text it is sufficient to use

$$z_s(\lambda) = \lambda b_{\min} + (1 - \lambda) b_{\max}, \tag{7.5}$$

with $0 \leqslant b_{\min} \leqslant b_{\max} \leqslant 1$. The value $b_{\min}$ is normally set to 0, but it can be set larger if only a portion of the range of the LJ potential is machine representable. The value $b_{\max}$ is determined experimentally, to minimize the error for free energy differences between two ends with the LJ interactions turned on and off, respectively. In our simulations, we used $b_{\max} = 1.0$.

## 7.3 Figures

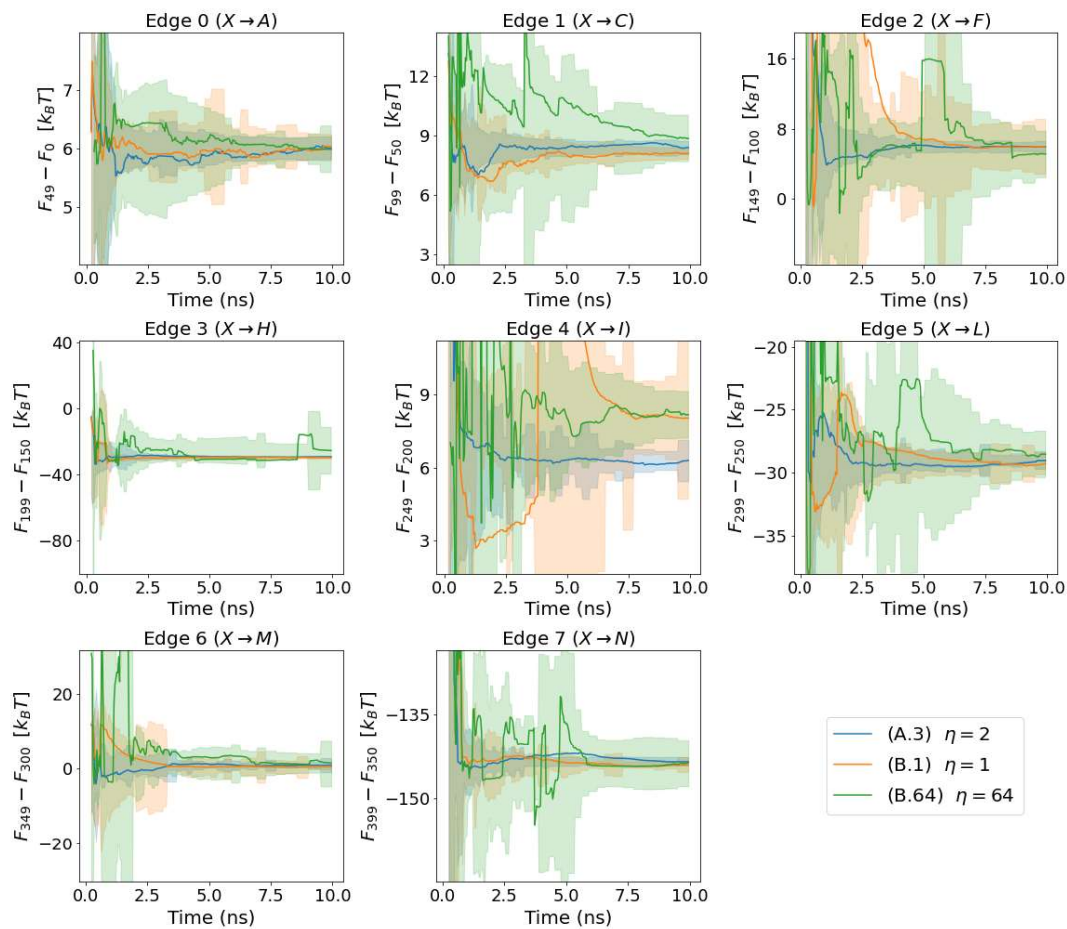

**Figure S1.a** Visit control accelerates convergence to a neighborhood of the equilibrium values early in the simulation. Plotted are the first 10 ns of simulations using η = 1, 2, and 64 ((B.1), (B.2), (B.64), respectively, with (B.2) and (A.3) referring to the same simulation). The free energy difference estimates are shown for each of the eight edges, along with their estimated error bars. The η = 1 simulations suffer from the exponential slowdown described in Section 1.4, an effect which is particularly visible in Edges 2, 4 and 5. On the other hand, the η = 64 estimates tend to oscillate around the η = 2 estimates, which is particularly visible in Edges 2, 3, 5, and 7.

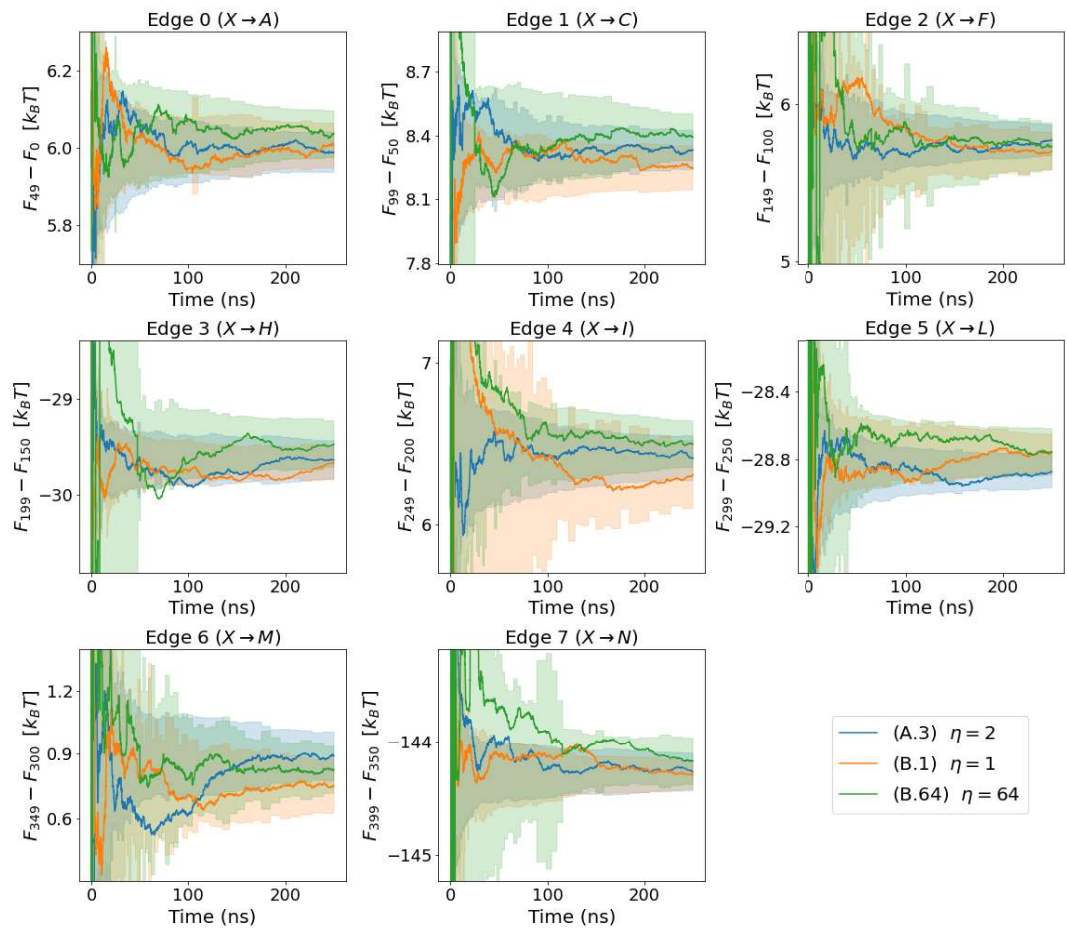

**Figure S1.b** Visit control promotes efficient sampling and estimation for a broad range of $\eta > 0$, with decreasing efficiency at the extremes (when $\eta$ is very small or very large). Plotted are simulations using $\eta = 1, 2,$ and 64 ((B.1), (B.2), (B.64), respectively). The free energy difference estimates are shown for each of the eight edges, along with their estimated error bars. The error bars for the simulations with $\eta = 1$ and $\eta = 64$ are both large early in the simulation, and the estimates themselves sometimes appear oscillatory, especially within the first 50 ns. Afterwards, all three estimates and their error bars become comparable. It is worth noting that even these choices of $\eta$ seem poor relative only to $\eta = 2$. The estimate for nearly every edge was within 1 $k_BT$ by the 10 ns mark (the exceptions being Edge 3 for $\eta = 64$ and Edge 4 for $\eta = 1$); both $\eta = 1$ and $\eta = 64$ fare far better than $\eta = 0$, which corresponds to disabling visit control (see Figure 8 in the main text).

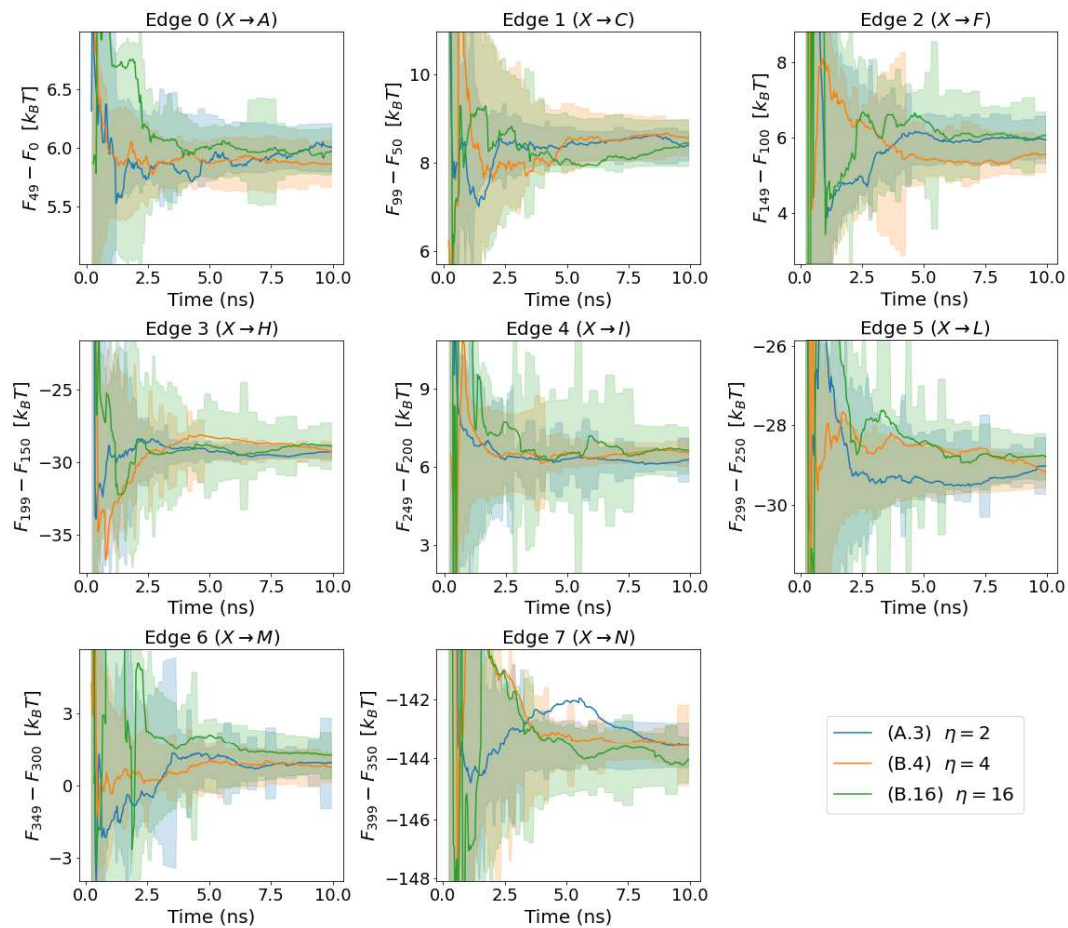

**Figure S2.a** Visit control comparably accelerates convergence for a broad range of $\eta > 0$. Plotted are simulations using $\eta = 2, 4$, and 16 ((B.2), (B.4), (B.16), respectively). The free energy difference estimates are shown for each of the eight edges, along with their estimated error bars. It is difficult to distinguish between the three values of $\eta$ after 5 ns, though we note that the most the most extreme value, $\eta = 16$, does appear to fluctuate more strongly prior to the 5 ns mark.

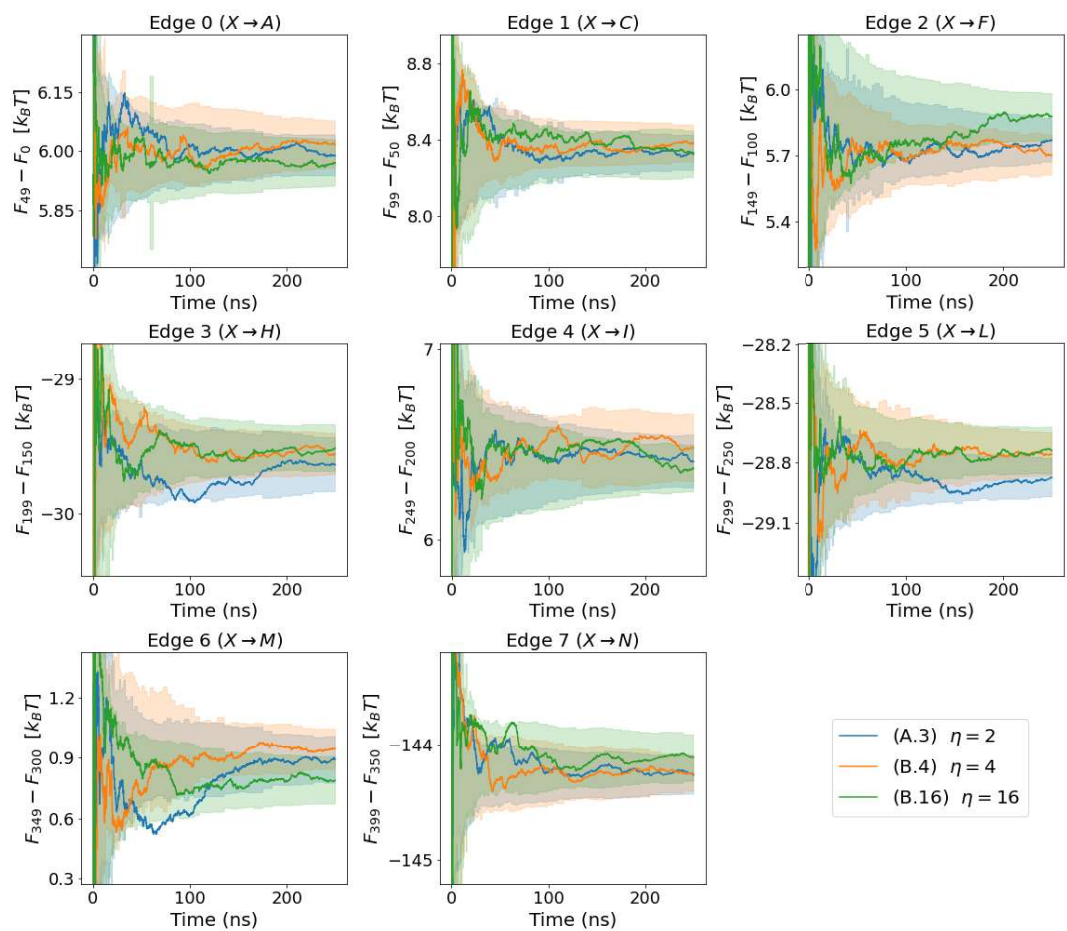

**Figure S2.b** Visit control comparably accelerates convergence for a broad range of $\eta > 0$. Plotted are simulations using $\eta = 2, 4$, and 16 ((B.2), (B.4), (B.16), respectively). The free energy difference estimates are shown for each of the eight edges, along with their estimated error bars. Though the estimates for $\eta = 16$ tend to fluctuate more early in the simulation, on longer timescales the three choices of $\eta$ are not distinguishable through their estimated error bars.

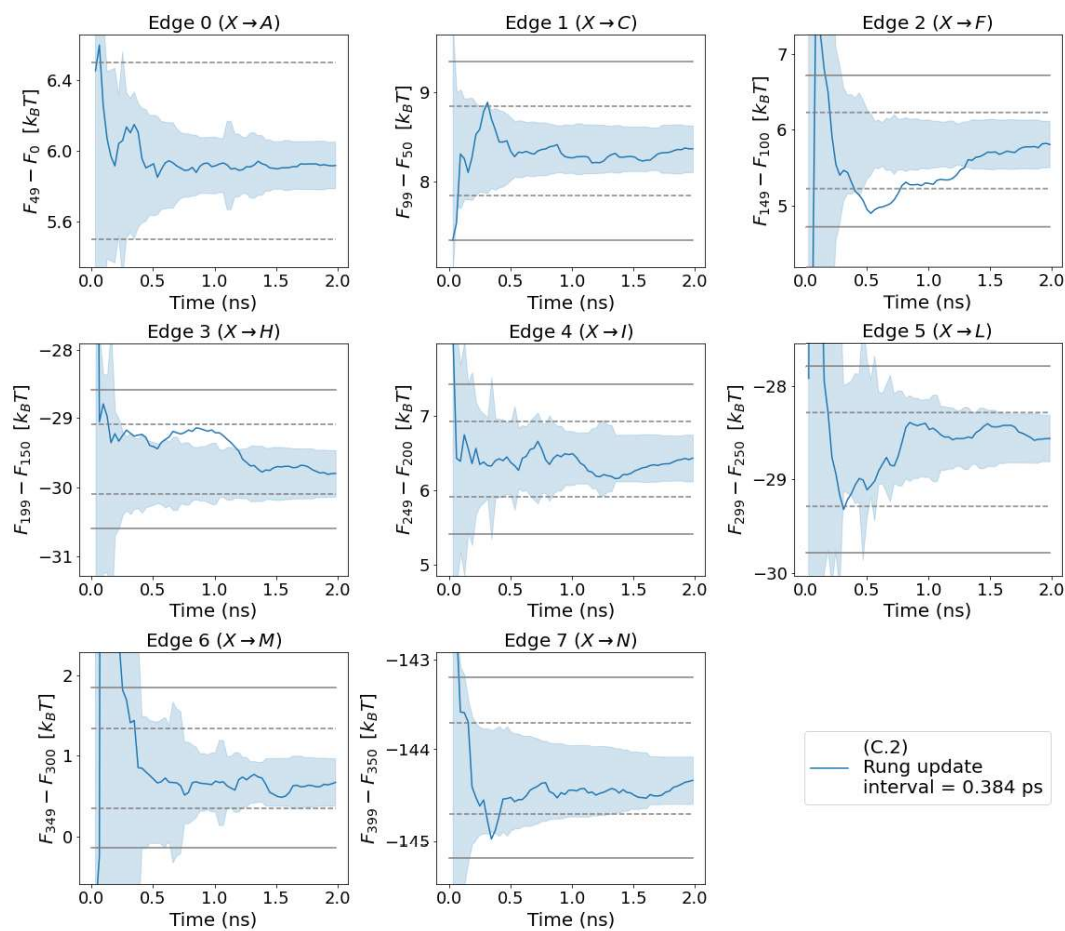

**Figure S3** When applied in practice, TSS provides useful estimates quickly. We performed long timescale simulations for the solvation free energy calculations in order to study the asymptotic behavior of the algorithm, but in practice, one would typically run much shorter simulations. TSS provides estimates that reach a small neighborhood of the equilibrium values quickly. The figure shows the results of simulation (C.2), using $\nu = 1$ and a rung update interval of 0.384 ps, over the course of the first 2 ns of simulation, with two such neighborhoods delineated. The first neighborhood, delineated by dashed lines, indicates the range within 0.5 $k_B T$ of the estimates at 250 ns; the second one, delineated by solid lines, indicates the range within 1 $k_B T$ of the estimates at 250 ns. For this simulation, all the estimates stay within the first neighborhood for the remainder of the simulation after the 1 ns mark, and within the second neighborhood after just 0.3 ns.

## 8 Gaussian model

We will study the model

$$H_\lambda(x) = \frac{1}{2}(x-\lambda)^2 \tag{8.1}$$

with $\Lambda = [0, L]$. We take $L \in \mathbb{Z}_{>0}$ and discretize $\Lambda$ as $\lambda_k = k$ for $k = 0, \ldots, K-1$, where $K = L+1$. For each $k \in [K]$, the free energy is

$$F_k^\star = -\log \int_\mathbb{R} e^{-H_k(x)} dx = -\log(\sqrt{2\pi}), \tag{8.2}$$

and so all the exact free energy differences $F_k^\star - F_j^\star$, $k, j \in [K]$ are equal to 0. We use a flat asymptotic distribution $\gamma_k = 1/(K-1)$ for $k = 1, \ldots, K-2$ and $\gamma_0 = \gamma_{K-1} = 1/(2(K-1))$.

This model is motivated by umbrella sampling, in which harmonic restraints are placed upon a so-called *collective variable* (also known as a *reaction coordinate*). The situation studied here is simpler, owing to the lack of any "orthogonal" degrees of freedom, and the independent sampling of $\mathcal{X}$. Nevertheless, it can provide valuable insight on the pre-asymptotic behavior of the estimates, and the scaling of the asymptotic variance with respect to the parameter $L$.

We make a few remarks about the simulations below. The plots shown below are from subsampled trajectories (i.e., not all frames are stored to disk). When error bars are displayed, they are centered around the true asymptotic free energy difference ($F_{K-1}^\star - F_0^\star = 0$), so we plot $\pm 2\sqrt{\mathrm{MSE}(t; 0, K-1)}$, where MSE is defined in Section 6, equations (6.3). All the examples in this section are available in a Python notebook on Github, `https://github.com/DEShawResearch/tss`.

### 8.1 SAMS

We first study TSS without visit control ($\eta = 0$), without history forgetting ($\alpha = 0$), with a single window, and with a single move per estimator update ($\nu = 0$). In such a case, TSS reduces to SAMS [22], with the difference that TSS uses a logarithmic rather than a linear form of the estimator update; but this difference is minor, and the algorithms are asymptotically equivalent. We ran two simulations in these conditions: one with $L = 7$ and $K = 8$ rungs, and another with $L = 15$ and 16 rungs.

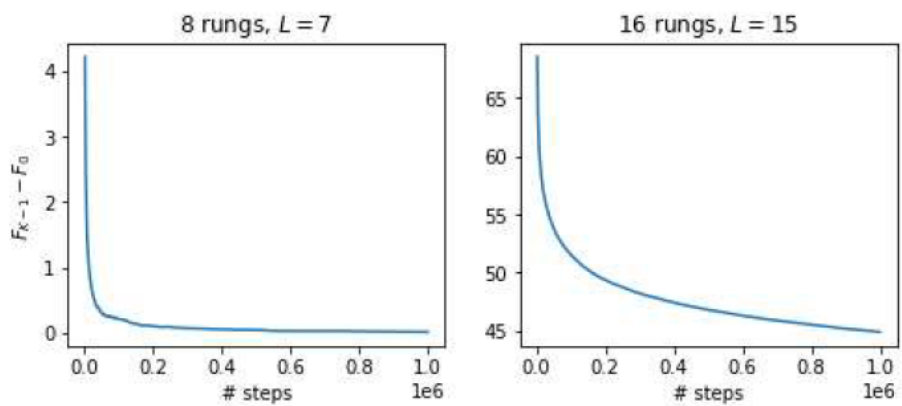

Figure 8.1.1: End-to-end free energy difference estimates for $L = 7, 15$.

While the estimator performs adequately for $L = 7$, for $L = 15$ it suffers from dramatically from the exponential slowdown alluded to in Sections 1.4 and 2.3 of the main text. A partial remedy could be sought by changing the gain in the stochastic approximation recursions from $O(t^{-1})$ to $O(t^{-\beta})$, for some fixed $\beta \in (0.5, 1]$, all $t \leqslant t_0$, and some choice of $0 < t_0 < \infty$. An appropriate $t_0$ depends on $L > 0$, however, and $L$ is unknown in practice.

The lack of convergence of the estimates in the $L = 15$ case appears to be associated with the rung process being trapped in a region of $\Lambda$.

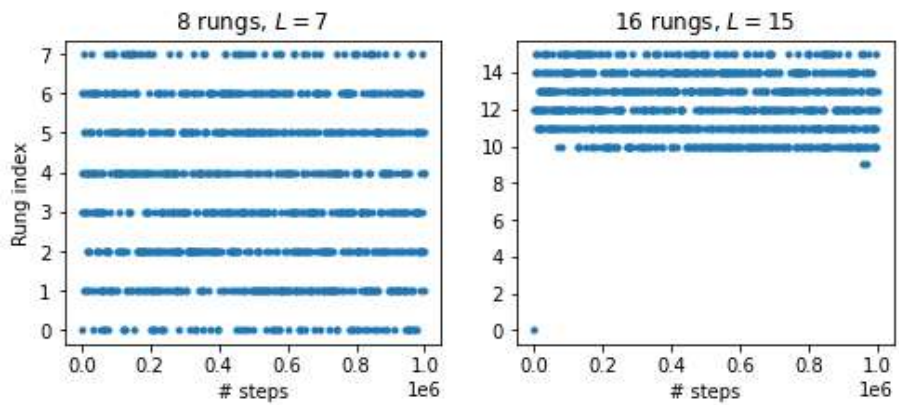

Figure 8.1.2: Subsampled rung trajectory for $L = 7, 15$.

## 8.2 History forgetting

We now investigate the use of history forgetting to alleviate the issue of dependence on initial conditions and the early history of the simulation. The use of history forgetting as a method of providing burn-in is mainly out of consideration for the complexity for the sampling in $\mathcal{X}$, but in the current example the samples $\mathcal{X}^t$ are being drawn conditionally independently from the distributions $\rho_{\mathcal{K}^{t-1}}(x)dx$, so there is no correlation to $\mathcal{X}^{t-1}$ as would normally be the case in practice (when a local move-based sampler such as MCMC or MD is used). However, there is still dependence on the early history of the process through $\{\mathcal{K}^t, t \geqslant 0\}$. Though the addition of history-forgetting will not remedy the convergence issues, it will at least alleviate the dependence on $(\mathcal{X}^t, \mathcal{K}^t)$ for $t$ early on in the simulation.

The setup below deviates from the setup in 8.1 only in the addition of the history-forgetting mechanism, using $\alpha = 0.19$ and 32 epochs.

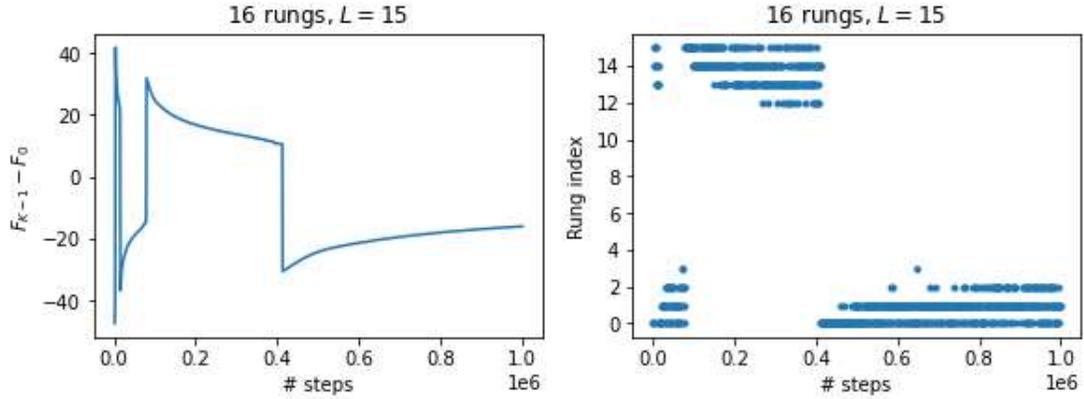

Figure 8.2.1: End-to-end free energy difference estimate and subsampled rung trajectory for $L = 15$.

The right-hand plot shows that the process is no longer tied to a fixed region of $\Lambda$, but now oscillates between different regions, in a way that is strongly correlated to the free energy estimates. The history forgetting mechanism does not by itself alleviate the exponential slowdown from the use of $O(t^{-1})$ gain.

## 8.3 Windows

To avoid the exponential slowdown when $L = 15$, we now use windows to decompose $\Lambda = [0, L]$ into windows $W_j \subset \Lambda$ such that $\{W_j, j \in [J]\}$ forms a double cover of $[0, L]$, as introduced in Sections 2.1 and 2.2 of the main text. Specifically, we use two window setups for $L = 15$, $K = 16$:

- (3 windows) : $W_1 = \{0, \ldots, 15\}$, $W_2 = \{0, \ldots, 7\}$, $W_3 = \{8, \ldots, 15\}$;

- (5 windows) : $W_1 = \{0, \ldots, 7\}$, $W_2 = \{8, \ldots, 15\}$, $W_3 = \{0, \ldots, 3\}$, $W_4 = \{4, \ldots, 11\}$, $W_5 = \{12, \ldots, 15\}$.

From Section 8.1, we know that even the featureless TSS estimator can handle $L = 7$, and thus we expect the second setup to accurately estimate the free energies within a relatively small number of steps. The plots below, which use these window setups (and also have the history forgetting parameter $\alpha = 0.19$ with 32 epochs), show this.

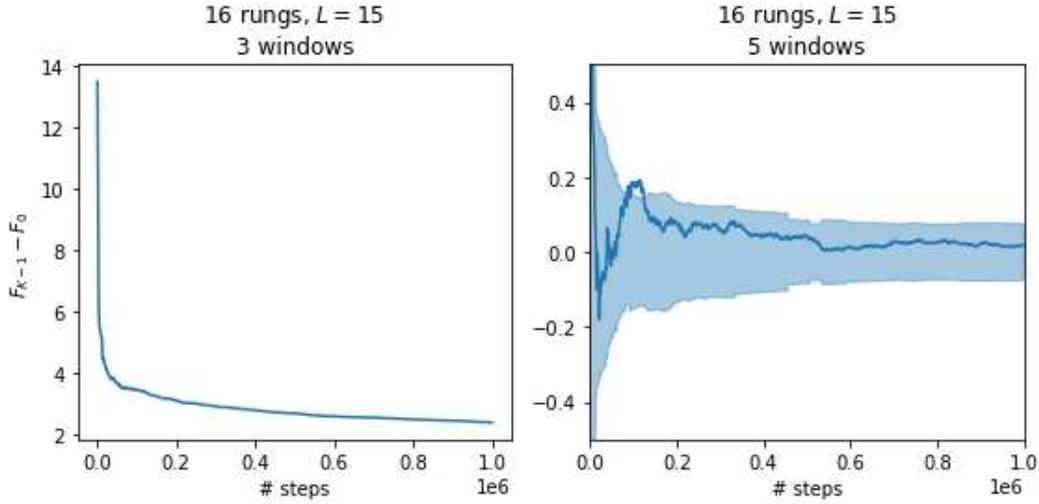

Figure 8.3.1: End-to-end free energy difference estimates for $L = 15$ using 3 and 5 windows.

We note that the rung dynamics properly explore $\Lambda$ even in the case of three windows, but that in the three-window case the free energy estimates tend to lag behind the sampling and suffer from the exponential slowdown. Thus proper exploration of $\Lambda$ is a prerequisite for proper estimation, but not a guarantee.

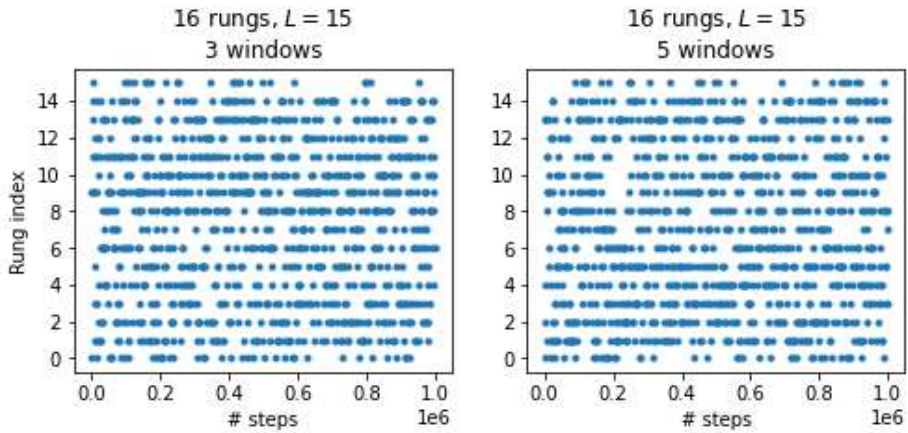

Figure 8.3.2: Subsampled rung trajectory for $L = 15$ using 3 and 5 windows.

## 8.4 Visit control

Windows handle the problem of convergence by decomposing the space $\Lambda$, but this approach is not entirely satisfactory for at least two reasons. First, determining an appropriate size appears to be a trial and error process. Second, restricting the window size can lower the statistical efficiency by restricting the rung movement. In order to make the free energy estimates less sensitive to the selection of window size, we introduce the visit control mechanism, as described in Section 1.3 of the main text.

We plot here $\eta = 4$ for $L = 15$ and $L = 63$ (using $\alpha = 0.19$, 32 epochs and a single window),

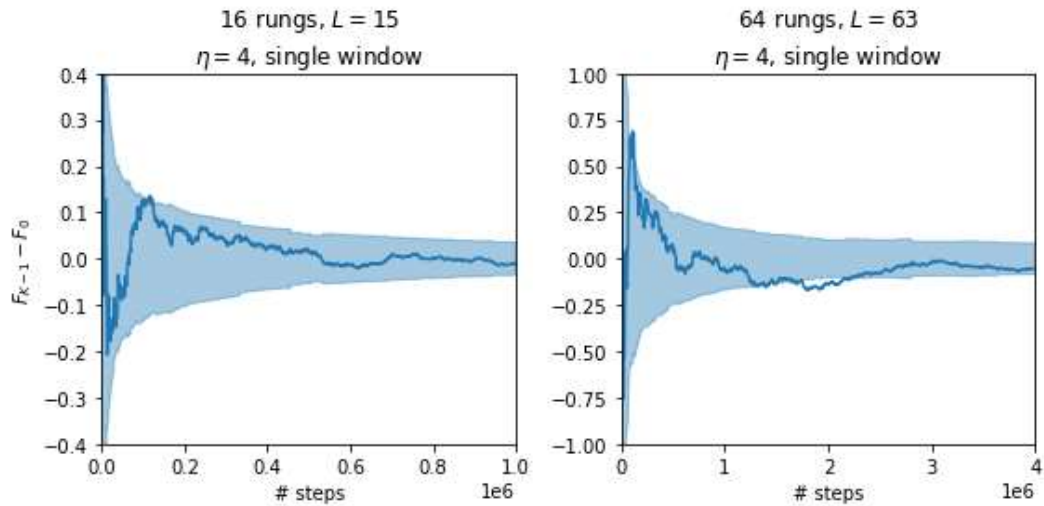

Figure 8.4.1: End-to-end free energy difference estimates for $L = 15, 63$ with visit control turned on ($\eta = 4$).

The visit control mechanism does a remarkable job of accelerating convergence of the free energy difference estimates to their equilibrium for both $L = 15$ and $L = 63$. In additional empirical studies (data not shown) using larger values of $L$ (up to $2,000$) the visit control mechanism remained effective at bringing estimates close to a neighborhood of the equilibrium value. (We have not studied the limit $L \to \infty$ since we couple the visit control mechanism with the windowing system in practical examples, which keeps $L$ bounded.)

## 8.5 Self-adjustment

An approach to improve the sampling within a window is to apply the transition kernel more frequently between estimator updates. To simplify the illustration of such an approach, we consider the case of only a single window, but the approach extends naturally to multiple windows. Within a window, the transition kernel draws an independent sample $X^{t+1}$ according to the density $\rho_{\mathcal{K}^t}$, followed by an update to $\mathcal{K}^{t+1}$ given $X^{t+1}$; this composed kernel was denoted $P_{\pi,Z} = P^k_{\pi,Z} \circ P^x_{\pi,Z}$ in Section 1.1 of the main text. To improve the sampling, we use instead the "high-frequency" kernel $P^{\nu}_{\pi,Z} = P_{\pi,Z} \circ \cdots \circ P_{\pi,Z}$ ($\nu$ times).

We demonstrate (using $\alpha = 0.19$, 32 epochs, $\eta = 4$ and a single window) the results of doing so for $L = 15$ using $L = 1, 32$, and for $L = 63$ using $\nu = 1, 32, 100$.

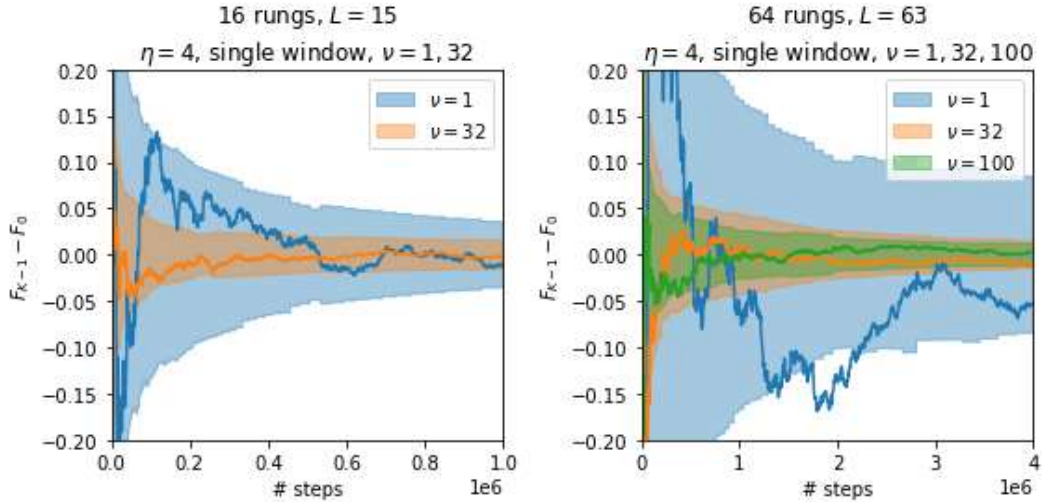

Figure 8.5.1

The variance reduction due to self-adjustment (large ν) is significant, the variance for $\nu = 32, L = 63$ being roughly 25 times smaller than that of $\nu = 1, L = 63$, and $\nu = 100, L = 63$ roughly 50 times smaller than that of $\nu = 1, L = 63$. Moreover, from the expressions for the asymptotic variance in Section 3, the asymptotic variance of $F_{K-1} - F_0$ is $O(L^2)$ for both MBAR and TSS for this model when $\nu = 1$. In the limit $\nu \to \infty$, however, the asymptotic variance of $F_{K-1} - F_0$ becomes $O(L)$ for TSS. In other words, the use of a very high-frequency kernel has the potential to change the scaling of the problem, suggesting that there may be regimes in practice where offline estimators such as MBAR cannot outperform on-the-fly estimators like TSS.

## 8.6 Summary

This section started with the SAMS algorithm, and examined the additional features implemented in TSS that make it a robust algorithm for use in computational chemistry. The history-forgetting automatically provides a method of performing burn-in, which means the system does not necessarily need to be well-equilibrated beforehand. Next, the windowing system introduces the notion of locality, and makes the algorithm highly scalable, while visit control decreases the sensitivity of the free energy estimates on the window size and accelerates convergence of the estimates by modifying the gain and biasing the rung dynamics. Finally, by using multiple rung moves between estimator updates, we can reduce the asymptotic variance of the free energy estimates without performing more updates to the free energy estimators.